\documentclass[10pt]{article}
\usepackage{amsfonts,graphicx,epsfig,color,fancyhdr}
\usepackage[letterpaper,asymmetric,left=1.5in,right=1.5in,top=1.5in,bottom=1.5in,
bindingoffset=0.0in]{geometry}

\title{Andreev's Theorem on hyperbolic polyhedra}

\author{Roland K. W. Roeder \footnote{rroeder@fields.utoronto.ca}
\footnote{Supported by a U.S. National Defense Science and Engineering
Fellowship and later by a NSF Integrative Graduate Research and Training (IGERT)
Fellowship.}, John H. Hubbard, and William D. Dunbar}

\input{psfig.sty}

\newtheorem{thm}{Theorem}[section]
\newtheorem{prop}[thm]{Proposition}
\newtheorem{cor}[thm]{Corollary}
\newtheorem{lem}[thm]{Lemma}
\newtheorem{slem}[thm]{Sub-lemma}

\newcommand{\Endproof}{$\Box$ \vspace{.1in}}
\newcommand{\Area}{{\rm Area}}

\begin{document}

\maketitle

\begin{center} \bf Abstract \end{center}
{\small
In 1970, E. M. Andreev published a classification of all three-dimensional
compact hyperbolic polyhedra (other than tetrahedra) having non-obtuse dihedral
angles \cite{AND,AND2}.  Given a combinatorial description of a polyhedron,
$C$, Andreev's Theorem provides five classes of linear inequalities, depending
on $C$, for the dihedral angles, which are necessary and sufficient conditions
for the existence of a hyperbolic polyhedron realizing $C$ with the assigned
dihedral angles.  Andreev's Theorem also shows that the resulting polyhedron is
unique, up to hyperbolic isometry. 

Andreev's Theorem is both an interesting statement about the geometry of
hyperbolic 3-dimensional space, as well as a fundamental tool used in the
proof for Thurston's Hyperbolization Theorem for 3-dimensional
Haken manifolds.  It is also remarkable to what level the proof of Andreev's
Theorem resembles (in a simpler way) the proof of Thurston.

We correct a fundamental error in Andreev's proof of existence and also provide a
readable new proof of the other parts of the proof of Andreev's Theorem, because
Andreev's paper has the reputation of being ``unreadable''.  
}

\begin{center} \bf R\'esum\'e \end{center}
{\small
E. M. Andreev a publi\'e en 1970 une classification des poly\`edres
hyperboliques compacts de dimension trois (autrement que les t\'etra\`edres) dont les angles di\`edres sont
non-obtus \cite{AND,AND2}.  Etant donn\'e une description combinatoire d'un
poly\`edre $C$, le Th\'eor\`eme d'Andreev dit que les angles di\`edres
possibles sont exactement d\'ecrits par cinq classes d'in\'egalit\'es
lin\'eaires. Le Th\'eor\`eme d'Andreev d\'emontre \'egalement que le
poly\`edre r\'esultant est alors unique \`a isom\'etrie hyperbolique pr\`es. 

D'une part, le Th\'eor\`eme de Andreev est \'evidemment un \'enonc\'e
int\'eressant de la g\'eom\'etrie de l'espace hyperbolique en dimension 3; 
d'autre part c'est un outil essentiel dans la preuve du Th\'eor\`eme
d'Hyperbolization de Thurston pour les vari\'et\'es Haken de dimension 3.  Il
est d'ailleurs remarquable \`a quel point la d\'emonstration d'Andreev rappelle
(en plus simple) la d\'emonstration de Thurston. 

La d\'emonstration d'Andreev contient une erreur importante. Nous corrigeons
ici cette erreur et nous fournissons aussi une nouvelle preuve lisible des
autres parties de la preuve, car le papier d'Andreev a la r\'eputation
d'\^etre ``illisible''.
}

\section{Statement of Andreev's Theorem}

Andreev's Theorem provides a complete characterization of compact
hyperbolic polyhedra having non-obtuse dihedral angles.  This
classification is essential for proving Thurston's Hyperbolization theorem
for Haken 3-manifolds and is also a particularly beautiful and interesting
result in its own right.  Complete and detailed proofs of Thurston's
Hyperbolization for Haken 3-manifolds are available written in English by
Jean-Pierre Otal \cite{OT} and in French by Michel Boileau \cite{BOI}.

In this paper, we prove Andreev's Theorem based on the main ideas from his
original proof \cite{AND}.  However, there is an error in Andreev's proof of
existence.  We explain this error in Section 6 and provide a correction.
Although the other parts of the proof are proven in much the same way as
Andreev proved them, we have re-proven them and re-written them to verify them
as well as to make the overall proof of Andreev's Theorem clearer.  This paper
is based on the doctoral thesis of the first author \cite{ROE}, although
certain proofs have been streamlined, especially in Sections 4 and 5.

The reader may also wish to consider the similar results of Rivin and Hodgson
\cite{RH,H}, Thurston \cite[Chapter 13]{T_NOTES}, Marden and Rodin \cite{MR},
Bowers and Stephenson \cite{STEVE}, Rivin
\cite{RIV_IDEAL2,RIV_IDEAL1,RIV_IDEAL3}, and Bao and Bonahon \cite{BAO}.  In
\cite{RH}, the authors prove a more general statement than Andreev's Theorem
and in \cite{H} Hodgson deduces Andreev's Theorem as a consequence of their
previous work.  The proof in \cite{RH} is similar to the one presented here,
except that the conditions classifying the polyhedra are written in terms of
measurements in the De Sitter space, the space dual to the hyperboloid model of
hyperbolic space.  Although a beautiful result, the main drawback of this proof
is that the last sections of the paper, which are necessary for their proof
that such polyhedra exist, are particularly hard to follow. 

Marden and Rodin \cite{MR} and Thurston \cite[Chapter 13]{T_NOTES} consider
configurations of circles with assigned overlap angles on the Riemann Sphere
and on surfaces of genus $g$ with $g > 0$.  Such a configuration on the Riemann
Sphere corresponds to a configuration of hyperbolic planes in the
Poincar\'e ball model of hyperbolic space.  Thus, there is a direct connection
between circle patterns and hyperbolic polyhedra.  The proof of Thurston
provides a classification of configurations of circles on surfaces of
genus $g > 0$.  The proof of Marden and Rodin \cite{MR} is an adaptation of
Thurston's circle packing theorem to the Riemann Sphere, resulting in a
theorem similar to Andreev's Theorem.  Although Thurston's statement has
analogous conditions to Andreev's classical conditions, Marden and Rodin
require that the sum of angles be less than $\pi$ for every triple of circles
for which each pair intersects.  This prevents the patterns of overlapping
circles considered in their theorem from corresponding to compact hyperbolic
polyhedra.

Bowers and Stephenson \cite{STEVE} prove a ``branched version'' of Andreev's
Theorem, also in terms of circle patterns on the Riemann Sphere.   Instead of
the continuity method used by Thurston and Marden-Rodin, Bowers and Stephenson
use ideas intrinsic to the famous Uniformization Theorem from complex analysis.
The unbranched version of their theorem provides a complete proof of Andreev's
Theorem, which provides an alternative to the proof presented here. 

Rivin has proven beautiful results on ideal hyperbolic polyhedra
having arbitrary dihedral angles \cite{RIV_IDEAL2,RIV_IDEAL1} (see also
Gu{\'e}ritaud \cite{GUE} for an alternative viewpoint, with exceptionally clear
exposition.)  Similar nice results are proven for hyperideal polyhedra by Bao
and Bonahon \cite{BAO}.  Finally, the papers of Vinberg on discrete groups of
reflections in hyperbolic space \cite{AVS,VIN,VINREFL,VINVOL,VS} are also
closely related, as well as the work of Bennett and Luo \cite{LUO} and
Schlenker \cite{SCH2,SCH1,SCH3}.

\vspace{.1in}

Let $E^{3,1}$ be $\mathbb{R}^4$ with the indefinite metric $\Vert {\bf x}
\Vert^2 = -x_0^2+x_1^2+x_2^2+x_3^2$.  The space of ${\bf x}$ for which this
indefinite metric vanishes is typically referred to as the lightcone, which we
denote by $C$.

In this paper, we work in the hyperbolic space $\mathbb{H}^3$ given by the
component of the subset of $E^{3,1}$ given by
$$\Vert {\bf x} \Vert^2 = -x_0^2+x_1^2+x_2^2+x_3^2 = -1$$

\noindent
having $x_0 > 0$, with the Riemannian metric induced by the indefinite
metric
$$-dx_0^2+dx_1^2+dx_2^2+dx_3^2.$$

Hyperbolic space $\mathbb{H}^3$ can be compactified by adding the set of rays
in $\{{\bf x} \in C \mbox{  :  } x_0 \ge 0 \},$
which clearly form a topological space $\partial \mathbb{H}^3$
homeomorphic to the sphere $\mathbb{S}^{2}$.  We will refer to points in
$\partial \mathbb{H}^3$ as {\it points at infinity} and refer to the
compactification as $\overline{\mathbb{H}^3}$.  For more details, see \cite[p. 66]{THURSTON_BOOK}.

The hyperplane orthogonal to a vector ${\bf v} \in
E^{3,1}$ intersects $\mathbb{H}^3$ if and only if $\langle{\bf v},{\bf
v}\rangle> 0$.  Let ${\bf v} \in E^{3,1}$ be a vector with $\langle{\bf
v},{\bf v}\rangle > 0$, and define
\begin{eqnarray*}
P_{\bf v} = \{{\bf w} \in \mathbb{H}^3 | \langle{\bf w},{\bf v}\rangle
= 0\} \mbox{      and      }
H_{\bf v} = \{{\bf w} \in \mathbb{H}^3 | \langle{\bf w},{\bf v}\rangle \leq
0 \}
\end{eqnarray*}

\noindent to be the hyperbolic plane orthogonal to ${\bf v}$ and the
corresponding closed half space, oriented so that ${\bf v}$ is the outward pointing normal.

If one normalizes $\langle{\bf v},{\bf v}\rangle = 1$ and $\langle{\bf w},{\bf
w}\rangle = 1$ the planes $P_{\bf v}$ and $P_{\bf w}$ in $\mathbb{H}^3$
intersect in a line if and only if $\langle{\bf v},{\bf w}\rangle^2 <
1$, in which case their dihedral angle is $\arccos(-\langle{\bf v},{\bf
w}\rangle)$.  They intersect in a single point at infinity if and only if
$\langle{\bf v},{\bf w}\rangle^2 = 1$; in this case their dihedral angle is
$0$.

A {\it hyperbolic polyhedron} is an intersection

$$P = \bigcap_{i=0}^N H_{\bf v_i} $$

\noindent
having non-empty interior.  Throughout this paper we will make the assumption
that ${\bf v}_1,\cdots,{\bf v}_N$ form a minimal set of vectors specifying $P$.  That
is, we assume that none of the half-spaces $H_{\bf v_i}$ contains the
intersection of all the others.

It is not hard to verify that if $H_{\bf v_i}, H_{\bf v_j},H_{\bf v_k}$
are three distinct halfspaces appearing in the definition of the polyhedron
$P$, then the three vectors ${\bf v_i}$, ${\bf v_j}$, and ${\bf v_k}$ will always be
linearly independent.  For example, if ${\bf v_i} = {\bf v_j} + {\bf v_k}$, then $H_{\bf v_i}$
would be a subset of $H_{\bf  v_j} \cap H_{\bf v_k}$, contradicting minimality.

\vspace{0.05in}

We will often use the Poincar\'e ball model of hyperbolic space, given by the
open unit ball in $\mathbb{R}^3$ with the metric

$$4\frac{dx_1^2+dx_2^2+dx_3^2}{(1 -\Vert {\bf x}\Vert^2)^2}$$

\noindent
and the upper half-space model of hyperbolic space,
given by the subset of $\mathbb{R}^3$ with $x_3 > 0$ equipped with the
metric
$$\frac{dx_1^2+dx_2^2+dx_3^2}{x_3^2}.$$

\noindent
Both of these models are isomorphic to $\mathbb{H}^3$.

Hyperbolic planes in these models correspond to Euclidean hemispheres and
Euclidean planes that intersect the boundary perpendicularly.  Furthermore,
these models are conformally correct, that is, the hyperbolic angle between a
pair of such intersecting hyperbolic planes is exactly the Euclidean angle
between the corresponding spheres or planes.

Below is an image of a hyperbolic polyhedron depicted in the Poincar\'e 
ball model with the sphere at infinity shown for reference.
It was displayed in the excellent computer program
Geomview \cite{GEO}.

\begin{center}
\includegraphics[scale = .45]{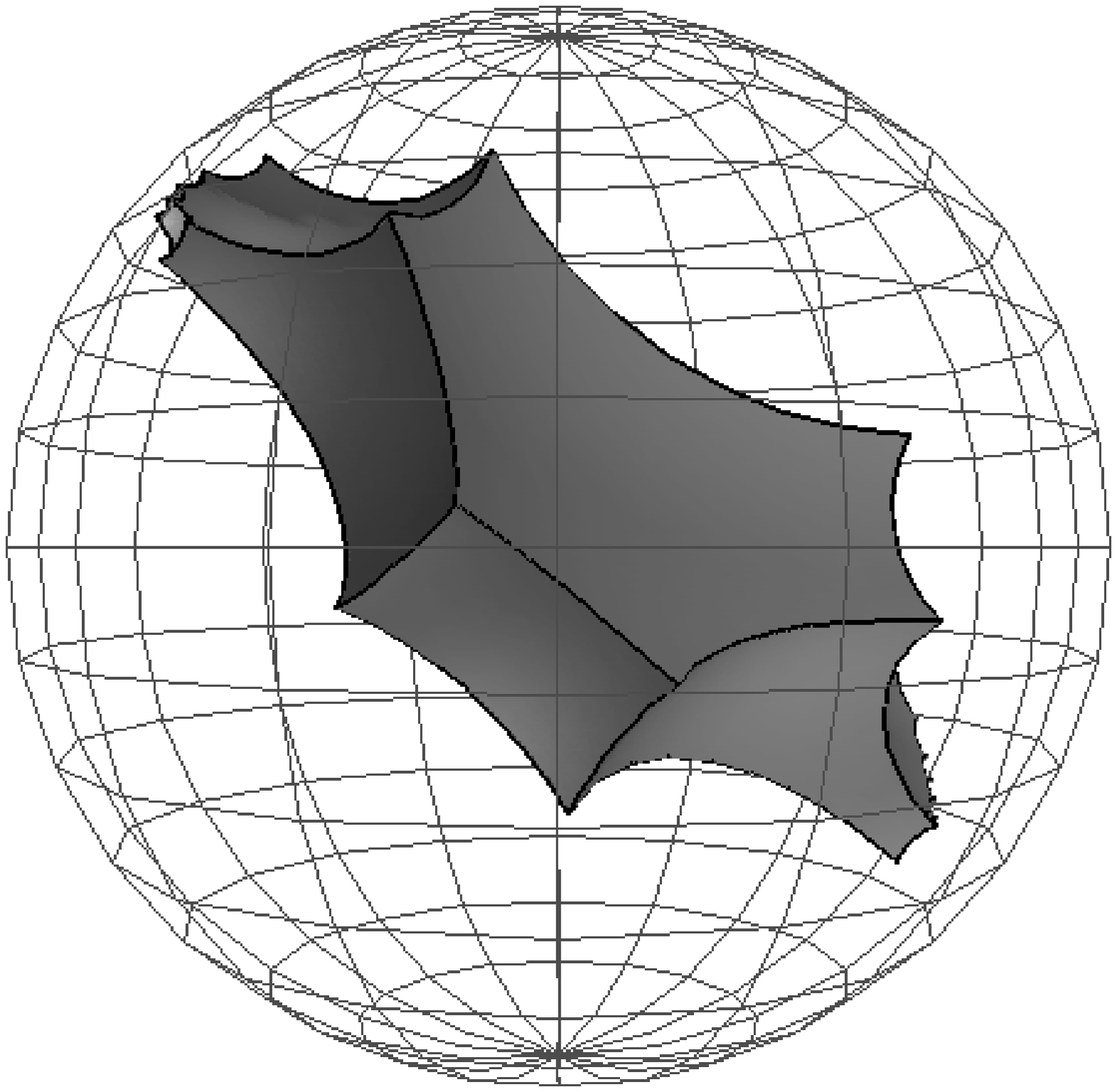}
\end{center}

\vspace{.2in}
\noindent
{\large \it Abstract polyhedra 
and Andreev's Theorem}
\vspace{.1in}

\noindent
Some elementary combinatorial facts about hyperbolic polyhedra are
essential before we can state Andreev's Theorem.  Notice that a compact
hyperbolic polyhedron $P$ is topologically a 3-dimensional ball, and its
boundary a 2-sphere $\mathbb{S}^2$.  The face structure of $P$ gives
$\mathbb{S}^2$ the structure of a cell complex $C$ whose faces correspond
to the faces of $P$, and so forth.

Considering only hyperbolic polyhedra with non-obtuse dihedral angles
simplifies the combinatorics of any such $C$: 

\begin{prop} \label{TRIVALENT}
(a)  A vertex of a non-obtuse hyperbolic polyhedron $P$ is the
intersection of exactly 3 faces. \newline (b) For such a $P$, we can
compute the angles of the faces in terms of the dihedral angles; these
angles are also $\leq \pi/2$. 

\end{prop}

\vspace{.05in}
\noindent
{\bf Proof:} 
Let $v$ be a finite vertex where $n$ faces of $P$ meet.  After an
appropriate isometry, we can assume that $v$ is the origin in the Poincar\'e
ball model, so that the faces at $v$ are subsets of Euclidean planes
through the origin.  A small sphere centered at the origin
will intersect $P$ in a spherical $n$-gon $Q$ whose angles are the
dihedral angles between faces.  Call these angles $\alpha_1,...,\alpha_n$. 
Re-scale $Q$ so that it lies on the sphere of unit radius, then
the Gauss-Bonnet formula gives $\alpha_1+\cdots+\alpha_n = \pi(n-2) +
\Area(Q)$. The restriction to $\alpha_i \leq \pi/2$ for all $i$ gives $n\pi/2
\geq \pi(n-2) + \Area(Q)$.  Hence, $n\pi/2 < 2\pi$.  We conclude that $n =
3$. 

The
edge lengths of $Q$ are precisely the angles in the faces at the origin. 
Supposing that $Q$ has angles $(\alpha_i,\alpha_j,\alpha_k)$ and edge
lengths $(\beta_i,\beta_j,\beta_k)$ with the edge $\beta_l$ opposite of
angle $\alpha_l$ for each $l$, The law of cosines in spherical geometry
states that:  

\begin{eqnarray} \label{TLC}
\cos(\beta_i) = \frac{\cos(\alpha_i)
+\cos(\alpha_j)\cos(\alpha_k)} {\sin(\alpha_j)\sin(\alpha_k)}.
\end{eqnarray}

\noindent
Hence, the face angles are calculable from the dihedral angles. 
They are non-obtuse, since  the right-hand side of the equation is positive
for $\alpha_i,\alpha_j,\alpha_k$ non-obtuse. 
(Equation (\ref{TLC}) will be used frequently
throughout this paper.)
\Endproof

The fundamental axioms of incidence place the following, obvious, further
restrictions on the complex $C$: 

\begin{itemize}

\item Every edge of $C$ belongs to exactly two faces.

\item A non-empty intersection of two faces is either an edge or a vertex.

\item Every face contains not fewer than three edges.

\end{itemize}

We will call any trivalent cell complex $C$ on $\mathbb{S}^2$ that satisfies
the three conditions above an {\it abstract polyhedron}.  Notice that since $C$
must be a trivalent cell complex on $\mathbb{S}^2$, its dual, $C^*$, has only
triangular faces.  The three other conditions above ensure that the dual
complex $C^*$ is a simplicial complex, which we embed in the same
$\mathbb{S}^2$ so that the vertex corresponding to any face of $C$ is an
element of the face, etc.  (Andreev refers to this dual complex as the {\it
scheme of the polyhedron}.)  The figure below shows an abstract  polyhedron $C$
drawn in the plane (i.e.  with one of the faces corresponding to the region
outside of the figure.) The dual complex is also shown, in dashed lines.

\vspace{.07in}
\begin{center}
\begin{picture}(0,0)%
\epsfig{file=comb_poly.pstex}%
\end{picture}%
\setlength{\unitlength}{3947sp}%
\begingroup\makeatletter\ifx\SetFigFont\undefined%
\gdef\SetFigFont#1#2#3#4#5{%
  \reset@font\fontsize{#1}{#2pt}%
  \fontfamily{#3}\fontseries{#4}\fontshape{#5}%
  \selectfont}%
\fi\endgroup%
\begin{picture}(4015,2919)(139,-2143)
\end{picture}%

\end{center}
\vspace{.07in}

 We call a simple closed curve $\Gamma$ formed of $k$ edges of $C^*$ a
{\it k-circuit} and if all of the endpoints of the edges of $C$
intersected by $\Gamma$ are distinct, we call such a circuit a {\it
prismatic k-circuit}.  The figure below shows the same abstract polyhedron
as above, except this time the prismatic 3-circuits are dashed, the prismatic
4-circuits are dotted, and the dual complex is not shown.

\vspace{.07in}
\begin{center}
\begin{picture}(0,0)%
\epsfig{file=comb_poly2.pstex}%
\end{picture}%
\setlength{\unitlength}{3947sp}%
\begingroup\makeatletter\ifx\SetFigFont\undefined%
\gdef\SetFigFont#1#2#3#4#5{%
  \reset@font\fontsize{#1}{#2pt}%
  \fontfamily{#3}\fontseries{#4}\fontshape{#5}%
  \selectfont}%
\fi\endgroup%
\begin{picture}(4152,2571)(80,-1853)
\end{picture}%
\end{center}
\vspace{.07in}

Before stating Andreev's Theorem, we prove two basic lemmas about abstract
polyhedra:

\begin{lem}\label{3CIRC}
If $\gamma$ is a 3-circuit that is not prismatic in an abstract polyhedron 
$C$ intersecting edges $e_1,e_2$, and $e_3$,
then edges $e_1,e_2$, and $e_3$ meet at a vertex.
\end{lem}

\noindent {\bf Proof:}
Since $\gamma$ is 3-circuit that is not prismatic, a pair of the edges meet at a
vertex.  We suppose that $e_1$ and $e_2$ meet at this vertex, which we label
$v_1$.  Since the vertices of $C$ are trivalent, there is some edge $e'$
meeting $e_1$ and $e_2$ at $v_1$.  We suppose that $e'$ is not the edge
$e_3$ to obtain a contradiction.  Moving $\gamma$ past the vertex $v_1$, we can
obtain a new circuit $\gamma'$ intersecting only the two edges $e_3$ and $e'$.

\vspace{.07in}
\begin{center}
\begin{picture}(0,0)%
\includegraphics{three_circ.pstex}%
\end{picture}%
\setlength{\unitlength}{3947sp}%
\begingroup\makeatletter\ifx\SetFigFont\undefined%
\gdef\SetFigFont#1#2#3#4#5{%
  \reset@font\fontsize{#1}{#2pt}%
  \fontfamily{#3}\fontseries{#4}\fontshape{#5}%
  \selectfont}%
\fi\endgroup%
\begin{picture}(4957,962)(649,-722)
\put(826,112){\makebox(0,0)[lb]{\smash{\SetFigFont{8}{9.6}{\familydefault}{\mddefault}{\updefault}{\color[rgb]{0,0,0}\small{$e_2$}}%
}}}
\put(856,-691){\makebox(0,0)[lb]{\smash{\SetFigFont{8}{9.6}{\familydefault}{\mddefault}{\updefault}{\color[rgb]{0,0,0}\small{$e_1$}}%
}}}
\put(895,-286){\makebox(0,0)[lb]{\smash{\SetFigFont{8}{9.6}{\familydefault}{\mddefault}{\updefault}{\color[rgb]{0,0,0}\small{$v_1$}}%
}}}
\put(1334,-198){\makebox(0,0)[lb]{\smash{\SetFigFont{8}{9.6}{\familydefault}{\mddefault}{\updefault}{\color[rgb]{0,0,0}\small{$e'$}}%
}}}
\put(2057, 88){\makebox(0,0)[lb]{\smash{\SetFigFont{8}{9.6}{\familydefault}{\mddefault}{\updefault}{\color[rgb]{0,0,0}\small{$\gamma$}}%
}}}
\put(2313,  7){\makebox(0,0)[lb]{\smash{\SetFigFont{8}{9.6}{\familydefault}{\mddefault}{\updefault}{\color[rgb]{0,0,0}\small{Properties of abstract}}%
}}}
\put(2547,-299){\makebox(0,0)[lb]{\smash{\SetFigFont{8}{9.6}{\familydefault}{\mddefault}{\updefault}{\color[rgb]{0,0,0}\small{polyhedra imply}}%
}}}
\put(4133,138){\makebox(0,0)[lb]{\smash{\SetFigFont{8}{9.6}{\familydefault}{\mddefault}{\updefault}{\color[rgb]{0,0,0}\small{$e_2$}}%
}}}
\put(4163,-665){\makebox(0,0)[lb]{\smash{\SetFigFont{8}{9.6}{\familydefault}{\mddefault}{\updefault}{\color[rgb]{0,0,0}\small{$e_1$}}%
}}}
\put(4202,-260){\makebox(0,0)[lb]{\smash{\SetFigFont{8}{9.6}{\familydefault}{\mddefault}{\updefault}{\color[rgb]{0,0,0}\small{$v_1$}}%
}}}
\put(5377,134){\makebox(0,0)[lb]{\smash{\SetFigFont{8}{9.6}{\familydefault}{\mddefault}{\updefault}{\color[rgb]{0,0,0}\small{$\gamma$}}%
}}}
\put(1865,-190){\makebox(0,0)[lb]{\smash{\SetFigFont{8}{9.6}{\familydefault}{\mddefault}{\updefault}{\color[rgb]{0,0,0}\small{$e_3$}}%
}}}
\put(4872,-173){\makebox(0,0)[lb]{\smash{\SetFigFont{8}{9.6}{\familydefault}{\mddefault}{\updefault}{\color[rgb]{0,0,0}\small{$e_3$}}%
}}}
\put(5035,-422){\makebox(0,0)[lb]{\smash{\SetFigFont{8}{9.6}{\familydefault}{\mddefault}{\updefault}{\color[rgb]{0,0,0}\small{$\gamma'$}}%
}}}
\put(1735,-463){\makebox(0,0)[lb]{\smash{\SetFigFont{8}{9.6}{\familydefault}{\mddefault}{\updefault}{\color[rgb]{0,0,0}\small{$\gamma'$}}%
}}}
\end{picture}
\end{center}
\vspace{.07in}

\noindent

The curve $\gamma'$ intersects only two edges, hence it only crosses
two faces of $C$.  However, this implies that these two faces of $C$
intersect along the two distinct edges $e'$ and $e_3$, contrary to fact
that two faces of an abstract polyhedron which intersect do so along a
single edge. \Endproof

\begin{lem}\label{4CIRC}
Let $C$ be an abstract polyhedron having no prismatic 3-circuits.  If
$\gamma$ is a 4-circuit which is not prismatic, then $\gamma$ separates
exactly two vertices of $C$ from the remaining vertices of $C$.
\end{lem}

\noindent {\bf Proof:}
Suppose that $\gamma$ crosses edges $e_1,e_2,e_3$, and $e_4$ of $C$.
Because $\gamma$ is not a prismatic 4-circuit, a pair of these edges
meet at a vertex.
Without loss of generality, we suppose that edges
$e_1$ and $e_2$ meet at this vertex, which we denote $v_1$.
Since $C$
is trivalent, there is some edge $e'$ meeting $e_1$ and $e_2$ at $v_1$.
Let $\gamma '$ be the 3-circuit intersecting edges $e_3, e_4$ and
$e'$, obtained by sliding $\gamma$ past the vertex $v_1$.  Since $C$
has no prismatic 3-circuits, $\gamma'$ is not prismatic, so by Lemma
\ref{3CIRC}, edges $e_3,e_4$, and $e'$ meet at another vertex $v_2$.
The entire configuration is shown in the diagram below.

\vspace{.07in}
\begin{center}
\begin{picture}(0,0)%
\includegraphics{four_circ.pstex}%
\end{picture}%
\setlength{\unitlength}{3947sp}%
\begingroup\makeatletter\ifx\SetFigFont\undefined%
\gdef\SetFigFont#1#2#3#4#5{%
  \reset@font\fontsize{#1}{#2pt}%
  \fontfamily{#3}\fontseries{#4}\fontshape{#5}%
  \selectfont}%
\fi\endgroup%
\begin{picture}(4647,1139)(649,-722)
\put(4434,-143){\makebox(0,0)[lb]{\smash{{\SetFigFont{8}{9.6}{\familydefault}{\mddefault}{\updefault}{\color[rgb]{0,0,0}\small{$e'$}}%
}}}}
\put(3946,165){\makebox(0,0)[lb]{\smash{{\SetFigFont{8}{9.6}{\familydefault}{\mddefault}{\updefault}{\color[rgb]{0,0,0}\small{$e_2$}}%
}}}}
\put(3976,-638){\makebox(0,0)[lb]{\smash{{\SetFigFont{8}{9.6}{\familydefault}{\mddefault}{\updefault}{\color[rgb]{0,0,0}\small{$e_1$}}%
}}}}
\put(5116,-562){\makebox(0,0)[lb]{\smash{{\SetFigFont{8}{9.6}{\familydefault}{\mddefault}{\updefault}{\color[rgb]{0,0,0}\small{$e_4$}}%
}}}}
\put(5161,248){\makebox(0,0)[lb]{\smash{{\SetFigFont{8}{9.6}{\familydefault}{\mddefault}{\updefault}{\color[rgb]{0,0,0}\small{$e_3$}}%
}}}}
\put(5296,-202){\makebox(0,0)[lb]{\smash{{\SetFigFont{8}{9.6}{\familydefault}{\mddefault}{\updefault}{\color[rgb]{0,0,0}\small{$\gamma$}}%
}}}}
\put(5034,-374){\makebox(0,0)[lb]{\smash{{\SetFigFont{8}{9.6}{\familydefault}{\mddefault}{\updefault}{\color[rgb]{0,0,0}\small{$\gamma'$}}%
}}}}
\put(4830,-233){\makebox(0,0)[lb]{\smash{{\SetFigFont{8}{9.6}{\familydefault}{\mddefault}{\updefault}{\color[rgb]{0,0,0}\small{$v_2$}}%
}}}}
\put(2581,-306){\makebox(0,0)[lb]{\smash{{\SetFigFont{8}{9.6}{\familydefault}{\mddefault}{\updefault}{\color[rgb]{0,0,0}\small{3-circuits implies}}%
}}}}
\put(2610, 21){\makebox(0,0)[lb]{\smash{{\SetFigFont{8}{9.6}{\familydefault}{\mddefault}{\updefault}{\color[rgb]{0,0,0}\small{No prismatic}}%
}}}}
\put(4023,-233){\makebox(0,0)[lb]{\smash{{\SetFigFont{8}{9.6}{\familydefault}{\mddefault}{\updefault}{\color[rgb]{0,0,0}\small{$v_1$}}%
}}}}
\put(826,112){\makebox(0,0)[lb]{\smash{{\SetFigFont{8}{9.6}{\familydefault}{\mddefault}{\updefault}{\color[rgb]{0,0,0}\small{$e_2$}}%
}}}}
\put(856,-691){\makebox(0,0)[lb]{\smash{{\SetFigFont{8}{9.6}{\familydefault}{\mddefault}{\updefault}{\color[rgb]{0,0,0}\small{$e_1$}}%
}}}}
\put(1996,-615){\makebox(0,0)[lb]{\smash{{\SetFigFont{8}{9.6}{\familydefault}{\mddefault}{\updefault}{\color[rgb]{0,0,0}\small{$e_4$}}%
}}}}
\put(2041,195){\makebox(0,0)[lb]{\smash{{\SetFigFont{8}{9.6}{\familydefault}{\mddefault}{\updefault}{\color[rgb]{0,0,0}\small{$e_3$}}%
}}}}
\put(2176,-255){\makebox(0,0)[lb]{\smash{{\SetFigFont{8}{9.6}{\familydefault}{\mddefault}{\updefault}{\color[rgb]{0,0,0}\small{$\gamma$}}%
}}}}
\put(1914,-427){\makebox(0,0)[lb]{\smash{{\SetFigFont{8}{9.6}{\familydefault}{\mddefault}{\updefault}{\color[rgb]{0,0,0}\small{$\gamma'$}}%
}}}}
\put(895,-286){\makebox(0,0)[lb]{\smash{{\SetFigFont{8}{9.6}{\familydefault}{\mddefault}{\updefault}{\color[rgb]{0,0,0}\small{$v_1$}}%
}}}}
\put(1246,-203){\makebox(0,0)[lb]{\smash{{\SetFigFont{8}{9.6}{\familydefault}{\mddefault}{\updefault}{\color[rgb]{0,0,0}\small{$e'$}}%
}}}}
\end{picture}%
\end{center}
\vspace{.07in}

\noindent
Therefore, the 4-circuit $\gamma$ separates the two
vertices $v_1$ and $v_2$ from the remaining vertices of $C$.
\Endproof

\begin{thm} {\bf Andreev's Theorem}

\noindent
Let $C$ be an abstract polyhedron with more than 4
faces and suppose that non-obtuse angles $\alpha_i$ are given corresponding
to each edge $e_i$ of $C$.  There is a 
compact hyperbolic polyhedron $P$ whose faces realize $C$
with dihedral angle $\alpha_i$ at each edge $e_i$ if and only if
the following five conditions all hold:

\setcounter{enumi}{0}
\begin{enumerate}

\item For each edge $e_i$, $\alpha_i > 0$.

\item Whenever 3 distinct edges $e_i,e_j,e_k$ meet at a vertex, 
$\alpha_i+\alpha_j+\alpha_k > \pi$.

\item Whenever $\Gamma$ is a prismatic 3-circuit intersecting edges
$e_i,e_j,e_k$, $\alpha_i+\alpha_j+\alpha_k < \pi$.

\item Whenever $\Gamma$ is a prismatic 4-circuit intersecting edges
$e_i,e_j,e_k,e_l$, then $\alpha_i+\alpha_j+\alpha_k+\alpha_l < 2\pi$.

\item Whenever there is a four sided face bounded by edges $e_1,$ $e_2,$
 $e_3,$ and $e_4$,
enumerated successively, with
edges $e_{12}, e_{23}, e_{34}, e_{41}$ entering the four vertices (edge
$e_{ij}$ connects to the ends of $e_i$ and $e_j$), then:
$$\alpha_1 + \alpha_3 + \alpha_{12} + \alpha_{23} + \alpha_{34} +
\alpha_{41} < 3\pi, \hspace{.2in} {\rm and}$$
$$\alpha_2 + \alpha_4 + \alpha_{12} + \alpha_{23} + \alpha_{34} +
\alpha_{41} < 3\pi.$$
\end{enumerate}

\noindent
Furthermore, this polyhedron is unique up to isometries of $\mathbb{H}^3$.
\end{thm}

In addition to the role that Andreev's theorem plays, as a bootstrap in the
proof of Thurston's hyperbolization theorem, it is worth noting that, in the
context of orbifolds, the former can be thought of as a very special case of
the latter (extended to Haken orbifolds as in \cite[Chapter 8]{BP} or
\cite{CHK}).  Consider closed $3$-orbifolds with underlying topological space a
$3$-ball, and with singular set equal to the boundary sphere.  That singular
set will consist of a trivalent graph, together with ``mirrors'' on the
complementary regions.  Each edge of the graph is labeled with an integer
$k>1$, corresponding to a dihedral angle of $\pi/k$.  The definition of a
$3$-orbifold implies that the angle sum at each vertex will satisfy condition
(2) in Andreev's theorem \cite[sections 6.1 and 6.3]{KAP}.  

Restrict the
combinatorics of the singular set, slightly more than in the statement of
Andreev's theorem: $C$ must be an abstract polyhedron with more than 5 faces.
Such an orbifold is Haken if and only if it is irreducible \cite[Proposition
13.5.2]{T_NOTES}.  Condition (3) in Andreev's theorem guarantees
irreducibility, and also, together with condition (4), guarantees that every
Euclidean 2-suborbifold is compressible.  Therefore, for Haken orbifolds of
this topological type, Andreev's theorem says precisely that having no
incompressible Euclidean 2-suborbifolds is equivalent to the existence of a
hyperbolic structure \cite[Section 13.6]{T_NOTES}; see also \cite[Section
6.4]{KAP}.

%Consider closed $3$-orbifolds with underlying 
%topological space a $3$-ball, and with singular set contained in the 
%boundary sphere.  That singular set will consist of a trivalent graph, 
%each edge of which is labeled with an integer greater than 1, together 
%with ``mirrors'' on the complementary regions.  

\vspace{.1in}

For a given $C$ let $E$ be the number of edges of $C$.  The subset
of $(0,\pi /2]^E$ satisfying these linear inequalities will be called the
{\it Andreev Polytope}, $A_C$.  Since $A_C$ is determined by linear
inequalities, it is convex.

Andreev's restriction to non-obtuse dihedral angles is emphatically necessary
to ensure that $A_C$ be convex.  Without this restriction, the corresponding
space of dihedral angles, $\Delta_C$, of compact (or finite volume) hyperbolic
polyhedra realizing a given $C$ is not convex \cite{DIAZ}.  In fact, the recent
work by D\'iaz \cite{DIAZ_ANDREEV} provides a detailed analysis of this space
of dihedral angles $\Delta_C$ for the class of abstract polyhedra $C$ obtained
from the tetrahedron by successively truncating vertices.  Her work nicely
illustrate the types of non-linear conditions that are necessary in a complete
analysis of the larger space of dihedral angles $\Delta_C$.

The work of Rivin \cite{RIV_IDEAL2, RIV_IDEAL1} shows that the space
of dihedral angles for ideal polyhedra forms a convex polytope, {\em without
the restriction to non-obtuse angles}.  (See also \cite{GUE}.)

Notice also that the hypothesis that the number of faces is greater than four
is also necessary because the space of non-obtuse dihedral angles for compact
tetrahedra is not convex \cite{ROE_TET}.  Conditions (1-5) remain necessary
conditions for compact tetrahedra, but they are no longer sufficient.

\begin{prop} \label{FIVE}
If $C$ is not the triangular prism, condition (5) of Andreev's Theorem is a
consequence of conditions (3) and (4).
\end{prop}

\noindent {\bf Proof:}
Given a quadrilateral face, if the four edges leading from it form a
prismatic 4-circuit, $\Gamma_1$, as depicted on the left hand side of the figure below, clearly condition (5) is a result of condition (4). 
Otherwise, at least one pair of the edges leading from it meet at a
vertex.  If only one pair meets at a point, we have the diagram below in 
the middle.  In this case, the curve $\Gamma_2$ can easily be shown to be a
prismatic 3-circuit, so that
$\alpha_{34} + \alpha_{41} + \beta < \pi$, so that condition (5) is
satisfied because $\alpha_{34}$ and $\alpha_{41}$ cannot both be
$\pi/2$. 

\begin{center}
\begin{picture}(0,0)%
\includegraphics{./prism_circ.pstex}%
\end{picture}%
\setlength{\unitlength}{4144sp}%
\begingroup\makeatletter\ifx\SetFigFont\undefined%
\gdef\SetFigFont#1#2#3#4#5{%
  \reset@font\fontsize{#1}{#2pt}%
  \fontfamily{#3}\fontseries{#4}\fontshape{#5}%
  \selectfont}%
\fi\endgroup%
\begin{picture}(3800,1237)(552,-681)
\put(2381,445){\makebox(0,0)[lb]{\smash{{\SetFigFont{10}{12.0}{\familydefault}{\mddefault}{\updefault}{\color[rgb]{0,0,0}$\beta$}%
}}}}
\put(3967,-174){\makebox(0,0)[lb]{\smash{{\SetFigFont{10}{12.0}{\familydefault}{\mddefault}{\updefault}{\color[rgb]{0,0,0}$e_0$}%
}}}}
\put(2683,-182){\makebox(0,0)[lb]{\smash{{\SetFigFont{10}{12.0}{\familydefault}{\mddefault}{\updefault}{\color[rgb]{0,0,0}$\Gamma_2$}%
}}}}
\put(2119,-632){\makebox(0,0)[lb]{\smash{{\SetFigFont{10}{12.0}{\familydefault}{\mddefault}{\updefault}{\color[rgb]{0,0,0}$\alpha_{4,1}$}%
}}}}
\put(2645,-632){\makebox(0,0)[lb]{\smash{{\SetFigFont{10}{12.0}{\familydefault}{\mddefault}{\updefault}{\color[rgb]{0,0,0}$\alpha_{3,4}$}%
}}}}
\put(1240,-179){\makebox(0,0)[lb]{\smash{{\SetFigFont{10}{12.0}{\familydefault}{\mddefault}{\updefault}{\color[rgb]{0,0,0}$\Gamma_1$}%
}}}}
\put(1165,-629){\makebox(0,0)[lb]{\smash{{\SetFigFont{10}{12.0}{\familydefault}{\mddefault}{\updefault}{\color[rgb]{0,0,0}$\alpha_{3,4}$}%
}}}}
\put(1165,234){\makebox(0,0)[lb]{\smash{{\SetFigFont{10}{12.0}{\familydefault}{\mddefault}{\updefault}{\color[rgb]{0,0,0}$\alpha_{2,3}$}%
}}}}
\put(639,-629){\makebox(0,0)[lb]{\smash{{\SetFigFont{10}{12.0}{\familydefault}{\mddefault}{\updefault}{\color[rgb]{0,0,0}$\alpha_{4,1}$}%
}}}}
\put(634,276){\makebox(0,0)[lb]{\smash{{\SetFigFont{10}{12.0}{\familydefault}{\mddefault}{\updefault}{\color[rgb]{0,0,0}$\alpha_{1,2}$}%
}}}}
\end{picture}%
\end{center}

\noindent
Otherwise, if two pairs of the edges leaving the quadrilateral face meet
at vertices, we have the diagram on the right-hand side.  The only way to
complete this diagram is with the edge labeled $e_0$, resulting in the triangular prism. \Endproof

Hence, we need only check condition (5) for the triangular prism,
which corresponds to the only five-faced $C$.  

Given some $C$, it may be a difficult problem to determine whether $A_C =
\emptyset$ and correspondingly, whether there are any hyperbolic polyhedra
realizing $C$ with non-obtuse dihedral angles.  In fact, for the abstract
polyhedron in the following figure, conditions (2) and (3) imply 
respectively that $\alpha_1+\cdots+\alpha_{12} > 4\pi$ and
$\alpha_1+\cdots+\alpha_{12} < 4\pi$.  So, for this $C$, we have $A_C =
\emptyset$.  However, for more complicated $C$, it can be significantly
harder to determine whether $A_C = \emptyset$.

\begin{center}
\begin{picture}(0,0)%
\epsfig{file=./truncated_cube.pstex}%
\end{picture}%
\setlength{\unitlength}{4144sp}%
\begingroup\makeatletter\ifx\SetFigFont\undefined%
\gdef\SetFigFont#1#2#3#4#5{%
  \reset@font\fontsize{#1}{#2pt}%
  \fontfamily{#3}\fontseries{#4}\fontshape{#5}%
  \selectfont}%
\fi\endgroup%
\begin{picture}(1862,1979)(455,-1172)
\put(1471,370){\makebox(0,0)[lb]{\smash{{\SetFigFont{12}{14.4}{\familydefault}{\mddefault}{\updefault}{\color[rgb]{0,0,0}\small{$\alpha_{10}$}}%
}}}}
\put(1846,464){\makebox(0,0)[lb]{\smash{{\SetFigFont{12}{14.4}{\familydefault}{\mddefault}{\updefault}{\color[rgb]{0,0,0}\small{$\alpha_{12}$}}%
}}}}
\put(1306,-164){\makebox(0,0)[lb]{\smash{{\SetFigFont{12}{14.4}{\familydefault}{\mddefault}{\updefault}{\color[rgb]{0,0,0}\small{$\alpha_9$}}%
}}}}
\put(455,-304){\makebox(0,0)[lb]{\smash{{\SetFigFont{12}{14.4}{\familydefault}{\mddefault}{\updefault}{\color[rgb]{0,0,0}\small{$\alpha_5$}}%
}}}}
\put(2317,-861){\makebox(0,0)[lb]{\smash{{\SetFigFont{12}{14.4}{\familydefault}{\mddefault}{\updefault}{\color[rgb]{0,0,0}\small{$\alpha_3$}}%
}}}}
\put(1428,651){\makebox(0,0)[lb]{\smash{{\SetFigFont{12}{14.4}{\familydefault}{\mddefault}{\updefault}{\color[rgb]{0,0,0}\small{$\alpha_{11}$}}%
}}}}
\put(1377,106){\makebox(0,0)[lb]{\smash{{\SetFigFont{12}{14.4}{\familydefault}{\mddefault}{\updefault}{\color[rgb]{0,0,0}\small{$\alpha_8$}}%
}}}}
\put(994, -6){\makebox(0,0)[lb]{\smash{{\SetFigFont{12}{14.4}{\familydefault}{\mddefault}{\updefault}{\color[rgb]{0,0,0}\small{$\alpha_7$}}%
}}}}
\put(777,-349){\makebox(0,0)[lb]{\smash{{\SetFigFont{12}{14.4}{\familydefault}{\mddefault}{\updefault}{\color[rgb]{0,0,0}\small{$\alpha_{6}$}}%
}}}}
\put(543,-590){\makebox(0,0)[lb]{\smash{{\SetFigFont{12}{14.4}{\familydefault}{\mddefault}{\updefault}{\color[rgb]{0,0,0}\small{$\alpha_4$}}%
}}}}
\put(1890,-863){\makebox(0,0)[lb]{\smash{{\SetFigFont{12}{14.4}{\familydefault}{\mddefault}{\updefault}{\color[rgb]{0,0,0}\small{$\alpha_{2}$}}%
}}}}
\put(1964,-1114){\makebox(0,0)[lb]{\smash{{\SetFigFont{12}{14.4}{\familydefault}{\mddefault}{\updefault}{\color[rgb]{0,0,0}\small{$\alpha_1$}}%
}}}}
\end{picture}%

\end{center}

\noindent
Luckily, there are special cases, including:

\begin{cor} \label{EX3APR}
If there are no prismatic 3-circuits in $C$, there exists a unique
hyperbolic polyhedron realizing C with dihedral angles $2\pi/5$.
\end{cor}

\noindent 
{\bf Proof:}
Since there are no prismatic 3-circuits in $C$, condition (3) of the
theorem is vacuous and clearly $\alpha_i = 2\pi/5$ satisfy conditions (1), (2), (4),
and (5).  \Endproof

\section{Setup of the Proof}\label{SETUP}

Let $C$ be a trivalent abstract polyhedron with $N$ faces. We say that a
hyperbolic polyhedron $P \subset \mathbb{H}^3$ {\it realizes $C$} if there
is a cellular homeomorphism from $C$ to $\partial P$ (i.e., a homeomorphism
mapping faces of $C$ to faces of $P$, edges of $C$ to edges of
$P$, and vertices of $C$ to vertices of $P$). We will call each isotopy
class of cellular homeomorphisms $\phi : C \rightarrow \partial P$ a {\it
marking} on $P$.

We will define ${\cal P}_C$ to be the set of pairs $(P,\phi)$ so that
$\phi$ is a marking with the equivalence relation that $(P,\phi) \sim
(P',\phi ')$ if there exists an isomorphism $\alpha :  \mathbb{H}^3
\rightarrow \mathbb{H}^3$ such that $\alpha(P) = P'$ and both $\phi '$
and $\alpha \circ \phi$ represent the same marking on $P'$.

\begin{prop}
The space ${\cal P}_C$ is a manifold of dimension $3N-6$ (perhaps empty).
\end{prop}

\noindent
{\bf Proof:}
Let ${\cal H}$ be the space of closed half-spaces of $\mathbb{H}^3$; clearly
${\cal H}$ is a 3-dimensional manifold.  Let ${\cal O}_C$ be the set of marked
hyperbolic polyhedra realizing $C$.  Using the marking to number the faces from
$1$ to $N$, an element of ${\cal O}_C$ is an $N$-tuple of half-spaces that
intersect in a polyhedron realizing $C$.  This induces a mapping from ${\cal
O}_C$ to ${\cal H}^N$ whose image is an open set.  We give ${\cal O}_C$ the
topology that makes this mapping from ${\cal O}_C$ into ${\cal H}^N$ a local
homeomorphism.  Since ${\cal H}^N$ is a $3N$-dimensional manifold, ${\cal O}_C$
must be a $3N$-dimensional manifold as well.

If $\alpha(P,\phi) = (P,\phi)$, we have that $\alpha \circ \phi$ is isotopic to
$\phi$ through cellular homeomorphisms.  Hence, the isomorphism $\alpha$ must
fix all vertices of $P$, and consequently restricts to the identity on all
edges and faces.  However, an isomorphism of $\mathbb{H}^3$ which fixes four
non-coplanar points must be the identity.  Therefore $Isom(\mathbb{H}^3)$ acts
freely on ${\cal O}_C$.  This quotient space of this action is ${\cal P}_C$,
hence   ${\cal P}_C$ is a manifold with dimension equal to $dim({\cal O}_C) -
dim(Isom(\mathbb{H}^3)) = 3N-6.$ \Endproof

An $m$-sided polygon $Q \subset \mathbb{H}^2$ with sides $s_i$ supported by lines
$l_i$ for $i=1,\dots,m$ is called a {\em parallelogram} if, after adjoining any
ideal endpoints in $\partial \mathbb{H}^2$ to these sides and lines, $s_i \cap s_j
= \emptyset$ implies $l_i \cap l_j = \emptyset$.  In other words, if two sides
of $Q$ don't meet in $\overline{\mathbb{H}^2}$, then their supporting lines have a
common perpendicular.  We then define ${\cal P}^1_C$ to be the subset of ${\cal
P}_C$ consisting of those polyhedra all of whose faces are parallelograms.  

Let's check that ${\cal P}_C^1$ is an open subset of ${\cal P}_C$.  If $P_{\bf
v_i},$ $i=1,2,3$ are three planes carrying faces of $P \in {\cal P}_C$, then
$\{{\bf v}_i \mbox{,  } i=1,2,3 \}$ will be a linearly independent set spanning
a subspace $V$.  Such a triple of planes has no common intersection point in
$\overline{\mathbb{H}^3}$ if and only if the metric is indefinite when
restricted to $V$, or equivalently, if and only if every vector orthogonal to
$V$ has positive inner product with itself.  This is an open condition on
triples of half-spaces; hence, it is an open condition on ${\cal P}_C$ to
require that the planes supporting three fixed faces of $P$ have no
intersection in $\overline{\mathbb{H}^3}$.  

Requiring that a single face of $P$ be a parallelogram is a finite intersection
of such open conditions, for triples formed of that face and two faces whose
intersections with that face form non-adjacent edges.  For $P$ to lie in ${\cal
P}^1_C$ is a further finite intersection over its faces, so ${\cal P}^1_C$ is
an open subset of ${\cal P}_C$, and hence ${\cal P}^1_C$ is a manifold of
dimension $3N-6$, as well.

In fact, we be most interested in the subset ${\cal P}_C^0$ of polyhedra with
non-obtuse dihedral angles.  Notice that ${\cal P}_C^0$ is not, {\it a priori},
a manifold or even a manifold with boundary.  However, as a consequence of
Proposition \ref{TRIVALENT} (b) and the fact that polygons with non-obtuse interior
angles are parallelograms, we have the inclusion ${\cal P}_C^0 \subset {\cal
P}_C^1$. 

Using the fact that the edge graph of $C$ is trivalent, one can check that
$E$, the number of edges of $C$, is the same as the dimension of ${\cal
P}_C^1$.  Since exactly three edges enter each vertex and each edge enters
exactly two vertices, $3V = 2E$.  The Euler characteristic gives $2=N - E +
V = N - E + 2E/3$ implying $E = 3(N-2)$, the dimension of ${\cal P}_C$ and ${\cal P}_C^1$.

Given any $P \in {\cal P}_C$, let $\alpha(P) =
(\alpha_1,\alpha_2,\alpha_3,...)$ be the $E$-tuple consisting of the dihedral
angles of $P$ at each edge (according to some fixed numbering of the edges of
$C$).  This map $\alpha$ is obviously continuous with respect to the topology
on ${\cal P}_C$, which it inherits from its manifold structure.

Our goal is to prove the following theorem, of which Andreev's Theorem is a 
consequence:

\begin{thm}\label{HOMEO}
For every abstract polyhedron $C$ having more than four faces, the mapping
$\alpha: {\cal P}_C^0 \rightarrow A_C$ is a homeomorphism.
\end{thm}

We will say that {\em Andreev's Theorem holds for $C$} if $\alpha: {\cal
P}_C^0 \rightarrow A_C$ is a homeomorphism for a specific abstract polyhedron
$C$.

\vspace{0.05in}
We begin the proof of Theorem \ref{HOMEO} by checking that $\alpha({\cal
P}_C^0) \subset A_C$ in Section \ref{SEC_INEQSAT}.  In Section \ref{SEC_INJ},
we prove that $\alpha$  restricted to ${\cal P}_C^1$ is injective, and in
Section \ref{SEC_PROP}, we prove that $\alpha$ restricted to ${\cal P}_C^0$ is
proper.  In the beginning of Section \ref{SEC_NONEMPT} we combine these results
to show that $\alpha: {\cal P}_C^0 \rightarrow A_C$ is a homeomorphism onto its
image and that this image is a component of $A_C$.  The remaining, and most
substantial part of Section \ref{SEC_NONEMPT}, is to show that $A_C \neq
\emptyset$ implies ${\cal P}_C^0 \neq \emptyset$.

\section{The inequalities are satisfied.} \label{SEC_INEQSAT}
\begin{prop}
\label{INEQSAT} Given $P \in {\cal P}_C^0$, the dihedral angles $\alpha(P)$
satisfy conditions (1-5).  \end{prop}

We will need the following two lemmas about the basic
properties of hyperbolic geometry.

\begin{lem} \label{INTERSECT}
Suppose that three planes $P_{\bf v_1},P_{\bf v_2},P_{\bf v_3}$ intersect
pairwise in $\mathbb{H}^3$ with non-obtuse dihedral angles $\alpha, \beta$, and
$\gamma$.   Then, $P_{\bf v_1},P_{\bf v_2},P_{\bf v_3}$ intersect at a vertex
in
$\overline{\mathbb{H}^3}$ if and only if
$\alpha+\beta+\gamma \geq \pi.$  The
planes intersect in $\mathbb{H}^3$ if and only if the inequality is strict. 
\end{lem}

\noindent {\bf Proof:}
The planes intersect in a point of $\overline{\mathbb{H}^3}$ if and only if the
inner product is either positive definite or semi-definite on the subspace $V$
spanned by $\{ {\bf v}_i \mbox{,  } i=1,2,3\}$.  In the former case the intersection
point is in $\mathbb{H}^3$, and in the latter case it is in $\partial \mathbb{H}^3$; in
both cases the point is determined by the orthogonal complement of $V$.  The
matrix describing the inner product on $V$ is

\begin{eqnarray*}
\left[
\begin{array}{ccc}
1 & \langle{\bf v_1},{\bf v_2}\rangle &  \langle{\bf v_1},{\bf v_3}\rangle \\
\langle{\bf v_1},{\bf v_2}\rangle & 1 &  \langle{\bf v_2},{\bf v_3}\rangle \\
 \langle{\bf v_1},{\bf v_3}\rangle &  \langle{\bf v_2},{\bf v_3}\rangle & 1\\
\end{array}
\right]
=
\left[ 
\begin{array}{ccc}
1 & -\cos\alpha & -\cos\beta \\
-\cos\alpha & 1 & -\cos\gamma \\
-\cos\beta & -\cos\gamma & 1 \\
\end{array}
\right]
\end{eqnarray*}
where $\alpha,\beta,$ and $\gamma$ are the dihedral angles between the pairs of
faces $(P_{\bf v_1},P_{\bf v_2})$, $(P_{\bf v_1},P_{\bf v_3}),$ and $(P_{\bf
v_2},P_{\bf v_3})$, respectively.

Since the principal minor is positive definite for $0 < \alpha \leq \pi/2$, it is enough to find out when the
determinant
$$1 -2\cos\alpha\cos\beta\cos\gamma -\cos^2\alpha-\cos^2\beta-\cos^2\gamma$$
is non-negative.

A bit of trigonometric trickery (we used complex exponentials) shows that the
expression above is equal to 
\begin{eqnarray}\label{COSEQN}
-4\cos\left(\frac{\alpha+\beta+\gamma}{2}\right)\cos\left(\frac{\alpha-\beta+\gamma}{2}\right)\cos\left(\frac{\alpha+\beta-\gamma}{2}\right)\cos\left(\frac{-\alpha+\beta+\gamma}{2}\right)
\end{eqnarray}

Let $\delta = \alpha+\beta+\gamma$.  When $\delta < \pi$, (\ref{COSEQN}) is strictly
negative; when $\delta = \pi$, (\ref{COSEQN}) is clearly zero; and when
$\delta > \pi$, (\ref{COSEQN}) is
strictly positive.  Hence the inner product on the space spanned by ${\bf v_1},{\bf
v_2},{\bf v_3 }$ is positive semidefinite if and only if $\delta \geq \pi$.  It is
positive definite if and only if $\delta > \pi$.

Then it is easy to see that the three planes $P_{\bf v_1},P_{\bf v_2},P_{\bf v_3}
\subset \mathbb{H}^3$ intersect at a point in $\overline{\mathbb{H}^3}$ if and only if they
intersect pairwise in $\mathbb{H}^3$ and the sum of the dihedral angles
$\delta \geq \pi$.  It is also clear that they intersect at a finite point
if and only if the inequality is strict. \Endproof

\begin{lem} \label{INTERSECT2}
Let $P_1,P_2,P_3 \subset \mathbb{H}^3$ be planes carrying
faces of a polyhedron $P$ that has all dihedral angles $\leq \pi/2$.
\newline (a) If  $P_1,P_2,P_3$ intersect at a point in
$\mathbb{H}^3$, then the point $p
= P_1 \cap P_2 \cap P_3$ is a vertex of $P$.

\noindent
(b) If $P_1,P_2,P_3$ intersect at a point in $\partial \mathbb{H}^3$, then $P$ is
not compact, and the point of intersection is in the closure of $P$.
\end{lem}

\noindent {\bf Proof:}
(a)  Consider what we see in the plane
$P_1$.  Let $H_i$ be the half-space bounded by $P_i$ which contains
the interior of $P$, and let $Q = P_1 \cap H_2 \cap H_3$.  If $p \notin P$,
then let $U$ be the component of $Q-P$ that contains $p$ in its closure.  This is a
non-convex polygon; let $p,p_1,...,p_k$ be its vertices.  The exterior angles of
$U$ at $p_1,...,p_k$ are the angles of the face of $P$ carried by $P_1$, hence
$\leq \pi/2$ by part (b) of Proposition \ref{TRIVALENT}.  See the following figure: 

\begin{center}
\begin{picture}(0,0)%
\epsfig{file=./funny_polygon.pstex}%
\end{picture}%
\setlength{\unitlength}{4144sp}%
\begingroup\makeatletter\ifx\SetFigFont\undefined%
\gdef\SetFigFont#1#2#3#4#5{%
  \reset@font\fontsize{#1}{#2pt}%
  \fontfamily{#3}\fontseries{#4}\fontshape{#5}%
  \selectfont}%
\fi\endgroup%
\begin{picture}(1945,2165)(983,-1633)
\put(984,-1183){\makebox(0,0)[lb]{\smash{{\SetFigFont{7}{8.4}{\familydefault}{\mddefault}{\updefault}{\color[rgb]{0,0,0}$p_1$}%
}}}}
\put(1800,-1606){\makebox(0,0)[lb]{\smash{{\SetFigFont{7}{8.4}{\familydefault}{\mddefault}{\updefault}{\color[rgb]{0,0,0}$P$}%
}}}}
\put(1776,-680){\makebox(0,0)[lb]{\smash{{\SetFigFont{7}{8.4}{\familydefault}{\mddefault}{\updefault}{\color[rgb]{0,0,0}$U$}%
}}}}
\put(1120,-1268){\makebox(0,0)[lb]{\smash{{\SetFigFont{7}{8.4}{\familydefault}{\mddefault}{\updefault}{\color[rgb]{0,0,0}$p_2$}%
}}}}
\put(1695,248){\makebox(0,0)[lb]{\smash{{\SetFigFont{7}{8.4}{\familydefault}{\mddefault}{\updefault}{\color[rgb]{0,0,0}$p$}%
}}}}
\put(2373,-1331){\makebox(0,0)[lb]{\smash{{\SetFigFont{7}{8.4}{\familydefault}{\mddefault}{\updefault}{\color[rgb]{0,0,0}$p_{k-1}$}%
}}}}
\put(2698,-1114){\makebox(0,0)[lb]{\smash{{\SetFigFont{7}{8.4}{\familydefault}{\mddefault}{\updefault}{\color[rgb]{0,0,0}$p_k$}%
}}}}
\end{picture}%

\end{center}

Suppose that $\alpha_1,...\alpha_k$ are the angles of $P$ at $p_1,...,p_k$, and let
$\alpha$ be the angle at $p$.  Then the Gauss-Bonnet formula tells us that:
$$(\pi-\alpha)+\alpha_1-((\pi-\alpha_2)+\cdots+(\pi-\alpha_{k-1})) +\alpha_k
-\Area(U) = 2\pi,$$
which can be rearranged to read
$$(\alpha_1+\alpha_k-\pi) -\alpha -\sum_{j=2}^{k-1} (\pi-\alpha_j) = \Area(U).$$
This is clearly a contradiction.  All of the terms on the left are
non-positive, and $\Area(U) > 0$.  If $p$ is at infinity (i.e., $\alpha = 0$),
this expression is still a contradiction, proving part (b). \Endproof

\noindent {\bf Proof of Proposition \ref{INEQSAT}:}
For condition (1), notice that if two adjacent faces intersect at dihedral angle $0$,
they intersect at a point at infinity.  If this were the case, $P$ would be non-compact.

For condition (2), let $x$ be a vertex of $P$.  Since $P$ is compact, $x \in
\mathbb{H}^3$ and by Lemma \ref{INTERSECT}, the sum of the dihedral angles
between the three planes intersecting at $x$ must be $> \pi$.

For condition (3), note first that by Lemma \ref{INTERSECT} if three faces
forming a 3-circuit have dihedral angles summing to a number $\geq \pi$, then
they meet in $\overline{\mathbb{H}^3}$.  If they meet at a point in
$\mathbb{H}^3$, by Lemma \ref{INTERSECT2}(a) this point is a vertex of
$P$, so these three faces do not form a prismatic 3-circuit.  Alternatively, if
the three planes meet in $\partial \mathbb{H}^3$ by  Lemma \ref{INTERSECT2}(b),
then $P$ is non-compact, contrary to assumption.  Hence, any three
faces forming a prismatic 3-circuit in $P$ must have dihedral angles summing to
$< \pi$.

For condition (4), let $H_{\bf v_1},H_{\bf v_2},H_{\bf v_3},H_{\bf v_4}$ be
half spaces corresponding to the faces which form a prismatic 4-circuit;
obviously condition (4) is satisfied unless all of the dihedral angles are
$\pi/2$, so we suppose that they are.  We will assume the normalization
$\langle{\bf v_i},{\bf v_i}\rangle = 1$ for each ${\bf i}$.  The Gram matrix $Q
= \left[\langle{\bf v_i},{\bf v_j}\rangle \right]_{i,j} =$

{\small
\begin{eqnarray*}
 \left[ 
\begin{array}{cccc}
1 & 0 & \langle{\bf v_1},{\bf v_3}\rangle & 0 \\
0 & 1 & 0 &  \langle{\bf v_2},{\bf v_4}\rangle \\
\langle{\bf v_3},{\bf v_1}\rangle & 0 & 1 & 0 \\
0 & \langle{\bf v_4},{\bf v_2}\rangle & 0 & 1 \\
\end{array}
\right]
\end{eqnarray*}
}

\noindent
has determinant $0$ if the ${\bf v}$'s are linearly dependent, and otherwise
represents the inner product of $E^{3,1}$ and hence has negative determinant.  In
both cases we have
$$\det Q = (1 - \langle{\bf v_1},{\bf v_3}\rangle^2)(1 - \langle{\bf v_2},{\bf
v_4}\rangle^2) \leq 0.$$
So  $\langle{\bf v_1},{\bf v_3}\rangle^2 \leq 1$ and $\langle{\bf v_2},{\bf
v_4}\rangle^2 \geq 1$ 
or vice versa (perhaps one or both are equalities).
This
means that one of the opposite pairs of faces of the 4-circuit intersect, perhaps
at a point at infinity.  We can suppose that this pair is $H_{\bf v_1}$ and $H_{\bf
v_3}$.

If  $H_{\bf v_1}$ and $H_{\bf v_3}$ intersect in $\mathbb{H}^3$, they do so with
positive dihedral angle.  Since   $H_{\bf v_2}$ intersects each $H_{\bf v_1}$ and
$H_{\bf v_3}$ orthogonally, the three faces pairwise intersect and
have dihedral angle sum $> \pi$.  By Lemmas \ref{INTERSECT} and \ref{INTERSECT2}
these three faces intersect at a point in $\mathbb{H}^3$ which is a vertex of $P$.
In this case, the 4-circuit $H_{\bf v_1},H_{\bf v_2},H_{\bf v_3},H_{\bf v_4}$ is
not prismatic.

Otherwise, $H_{\bf v_1}$ and $H_{\bf v_3}$ intersect at a point at infinity.  In
this case, since $H_{\bf v_2}$ intersects each $H_{\bf v_1}$ and
$H_{\bf v_3}$ with dihedral angle $\pi/2$ the three faces intersect at this point
at infinity by Lemma \ref{INTERSECT} and then by Lemma \ref{INTERSECT2} $P$ is not
compact, contrary to assumption.

Hence, if $H_{\bf v_1},H_{\bf v_2},H_{\bf v_3},H_{\bf v_4}$ forms a prismatic
4-circuit, the sum of the dihedral angles cannot be $2\pi$.  

For condition (5), suppose that the quadrilateral is formed by edges $e_1,
e_2, e_3, e_4$.  Violation of one of the inequalities would give that the
dihedral angles at each of the edges $e_{ij}$ leading to the quadrilateral
is $\pi/2$ and that the dihedral angles at two of the opposite edges of
the quadrilateral are $\pi/2$.   See the diagram below:

\begin{center}
\begin{picture}(0,0)%
\includegraphics{./quadrilateral.pstex}%
\end{picture}%
\setlength{\unitlength}{4144sp}%
\begingroup\makeatletter\ifx\SetFigFont\undefined%
\gdef\SetFigFont#1#2#3#4#5{%
  \reset@font\fontsize{#1}{#2pt}%
  \fontfamily{#3}\fontseries{#4}\fontshape{#5}%
  \selectfont}%
\fi\endgroup%
\begin{picture}(1332,1312)(496,-986)
\put(752,-612){\makebox(0,0)[lb]{\smash{{\SetFigFont{9}{10.8}{\familydefault}{\mddefault}{\updefault}{\color[rgb]{0,0,0}\small{$\pi/2$}}%
}}}}
\put(843,-142){\makebox(0,0)[lb]{\smash{{\SetFigFont{9}{10.8}{\familydefault}{\mddefault}{\updefault}{\color[rgb]{0,0,0}\small{$\gamma$}}%
}}}}
\put(1399,-629){\makebox(0,0)[lb]{\smash{{\SetFigFont{9}{10.8}{\familydefault}{\mddefault}{\updefault}{\color[rgb]{0,0,0}\small{$\beta$}}%
}}}}
\put(1191,-941){\makebox(0,0)[lb]{\smash{{\SetFigFont{9}{10.8}{\familydefault}{\mddefault}{\updefault}{\color[rgb]{0,0,0}\small{$\pi/2$}}%
}}}}
\put(1191,205){\makebox(0,0)[lb]{\smash{{\SetFigFont{9}{10.8}{\familydefault}{\mddefault}{\updefault}{\color[rgb]{0,0,0}\small{$\pi/2$}}%
}}}}
\put(496,-281){\makebox(0,0)[lb]{\smash{{\SetFigFont{9}{10.8}{\familydefault}{\mddefault}{\updefault}{\color[rgb]{0,0,0}\small{$\pi/2$}}%
}}}}
\put(1608,-281){\makebox(0,0)[lb]{\smash{{\SetFigFont{9}{10.8}{\familydefault}{\mddefault}{\updefault}{\color[rgb]{0,0,0}\small{$\pi/2$}}%
}}}}
\put(1364,-136){\makebox(0,0)[lb]{\smash{{\SetFigFont{9}{10.8}{\familydefault}{\mddefault}{\updefault}{\color[rgb]{0,0,0}\small{$\pi/2$}}%
}}}}
\end{picture}%

\end{center}

Each vertex of this quadrilateral had three incident edges labeled $e_i$, $e_j$, and $e_{ij}$.
Violation of the inequality gives that $\alpha_{ij} = \pi/2$ and
either $\alpha_i = \pi/2$ or $\alpha_j = \pi/2$.  Using Equation (\ref{TLC}) from section 1, we
see that each face angle in the quadrilateral must be $\pi/2$.
So, we have that each of the face angles of the quadrilateral is $\pi/2$,
which is a contradiction to the Gauss-Bonnet Theorem.  Hence both of the
inequalities in condition (5) must be satisfied. 

This was the last step
in proving Proposition \ref{INEQSAT}.
\Endproof

\section{The mapping $\alpha: {\cal P}_C^1 \rightarrow \mathbb{R}^E$ is injective.} \label{SEC_INJ}

Recall from section \ref{SETUP} 
an $m$-sided polygon $Q \subset \mathbb{H}^2$ with sides $s_i$ supported by lines
$l_i$ for $i=1,\dots,m$ is called a {\em parallelogram} if, after adjoining any
ideal endpoints in $\partial \mathbb{H}^2$ to these sides and lines, $s_i \cap s_j
= \emptyset$ implies $l_i \cap l_j = \emptyset$. 

In section \ref{SETUP} that we defined ${\cal P}^1_C$ to be the
space of polyhedra realizing $C$ whose faces are parallelograms.
We then checked that ${\cal P}^1_C$ is an open subset of ${\cal P}_C$ and
that ${\cal P}^0_C \subset {\cal P}^1_C$.  The goal of this section is to prove:

\begin{prop} \label{INJECTIVE}
The mapping $\alpha : {\cal P}^1_C \rightarrow \mathbb{R}^E$ is injective.
\end{prop}

\noindent {\bf Proof:}
Suppose that $P, P' \in {\cal P}^1_C$ are two polyhedra such that
$\alpha(P) = \alpha(P')$.  We can label each edge $e$ of $C$ by $-,0,$
or $+$ depending on whether the length of $e$ in $P'$ is less than, equal
to, or greater than the length of $e$ in $P$.

We will prove that if $\alpha(P') = \alpha(P')$ then each pair of
corresponding edges has the same length.  This gives that the faces of $P$
and $P'$ are congruent since the face angles are determined by the dihedral
angles (Proposition \ref{TRIVALENT}).  Then, since $P$ and $P'$
have congruent faces and the same dihedral angles, they are themselves congruent.

Each edge of $C$ corresponds to a unique edge of the dual complex, $C^*$, which
we label with $-,0$, or $+$, accordingly.  Consider the graph $\cal{G}$
consisting of the edges of $C^*$ labeled either $+$ or $-$, but not $0$,
together with the vertices incident to these edges.  Since $C^*$ is a
simplicial complex on $\mathbb{S}^2$, $\cal{G}$ is a simple planar graph.
(Here, simple means that there is at most one edge between any distinct pair of
vertices and no edges from a vertex to itself.)  We assume that
$\cal{G}$ is non-empty, in order to find a contradiction.

\begin{prop} \label{CAUCHY}
Let $\cal{G}$ be a simple planar graph whose edges are labeled with $+$ and
$-$.  There is a vertex of $\cal{G}$ with at most two sign changes when
following the cyclic order of the edges meeting at that vertex.
\end{prop}

\noindent Proposition \ref{CAUCHY} provides the global statement necessary for
Cauchy's rigidity theorem on Euclidean polyhedra, see \cite[Chapter 12]{PFB},
and also the global statement necessary here.  The proof, see \cite[page
68]{PFB}, is a clever, yet elementary counting argument, combined with Euler's
Formula.

Therefore, at some vertex of $\cal{G}$, there are either zero or two changes of
sign when following the cyclic order of edges meeting at that vertex.  In the
case that there are zero changes in sign, we may assume, without loss in
generality that all of the signs at this vertex are $+$'s, by switching the
roles of $P$ and $P'$, if necessary.   Thus, in either case, there is a face $F$
of $C$ not marked entirely with $0$'s so that, once the edges labeled $-$ are
removed from $\partial F$, the edges labeled $+$ all lie in the same component
of what remains.

Let $Q$ and $Q'$ be the faces in $P$ and $P'$ corresponding to $F$.  By the
assumption that $P, P' \in {\cal P}^1_C$, $Q$ and $Q'$ are parallelograms and
because $P$ and $P'$ have the same dihedral angles, $Q$ and $Q'$ must have the
same face angles.  We will now show that $Q$ and $Q'$ cannot have side lengths
differing according to the distribution of $+$'s and $-$'s on $\partial F$ that
was deduced above.

\vspace{.1in}
The following lemma, from \cite[page 422]{AND} but with a new proof, shows that
(roughly speaking) stretching edges in a piece of the boundary will pull apart
the two edges at the ends of that piece.  It is important to keep in mind that
the parallelograms $R$ and $R'$ need not be compact and need not have finite
volume, since there are no restrictions on whether the first and $m$-th sides
intersect.

\begin{lem}{\bf (Andreev's Auxiliary Lemma)} \label{AUX}
Let $R$ and $R'$ be $m$-sided parallelograms further assume that
$R$ and $R'$ have finite vertices $A_i = s_i \cap s_{i+1}$ and $A_i' = s_i' \cap s_{i+1}'$ for $i=1,\cdots m-1$.
If
\begin{itemize}

\item The interior angle at vertex $A_i$ and vertex $A_i'$ are equal for $i=2,\cdots m-1$, and
\item $|s_j| \leq |s_j'|$ for $j=2,\cdots,m-1$,
\end{itemize}

with at least one of the inequalities strict.  Then
$\langle {\bf v}_1, {\bf v}_m \rangle > \langle {\bf v}_1', {\bf v}_m' \rangle$,
where ${\bf v}_i$ and ${\bf v}_m'$ are the outward pointing normal to the
edge $s_i$ of $R$ and $R'$, respectively.
\end{lem}

\noindent
{\bf Proof:}
We will prove the lemma first in the case where the side lengths differ only at
one side $|s_j| < |s_j'|$ and then  observe that the resulting polygon again
satisfies the hypotheses of the lemma so that one can repeat as necessary for
each pair of sides that differ in length.

We can situate side $s_j$ on the line $x_2=0$ centered at $(1,0,0)$ within the
upper sheet of the hyperboloid $-x_0^2 +x_1^2+x_2^2 = -1$ and assume that $R$
is entirely ``above'' this line, that is at points with $x_2 \geq 0$. We also
assume that the sides of $R$ are labeled counterclockwise, i.e., $s_{i+1}$ is
counterclockwise from $s_i$ for each $i$.

Applying the isometry:
\begin{eqnarray*}
I(t) = \left[\begin{array}{ccc} \cosh(t) & \sinh(t) & 0 \\ \sinh(t) & \cosh(t) & 0 \\ 0 & 0 & 1 \end{array} \right]
\end{eqnarray*}

\noindent
to the sides $s_i$ with index $i > j$ for $t > 0$ performs the desired deformation of $R$.

One can check that
\begin{eqnarray*}
\frac{d}{dt} \langle {\bf v}_1, I(t) {\bf v}_m \rangle  = \sinh(t) (v_{m1}v_{11}-v_{m0} v_{10})-\cosh(t)(v_{m1} v_{10}-v_{m0} v_{11}),
\end{eqnarray*}
\noindent
if we write ${\bf v_1} = (v_{10},v_{11},v_{12})$ and ${\bf v_m} = (v_{m0},v_{m1},v_{m2})$.
Since $(1,0,0)$ is in the interior of $s_j$ we must have that $\langle (1,0,0),
{\bf v}_1 \rangle < 0$, which is equivalent to $v_{10} > 0$.  For the same
reason we also have $v_{m0} > 0$.

We will first check that this derivative $\frac{d}{dt} \langle {\bf v}_1, I(t) {\bf v}_m \rangle$ is negative at $t=0$, or equivalently that:
\begin{eqnarray*}
0 < \det\left[\begin{array}{ccc}
v_{10} & 0 & v_{m0} \\
v_{11} & 0 & v_{m1} \\
v_{12} & -1 & v_{m2} \end{array} \right]
= v_{m1}v_{10} - v_{m0}v_{11}.
\end{eqnarray*}

Imagine the three column vectors, in order, in a right-handed coordinate system
with the $x_0$ axis pointing up, the $x_1$ axis pointing forward and the $x_2$
axis pointing to the right.  Because of the choice of orientation made above, the duals
to the geodesics carrying $s_1,s_i,s_m$, when viewed from ``above'' and in that
order, will also turn counter-clockwise.  Since these duals all have
non-negative $x_0$ coordinates, and two of these are positive, they form a
right-handed frame, and the determinant is therefore positive.

\vspace{0.1in}
We now check that $\frac{d}{dt} \langle {\bf v}_1, I(t) {\bf v}_m \rangle < 0$ for an arbitrary $t > 0$.
Because $v_{m1} v_{10} - v_{m0} v_{11} > 0$ this is equivalent to:
\begin{eqnarray}\label{WANT}
\frac{v_{m1}v_{11}-v_{m0} v_{10}}{v_{m1} v_{10} - v_{m0} v_{11}}  < \frac{\cosh(t)}{\sinh(t)}.
\end{eqnarray}
\vspace{0.1in}
Furthermore, since $\frac{\cosh(t)}{\sinh(t)} > 1$ it is sufficient to show that $\frac{v_{m1}v_{11}-v_{m0} v_{10}}{v_{m1} v_{10} - v_{m0} v_{11}} \leq 1$.

Since $R$ is a parallelogram oriented counter-clockwise, $l_1$ and $l_j$ cannot
intersect in $\mathbb{H}^2$ to the right of $O$, since only $l_{j+1}$ can
intersect $l_j$ there; nor can $l_1$ and $l_j$ be asymptotic at the right ideal
endpoint of $l_j$ (represented by the vector $(1,1,0) \in E^{2,1}$).  This
means that we never intersect the boundary of the half-space in $E^{2,1}$
corresponding to ${\bf v}_1$ when we move in a straight line from
$(1,0,0)$ to $(1,1,0)$, forcing the latter point to lie in the interior of that
half-space.  In other words, $0 > \langle (1,1,0),(v_{10},v_{11},v_{12})\rangle
= -v_{10} + v_{11}.$

The following diagrams (radially projecting $E^{2,1}$ to $\{x_0 = 1\}$)
illustrate some of the ways $l_1$, $l_j$, and $l_m$ can be arranged, but are
not intended to be a comprehensive list.  Configuration  1 is allowed (but will
only occur if $2 < j < m-1$).  Configuration 2 is only allowed when $m=j+1$.
Configuration 3 is ruled out by the intersection of $l_1$ with $l_j$ to the
right of $O$, violating the parallelogram condition, and Configuration 4 is
forbidden by the orientation condition.

\begin{center}
\begin{picture}(0,0)%
\epsfig{file=./andreev_inj.pstex}%
\end{picture}%
\setlength{\unitlength}{3947sp}%
\begingroup\makeatletter\ifx\SetFigFont\undefined%
\gdef\SetFigFont#1#2#3#4#5{%
  \reset@font\fontsize{#1}{#2pt}%
  \fontfamily{#3}\fontseries{#4}\fontshape{#5}%
  \selectfont}%
\fi\endgroup%
\begin{picture}(4988,4421)(211,-3699)
\put(921,-1889){\makebox(0,0)[lb]{\smash{{\SetFigFont{8}{9.6}{\familydefault}{\mddefault}{\updefault}{\color[rgb]{0,0,0}$l_m$}%
}}}}
\put(650,-1403){\makebox(0,0)[lb]{\smash{{\SetFigFont{8}{9.6}{\familydefault}{\mddefault}{\updefault}{\color[rgb]{0,0,0}Configuration 1}%
}}}}
\put(723,-3661){\makebox(0,0)[lb]{\smash{{\SetFigFont{8}{9.6}{\familydefault}{\mddefault}{\updefault}{\color[rgb]{0,0,0}Configuration 3}%
}}}}
\put(1066,471){\makebox(0,0)[lb]{\smash{{\SetFigFont{8}{9.6}{\familydefault}{\mddefault}{\updefault}{\color[rgb]{0,0,0}$l_m$}%
}}}}
\put(814,-348){\makebox(0,0)[lb]{\smash{{\SetFigFont{8}{9.6}{\familydefault}{\mddefault}{\updefault}{\color[rgb]{0,0,0}$l_j$}%
}}}}
\put(3491,-1373){\makebox(0,0)[lb]{\smash{{\SetFigFont{8}{9.6}{\familydefault}{\mddefault}{\updefault}{\color[rgb]{0,0,0}Configuration 2}%
}}}}
\put(4177,455){\makebox(0,0)[lb]{\smash{{\SetFigFont{8}{9.6}{\familydefault}{\mddefault}{\updefault}{\color[rgb]{0,0,0}$l_m$}%
}}}}
\put(3699,-317){\makebox(0,0)[lb]{\smash{{\SetFigFont{8}{9.6}{\familydefault}{\mddefault}{\updefault}{\color[rgb]{0,0,0}$l_j$}%
}}}}
\put(3395,-3668){\makebox(0,0)[lb]{\smash{{\SetFigFont{8}{9.6}{\familydefault}{\mddefault}{\updefault}{\color[rgb]{0,0,0}Configuration 4}%
}}}}
\put(3880,-2871){\makebox(0,0)[lb]{\smash{{\SetFigFont{8}{9.6}{\familydefault}{\mddefault}{\updefault}{\color[rgb]{0,0,0}$O$}%
}}}}
\put(3901,-1833){\makebox(0,0)[lb]{\smash{{\SetFigFont{8}{9.6}{\familydefault}{\mddefault}{\updefault}{\color[rgb]{0,0,0}$l_1$}%
}}}}
\put(3717,-2625){\makebox(0,0)[lb]{\smash{{\SetFigFont{8}{9.6}{\familydefault}{\mddefault}{\updefault}{\color[rgb]{0,0,0}$l_j$}%
}}}}
\put(946,-2651){\makebox(0,0)[lb]{\smash{{\SetFigFont{8}{9.6}{\familydefault}{\mddefault}{\updefault}{\color[rgb]{0,0,0}$l_j$}%
}}}}
\put(3631,483){\makebox(0,0)[lb]{\smash{{\SetFigFont{8}{9.6}{\familydefault}{\mddefault}{\updefault}{\color[rgb]{0,0,0}$l_1$}%
}}}}
\put(1116,-2887){\makebox(0,0)[lb]{\smash{{\SetFigFont{8}{9.6}{\familydefault}{\mddefault}{\updefault}{\color[rgb]{0,0,0}$O$}%
}}}}
\put(1049,-587){\makebox(0,0)[lb]{\smash{{\SetFigFont{8}{9.6}{\familydefault}{\mddefault}{\updefault}{\color[rgb]{0,0,0}$O$}%
}}}}
\put(3916,-567){\makebox(0,0)[lb]{\smash{{\SetFigFont{8}{9.6}{\familydefault}{\mddefault}{\updefault}{\color[rgb]{0,0,0}$O$}%
}}}}
\put(405,-2066){\makebox(0,0)[lb]{\smash{{\SetFigFont{8}{9.6}{\familydefault}{\mddefault}{\updefault}{\color[rgb]{0,0,0}$l_1$}%
}}}}
\put(3413,-1905){\makebox(0,0)[lb]{\smash{{\SetFigFont{8}{9.6}{\familydefault}{\mddefault}{\updefault}{\color[rgb]{0,0,0}$l_m$}%
}}}}
\put(596,362){\makebox(0,0)[lb]{\smash{{\SetFigFont{8}{9.6}{\familydefault}{\mddefault}{\updefault}{\color[rgb]{0,0,0}$l_1$}%
}}}}
\end{picture}%

\end{center}

An analogous argument shows $l_j$ and $l_m$ cannot intersect in $H^2$ to the
left of $O$ and that they cannot be asymptotic at $(1,-1,0)$, so $(1,-1,0)$ is
contained in the interior of the half-space dual to
$(v_{m0},v_{m1},v_{m2})$, or in other words, $-v_{m0} - v_{m1} < 0$.

Combining these two observations yields $0 > (v_{m1} + v_{m0})(v_{11} - v_{10})
= v_{m1}v_{11}  +v_{m0} v_{11}  - v_{m1}v_{10}  - v_{m0}v_{10} $ which is
equivalent to $v_{m1}v_{11} - v_{m0}v_{10}  < v_{m1}v_{10}  - v_{m0}v_{11}$.
In combination with the fact that the right hand side of this inequality is
positive, this shows that Equation (\ref{WANT}) holds.

\vspace{.1in}

For $i=1,\cdots m-1$ the adjacent sides $s_i$ and $s_{i+1}$ of $R$ continue to
intersect at finite vertices with the same interior angles as before this
deformation.  Applying what we have just proved to an appropriate sub-polygon
of $R$ we can see that $\langle  {\bf v}_k, {\bf v}_l \rangle$ is
non-increasing for pairs of sides of $s_k$ and $s_l$ that did not intersect
before this deformation.  Because these sides satisfied $\langle  {\bf v}_k,
{\bf v}_l \rangle < 0$ before the deformation, they continue to do so, and the
resulting polygon satisfies the hypotheses of Lemma \ref{AUX}.  Hence, one can
increase the lengths of sides $s_j$ sequentially, in order to prove Lemma
\ref{AUX} in full generality.

\Endproof

We continue the proof of Proposition \ref{INJECTIVE}.  Suppose that $F$ has $n$
sides.   We can renumber the sides of $F$ so that the second through $(m-1)$-st
sides of $F$ are labeled with $+$'s and $0$'s and at least one of them is
labeled with a $+$, so that the first and $m$-th sides are arbitrarily labeled,
and, if $m < n$, so that the remaining sides are all labeled with $-$'s and
$0$'s.  There is usually more than one way to do this (often with differing
values of $m$), any of which will suffice.  See the diagrams below for
examples.  In the former, there are many alternate choices; in the latter,
there is essentially one choice, up to combinatorial symmetry.

\begin{center}
\begin{picture}(0,0)%
\epsfig{file=./plus_minus.pstex}%
\end{picture}%
\setlength{\unitlength}{3947sp}%
\begingroup\makeatletter\ifx\SetFigFont\undefined%
\gdef\SetFigFont#1#2#3#4#5{%
  \reset@font\fontsize{#1}{#2pt}%
  \fontfamily{#3}\fontseries{#4}\fontshape{#5}%
  \selectfont}%
\fi\endgroup%
\begin{picture}(5032,2069)(376,-4371)
\put(4273,-4310){\makebox(0,0)[lb]{\smash{{\SetFigFont{8}{9.6}{\familydefault}{\mddefault}{\updefault}{\color[rgb]{0,0,0}$m=n$}%
}}}}
\put(4230,-2508){\makebox(0,0)[lb]{\smash{{\SetFigFont{8}{9.6}{\familydefault}{\mddefault}{\updefault}{\color[rgb]{0,0,0}$+$}%
}}}}
\put(5171,-3644){\makebox(0,0)[lb]{\smash{{\SetFigFont{8}{9.6}{\familydefault}{\mddefault}{\updefault}{\color[rgb]{0,0,0}$+$}%
}}}}
\put(3671,-3494){\makebox(0,0)[lb]{\smash{{\SetFigFont{8}{9.6}{\familydefault}{\mddefault}{\updefault}{\color[rgb]{0,0,0}$+$}%
}}}}
\put(3596,-2819){\makebox(0,0)[lb]{\smash{{\SetFigFont{8}{9.6}{\familydefault}{\mddefault}{\updefault}{\color[rgb]{0,0,0}$+$}%
}}}}
\put(5021,-2819){\makebox(0,0)[lb]{\smash{{\SetFigFont{8}{9.6}{\familydefault}{\mddefault}{\updefault}{\color[rgb]{0,0,0}$+$}%
}}}}
\put(4271,-2669){\makebox(0,0)[lb]{\smash{{\SetFigFont{8}{9.6}{\familydefault}{\mddefault}{\updefault}{\color[rgb]{0,0,0}$2$}%
}}}}
\put(3821,-2894){\makebox(0,0)[lb]{\smash{{\SetFigFont{8}{9.6}{\familydefault}{\mddefault}{\updefault}{\color[rgb]{0,0,0}$3$}%
}}}}
\put(3896,-3419){\makebox(0,0)[lb]{\smash{{\SetFigFont{8}{9.6}{\familydefault}{\mddefault}{\updefault}{\color[rgb]{0,0,0}$4$}%
}}}}
\put(4346,-3719){\makebox(0,0)[lb]{\smash{{\SetFigFont{8}{9.6}{\familydefault}{\mddefault}{\updefault}{\color[rgb]{0,0,0}$5$}%
}}}}
\put(4391,-3940){\makebox(0,0)[lb]{\smash{{\SetFigFont{8}{9.6}{\familydefault}{\mddefault}{\updefault}{\color[rgb]{0,0,0}$+$}%
}}}}
\put(4807,-2967){\makebox(0,0)[lb]{\smash{{\SetFigFont{8}{9.6}{\familydefault}{\mddefault}{\updefault}{\color[rgb]{0,0,0}$1$}%
}}}}
\put(4470,-3451){\makebox(0,0)[lb]{\smash{{\SetFigFont{8}{9.6}{\familydefault}{\mddefault}{\updefault}{\color[rgb]{0,0,0}$m=n=6$}%
}}}}
\put(1501,-2386){\makebox(0,0)[lb]{\smash{{\SetFigFont{8}{9.6}{\familydefault}{\mddefault}{\updefault}{\color[rgb]{0,0,0}$+$}%
}}}}
\put(751,-2536){\makebox(0,0)[lb]{\smash{{\SetFigFont{8}{9.6}{\familydefault}{\mddefault}{\updefault}{\color[rgb]{0,0,0}$+$}%
}}}}
\put(376,-2836){\makebox(0,0)[lb]{\smash{{\SetFigFont{8}{9.6}{\familydefault}{\mddefault}{\updefault}{\color[rgb]{0,0,0}$0$}%
}}}}
\put(376,-3361){\makebox(0,0)[lb]{\smash{{\SetFigFont{8}{9.6}{\familydefault}{\mddefault}{\updefault}{\color[rgb]{0,0,0}$+$}%
}}}}
\put(676,-3886){\makebox(0,0)[lb]{\smash{{\SetFigFont{8}{9.6}{\familydefault}{\mddefault}{\updefault}{\color[rgb]{0,0,0}$0$}%
}}}}
\put(601,-2911){\makebox(0,0)[lb]{\smash{{\SetFigFont{8}{9.6}{\familydefault}{\mddefault}{\updefault}{\color[rgb]{0,0,0}$3$}%
}}}}
\put(826,-3661){\makebox(0,0)[lb]{\smash{{\SetFigFont{8}{9.6}{\familydefault}{\mddefault}{\updefault}{\color[rgb]{0,0,0}$5$}%
}}}}
\put(901,-2686){\makebox(0,0)[lb]{\smash{{\SetFigFont{8}{9.6}{\familydefault}{\mddefault}{\updefault}{\color[rgb]{0,0,0}$2$}%
}}}}
\put(1555,-2579){\makebox(0,0)[lb]{\smash{{\SetFigFont{8}{9.6}{\familydefault}{\mddefault}{\updefault}{\color[rgb]{0,0,0}$1$}%
}}}}
\put(2331,-2753){\makebox(0,0)[lb]{\smash{{\SetFigFont{8}{9.6}{\familydefault}{\mddefault}{\updefault}{\color[rgb]{0,0,0}$-$}%
}}}}
\put(1347,-4023){\makebox(0,0)[lb]{\smash{{\SetFigFont{8}{9.6}{\familydefault}{\mddefault}{\updefault}{\color[rgb]{0,0,0}$+$}%
}}}}
\put(1268,-3791){\makebox(0,0)[lb]{\smash{{\SetFigFont{8}{9.6}{\familydefault}{\mddefault}{\updefault}{\color[rgb]{0,0,0}$m=6$}%
}}}}
\put(1874,-2836){\makebox(0,0)[lb]{\smash{{\SetFigFont{8}{9.6}{\familydefault}{\mddefault}{\updefault}{\color[rgb]{0,0,0}$n=9$}%
}}}}
\put(630,-3341){\makebox(0,0)[lb]{\smash{{\SetFigFont{8}{9.6}{\familydefault}{\mddefault}{\updefault}{\color[rgb]{0,0,0}$4$}%
}}}}
\put(1881,-3632){\makebox(0,0)[lb]{\smash{{\SetFigFont{8}{9.6}{\familydefault}{\mddefault}{\updefault}{\color[rgb]{0,0,0}$7$}%
}}}}
\put(2098,-3359){\makebox(0,0)[lb]{\smash{{\SetFigFont{8}{9.6}{\familydefault}{\mddefault}{\updefault}{\color[rgb]{0,0,0}$8$}%
}}}}
\put(2340,-3438){\makebox(0,0)[lb]{\smash{{\SetFigFont{8}{9.6}{\familydefault}{\mddefault}{\updefault}{\color[rgb]{0,0,0}$0$}%
}}}}
\put(2033,-3791){\makebox(0,0)[lb]{\smash{{\SetFigFont{8}{9.6}{\familydefault}{\mddefault}{\updefault}{\color[rgb]{0,0,0}$-$}%
}}}}
\put(1142,-4340){\makebox(0,0)[lb]{\smash{{\SetFigFont{8}{9.6}{\familydefault}{\mddefault}{\updefault}{\color[rgb]{0,0,0}$m < n$}%
}}}}
\end{picture}%
\end{center}

Because $Q$ and $Q'$ are parallelograms, the (possibly non-compact) polygons
bounded by the union of sides $s_1,\cdots,s_m$ from $Q$ and the union
of the sides $s_1',\cdots,s_m'$ from $Q'$ satisfy the hypotheses of Lemma \ref{AUX}.
Hence, if we denote the outward pointing normals to $s_1$ and $s_m$ in $Q$ by
${\bf v}_1$ and ${\bf v}_m$ and in $Q'$ by ${\bf v}_1'$ and ${\bf v}_m '$, the
lemma guarantees that $\langle {\bf v}_1', {\bf v}_m' \rangle > \langle {\bf
v}_1, {\bf v}_m \rangle$.

However, either $n=m$ so that the first and $m$-th sides of $F$  meet at a
vertex giving $\langle {\bf v}_1', {\bf v}_m' \rangle = \langle {\bf v}_1, {\bf
v}_m \rangle$, or, if $n > m$ all of the sides in $F$ with index greater than $m$ are
labeled $0$ or $-$.  Applying Lemma \ref{AUX} to the polygons bounded by the
sides $s_{m+1},\cdots,s_n$ from $Q$ and $s_{m+1}',\cdots,s_n'$ from $Q'$ we find
that $\langle {\bf v}_1', {\bf v}_m' \rangle \leq \langle {\bf v}_1, {\bf v}_m
\rangle$.  In both cases we obtain a contradiction.

\vspace{.05in}
We have not used any restriction on the dihedral angles, only the restriction
that $P_1$ and $P_2$ are have parallelogram faces, so we shown that $\alpha :
{\cal P}^1_C \rightarrow \mathbb{R}^E$ is injective.  \Endproof Proposition
\ref{INJECTIVE}.

Because ${\cal P}^0_C \subset {\cal P}^1_C$, it follows immediately that:
\begin{cor}
$\alpha : {\cal P}^0_C \rightarrow A_C$ is injective.
\end{cor}

\noindent
This gives the uniqueness part of Andreev's Theorem.

\section{The mapping $\alpha: {\cal P}_C^0 \rightarrow A_C$ is proper.} \label{SEC_PROP}

In this section, we prove that the mapping $\alpha: {\cal P}_C^0
\rightarrow A_C$ is a proper map.  In fact, we will prove a more
general statement (Proposition \ref{GENERAL}) which will be useful later in the paper.

\begin{lem} \label{FACEANG}
Let $\cal{F}$ be a face of a hyperbolic polyhedron $P$ with non-obtuse
dihedral angles. If a face angle of $\cal{F}$ equals $\pi/2$ at the
vertex $v$, then the dihedral angle of the edge opposite the face angle
 (the edge that enters $v$ and is not in $F$) is $\pi/2$ and the
dihedral angle of one of the two edges in $\cal{F}$ that enters $v$ is
$\pi/2$.
\end{lem}

\noindent {\bf Proof:}
This will follow from Equation (\ref{TLC}) in Proposition \ref{TRIVALENT}, which one
can use to calculate face angles from the dihedral angles at a vertex.  In
Equation (\ref{TLC}), if $\beta_i = \pi/2$ we have:

$$0 = \frac{\cos(\alpha_i)
+\cos(\alpha_j)\cos(\alpha_k)} {\sin(\alpha_j)\sin(\alpha_k)},$$

\noindent where $\alpha_i$ is the dihedral angle opposite the face
angle $\beta_i$ and $\alpha_j,\alpha_k$ are the dihedral angles of the
other two edges entering $v$.  Both $\cos(\alpha_i) \geq 0$ and
$\cos(\alpha_j)\cos(\alpha_k)  \geq 0$ for non-obtuse
$\alpha_i,\alpha_j,$ and $\alpha_k$, so that $\cos(\alpha_i) = 0$ and
$\cos(\alpha_j)\cos(\alpha_k)=0$.  Hence $\alpha_i = \pi/2$ and either
$\alpha_j = \pi/2$ or $\alpha_k = \pi/2$.  \Endproof

\begin{lem} \label{NORMALIZE}
Given three points $v_1,v_2,v_3$ that form a non-obtuse, non-degenerate
triangle in the Poincar\'e model of $\mathbb{H}^3$, there is a unique
isometry taking $v_1$ to a positive point on the $x$-axis, $v_2$ to a
positive point on the $y$-axis, and $v_3$ to a positive point on the
$z$-axis.
\end{lem}

\noindent {\bf Proof:}
The points $v_1,v_2,$ and $v_3$ form a triangle $T$ in a plane $P_T$.
It is sufficient to show that there is a plane $Q_T$ in the Poincar\'e
ball model that intersects the positive octant in a triangle isomorphic
to $T$.  The isomorphism taking $v_1,v_2,$ and $v_3$ to the $x,y$, and
$z$-axes will then be the one that takes the plane $P_T$ to the plane
$Q_T$ and the triangle $T$ to the intersection of $Q_T$ with the
positive octant.

Let $s_1,s_2,$ and $s_3$ be the side lengths of $T$.  The plane $Q_T$
must intersect the $x,y,$ and $z$-axes distances $a_1,a_2,$ and $a_3$
satisfying the hyperbolic Pythagorean theorem:

$$\cosh(s_1) = \cosh(a_2) \cosh(a_3)$$
$$\cosh(s_2) = \cosh(a_3) \cosh(a_1)$$
$$\cosh(s_3) = \cosh(a_1) \cosh(a_2)$$

These equations can be solved for $\left(\cosh^2(a_1),\cosh^2(a_2),\cosh^2(a_3)\right)$,
obtaining
$$\left( \frac{\cosh(s_2) \cosh(s_3)}{\cosh(s_1)},
\frac{\cosh(s_3) \cosh(s_1)}{\cosh(s_2)}, \frac{\cosh(s_1)
\cosh(s_2)}{\cosh(s_3)} \right),$$
The only concern in solving for $a_i$ is that each of these expressions
is $\geq 1$. However, this follows from the triangle $T$ being
non-obtuse, using the hyperbolic law of cosines. \Endproof

\vspace{.07in}
All of the results in this chapter are corollaries to the
following proposition:

\begin{prop}\label{GENERAL}
Given a sequence of compact polyhedra $P_i$ realizing $C$ with
$\alpha(P_i) = {\bf a}_i \in A_C$.  If ${\bf a}_i$ converges to ${\bf
a} \in \overline{A_C}$, satisfying conditions (1,3-5), then there
exists a polyhedron $P_0$ realizing $C$ with
dihedral angles ${\bf a}$.
\end{prop}

\noindent {\bf Proof: }
Throughout this proof, we will denote the vertices of $P_i$ by
$v_1^i,\cdots,v_n^i$.  Let $v_a^i,v_b^i$, and $v_c^i$ be three vertices
on the same face of $P_i$, which will form a triangle with non-obtuse
angles for all $i$.  According to Lemma \ref{NORMALIZE}, for each $i$
we can normalize $P_i$ in the Poincar\'e ball model so that $v_a^i$ is
on the $x$-axis, $v_b^i$ is on the $y$-axis, and $v_c^i$ is on the
$z$-axis.

For each $i$, the vertices $v^i_1,\cdots,v^i_n$ are in $\overline{
\mathbb{H}^3}$, which is a compact space under the Euclidean metric.
Therefore, by taking a subsequence, if necessary, we can assume that
each vertex of $P_i$ converges to some point in
$\overline{\mathbb{H}^3}$. We denote the collection of all of these
limit points of $v^i_1,\cdots,v^i_n$ by ${\cal A}_1,\cdots,{\cal A}_q
\in \overline{\mathbb{H}^3}$.  Let $P_0$ be the convex hull of ${\cal
A}_1,\cdots,{\cal A}_q$. Since each $P_i$ realizes $C$, if we can show
that each ${\cal A}_m$ is the limit of a single vertex of $P_i$, then
$P_0$ will realize $C$ and have dihedral angles ${\bf a}$.

In summary, we must show that no more than one vertex converges to each
${\cal A}_m$, using that ${\bf a}$ satisfies conditions (1), (3), (4),
and (5).

\vspace{.07in}\noindent
{\em \large We first check that if ${\cal A}_m \in \partial
\mathbb{H}^3$, then there is a single vertex of $P_i$ converging to
${\cal A}_m$: }
\vspace{.05in}

Therefore, suppose that there are $k>1$ vertices of $P_i$ converging to
${\cal A}_m$ to find a contradiction to the fact that ${\bf a}$
satisfies conditions (1), and (3-5).  Without loss of generality, we
assume that $v_1^i,\cdots,v_k^i$ converge to ${\cal A}_m$ and
$v_{k+1}^i,\cdots,v_n^i$ converge to other points ${\cal A}_j$ for $j
\neq m$.

Since ${\cal A}_m$ is at infinity, our normalization (restricting
$v_a^i,v_b^i,$ and $v_c^i$ to the $x,y,$ and $z$-axes respectively)
ensures that at least two of these vertices ($v_a^i,v_b^i,$ and
$v_c^i$) do not converge to ${\cal A}_m$. This fact will be essential
throughout this part of the proof.

Doing a Euclidean rotation of the entire Poincar\'e ball, we can assume
that ${\cal A}_m$ is at the north pole of the sphere, without changing
the fact that there are at least two vertices of $P_i$ that do not
converge to ${\cal A}_m$.

We will do a sequence of normalizations of the position of $P_i$ in the
Poincar\'e ball model to study the geometry of $P_i$ near ${\cal A}_m$
as $i$ increases.

For all sufficiently large $i$, there is some hyperbolic plane $Q$,
which is both perpendicular to the $z$-axis and has ${\cal A}_m$ and
$v_1^i,\cdots,v_k^i$ on one side of $Q$ and the remaining ${\cal A}_j$
($j \neq m$)  and all of the vertices $v_{k+1}^i,\cdots,v_n^i$ that do
not converge to ${\cal A}_m$ on the other side of $Q$. This is possible
because ${\cal A}_1,\cdots,{\cal A}_q$ are distinct points in
$\overline{\mathbb{H}^3}$, and because ${\cal A}_m \in \partial
\mathbb{H}^3$.

For each $i$, let $R_i$ be the hyperbolic plane which intersects the
$z$-axis perpendicularly, and at the point farthest from the origin
such that the closed half-space toward ${\cal A}_m$ contains all
vertices which will converge to ${\cal A}_m$ (see figure below).  Let
$D_i$ denote the distance from $R_i$ to $Q$ along the $z$-axis; as $i
\rightarrow \infty$ and the vertices $v^i_1,\cdots,v^i_k$ tend to
${\cal A}_m \in \partial \mathbb{H}^3$, $D_i \rightarrow \infty$.  Let
$S_i$ be the hyperbolic plane intersecting the $z$-axis
perpendicularly, halfway between $R_i$ and $Q$.

\vspace{.07in}
\begin{center}
\begin{picture}(0,0)%
\includegraphics{./normalization.pstex}%
\end{picture}%
\setlength{\unitlength}{4144sp}%
\begingroup\makeatletter\ifx\SetFigFont\undefined%
\gdef\SetFigFont#1#2#3#4#5{%
  \reset@font\fontsize{#1}{#2pt}%
  \fontfamily{#3}\fontseries{#4}\fontshape{#5}%
  \selectfont}%
\fi\endgroup%
\begin{picture}(4577,2060)(1327,-2093)
\put(4403,-2069){\makebox(0,0)[lb]{\smash{{\SetFigFont{6}{7.2}{\familydefault}{\mddefault}{\updefault}{\color[rgb]{0,0,0}$Q$}%
}}}}
\put(2243,-90){\makebox(0,0)[lb]{\smash{{\SetFigFont{6}{7.2}{\familydefault}{\mddefault}{\updefault}{\color[rgb]{0,0,0}${\cal A}_m$}%
}}}}
\put(2552,-480){\makebox(0,0)[lb]{\smash{{\SetFigFont{6}{7.2}{\familydefault}{\mddefault}{\updefault}{\color[rgb]{0,0,0}$v_{k+1}^i$}%
}}}}
\put(2039,-1691){\makebox(0,0)[lb]{\smash{{\SetFigFont{6}{7.2}{\familydefault}{\mddefault}{\updefault}{\color[rgb]{0,0,0}$v_{k+4}^i$}%
}}}}
\put(2054,-113){\makebox(0,0)[lb]{\smash{{\SetFigFont{6}{7.2}{\familydefault}{\mddefault}{\updefault}{\color[rgb]{0,0,0}$R_i$}%
}}}}
\put(1907,-140){\makebox(0,0)[lb]{\smash{{\SetFigFont{6}{7.2}{\familydefault}{\mddefault}{\updefault}{\color[rgb]{0,0,0}$S_i$}%
}}}}
\put(1763,-205){\makebox(0,0)[lb]{\smash{{\SetFigFont{6}{7.2}{\familydefault}{\mddefault}{\updefault}{\color[rgb]{0,0,0}$Q$}%
}}}}
\put(1529,-1032){\makebox(0,0)[lb]{\smash{{\SetFigFont{6}{7.2}{\familydefault}{\mddefault}{\updefault}{\color[rgb]{0,0,0}$v_{k+5}^i$}%
}}}}
\put(2482,-942){\makebox(0,0)[lb]{\smash{{\SetFigFont{6}{7.2}{\familydefault}{\mddefault}{\updefault}{\color[rgb]{0,0,0}$v_{k+2}^i$}%
}}}}
\put(2635,-1244){\makebox(0,0)[lb]{\smash{{\SetFigFont{6}{7.2}{\familydefault}{\mddefault}{\updefault}{\color[rgb]{0,0,0}$v_{k+3}^i$}%
}}}}
\put(1917,-598){\makebox(0,0)[lb]{\smash{{\SetFigFont{6}{7.2}{\familydefault}{\mddefault}{\updefault}{\color[rgb]{0,0,0}$v_n^i$}%
}}}}
\put(4871,-94){\makebox(0,0)[lb]{\smash{{\SetFigFont{6}{7.2}{\familydefault}{\mddefault}{\updefault}{\color[rgb]{0,0,0}${\cal A}_m$}%
}}}}
\put(4101,-1104){\makebox(0,0)[lb]{\smash{{\SetFigFont{6}{7.2}{\familydefault}{\mddefault}{\updefault}{\color[rgb]{0,0,0}$S_i$}%
}}}}
\put(4423,-167){\makebox(0,0)[lb]{\smash{{\SetFigFont{6}{7.2}{\familydefault}{\mddefault}{\updefault}{\color[rgb]{0,0,0}$R_i$}%
}}}}
\end{picture}%

\end{center}
\vspace{.07in}

For each $i$, we normalize the polyhedra $P_i$ by translating the plane
$S_i$ along the $z$-axis to the equatorial plane, $H = \{z=0\}$.  We consider
this a change of viewpoint, so the translated points and planes retain
their former names. Hence, under this normalization, we have planes
$R_i$ and $Q$ both perpendicular to the $z$-axis, and hyperbolic
distance $D_i/2$ from the origin.  Each vertex $v_1^i,\cdots,v_k^i$ is
bounded above the plane $R_i$, and all of the vertices
$v_{k+1}^i,\cdots,v_n^i$ of $P_i$ that do not converge to ${\cal A}_m$
are bounded below the plane $Q$.

As $i$ tends to infinity, $R_i$ and $Q$ intersect the $z$-axis at
arbitrarily large hyperbolic distances from the origin.  Denote by
${\cal N}$ and by ${\cal S}$ the half-spaces that these planes bound
away from the origin. Given arbitrarily small (Euclidean) neighborhoods
of the pole, for all sufficiently large $i$ the half-spaces will be
contained in these neighborhoods. Hence, any edge running from ${\cal
N}$ to ${\cal S}$ will intersect $H$ almost
orthogonally, and close to the origin, as illustrated in the figure
below.  Let $e^i_1,\cdots,e^i_l$ denote the collection of such edges.

\vspace{.07in}
\begin{center}
\begin{picture}(0,0)%
\includegraphics{./non_compact.pstex}%
\end{picture}%
\setlength{\unitlength}{4144sp}%
\begingroup\makeatletter\ifx\SetFigFont\undefined%
\gdef\SetFigFont#1#2#3#4#5{%
  \reset@font\fontsize{#1}{#2pt}%
  \fontfamily{#3}\fontseries{#4}\fontshape{#5}%
  \selectfont}%
\fi\endgroup%
\begin{picture}(2392,2390)(580,-2129)
\put(1015,-861){\makebox(0,0)[lb]{\smash{{\SetFigFont{8}{9.6}{\familydefault}{\mddefault}{\updefault}{\color[rgb]{0,0,0}$P_i \cap H$}%
}}}}
\put(2299,-956){\makebox(0,0)[lb]{\smash{{\SetFigFont{8}{9.6}{\familydefault}{\mddefault}{\updefault}{\color[rgb]{0,0,0}$H$}%
}}}}
\put(1725,-2041){\makebox(0,0)[lb]{\smash{{\SetFigFont{8}{9.6}{\familydefault}{\mddefault}{\updefault}{\color[rgb]{0,0,0}$\cal{S}$}%
}}}}
\put(1718, 91){\makebox(0,0)[lb]{\smash{{\SetFigFont{8}{9.6}{\familydefault}{\mddefault}{\updefault}{\color[rgb]{0,0,0}$\cal{N}$}%
}}}}
\end{picture}%

\end{center}
\vspace{.07in}

By assumption, there are $k \geq 2$ vertices in ${\cal N}$ and due to
our normalization, there are two or more vertices in ${\cal S}$.

Thus the intersection $P_i \cap H$ will be almost a Euclidean polygon
and its angles will be almost the dihedral angles
$\alpha(e_1^i),...\alpha(e_l^i)$; in particular, for $i$ sufficiently
large they will be at most only slightly larger than $\pi/2$.  This
implies that $l_i$, the number of such edges, is 3 or 4 for $i$
sufficiently large, as a Euclidean polygon with at least 5 faces has at
least one angle $\geq 3\pi/5$.

If $l=3$, the edges $e_1^i,e_2^i,e_3^i$ intersecting $H$ are a
prismatic circuit because there are two or more vertices in both ${\cal
N}$ and in ${\cal S}$.  The sum
$\alpha(e_1^i)+\alpha(e_2^i)+\alpha(e_3^i)$ tends to $\pi$ as $i$ tends
to infinity, hence ${\bf a}$ cannot satisfy condition (3), contrary to
assumption.

Similarly, if $l_i = 4$, and if $e_1^i,e_2^i,e_3^i,e_4^i$ to form a
prismatic 4-circuit, then the corresponding sum of dihedral angles
tends to $2\pi$ violating condition (4).

So we are left with the possibility that $l = 4$, and that
$e_1^i,e_2^i,e_3^i,e_4^i$ do not form a prismatic 4-circuit.  In this
case, a pair of these edges meet at a vertex which may be in either
${\cal N}$, as shown in the diagram below, or in ${\cal S}$.  Without loss of 
generality, we assume that $e_1^i$ and $e_2^i$ meet at this vertex,
which we call $x^i$ and we assume that $x^i \in {\cal N}$.  
Because we assume that there are two or more
vertices converging to ${\cal A}_m$, there most be some edge $e_j^i$
($j \neq 1,2,3,4$) meeting $e_1^i$ and $e_2^i$ at $x^i$.  We denote by
$f^i$ the face of $P_i$ containing $e_1^i$, $e_2^i$ and $x^i$. An
example of this situation is drawn in the diagram below, although the
general situation can be more complicated.

\begin{center}
\begin{picture}(0,0)%
\includegraphics{./non_compact3.pstex}%
\end{picture}%
\setlength{\unitlength}{4144sp}%
\begingroup\makeatletter\ifx\SetFigFont\undefined%
\gdef\SetFigFont#1#2#3#4#5{%
  \reset@font\fontsize{#1}{#2pt}%
  \fontfamily{#3}\fontseries{#4}\fontshape{#5}%
  \selectfont}%
\fi\endgroup%
\begin{picture}(2388,2386)(580,-2126)
\put(1289,-767){\makebox(0,0)[lb]{\smash{{\SetFigFont{8}{9.6}{\familydefault}{\mddefault}{\updefault}{\color[rgb]{0,0,0}$f^i$}%
}}}}
\put(1718,-2010){\makebox(0,0)[lb]{\smash{{\SetFigFont{8}{9.6}{\familydefault}{\mddefault}{\updefault}{\color[rgb]{0,0,0}$\cal{S}$}%
}}}}
\put(2209,-990){\makebox(0,0)[lb]{\smash{{\SetFigFont{8}{9.6}{\familydefault}{\mddefault}{\updefault}{\color[rgb]{0,0,0}$H$}%
}}}}
\put(1683,-345){\makebox(0,0)[lb]{\smash{{\SetFigFont{8}{9.6}{\familydefault}{\mddefault}{\updefault}{\color[rgb]{0,0,0}$e_1^i$}%
}}}}
\put(1451,-299){\makebox(0,0)[lb]{\smash{{\SetFigFont{8}{9.6}{\familydefault}{\mddefault}{\updefault}{\color[rgb]{0,0,0}$e_2^i$}%
}}}}
\put(1674,-51){\makebox(0,0)[lb]{\smash{{\SetFigFont{8}{9.6}{\familydefault}{\mddefault}{\updefault}{\color[rgb]{0,0,0}$e_j^i$}%
}}}}
\put(1536,101){\makebox(0,0)[lb]{\smash{{\SetFigFont{8}{9.6}{\familydefault}{\mddefault}{\updefault}{\color[rgb]{0,0,0}$x^i$}%
}}}}
\end{picture}%

\end{center}

Since the sum of the dihedral angles along this 4-circuit limits to
$2\pi$, and each dihedral angle is non-obtuse, each of the dihedral
angles $\alpha(e_1^i), \cdots,\alpha(e_4^i)$ limits to $\pi/2$. One can
use Equation (\ref{TLC})  to check that the dihedral angle $\alpha(e_j^i)$
will converge to the face angle $\beta^i$ in the face $f^i$ at vertex
$v^i$.  This is because the right-hand side of this equation limits to
$\cos(\alpha(e_j^i))$, while the right hand side of the equation equals
$\cos(\beta_i)$. Then, as $i$ goes to infinity, the neighborhood ${\cal
N}$ (containing $x^i_1$) converges to the north pole while the
neighborhood ${\cal S}$ (containing the other two vertices in $f^i$
which form the face angle $\beta_i$ at $x^i_1$) converges to the south
pole.  This forces the face angle $\beta^i$ to tend to zero, and hence
the dihedral angle $\alpha(e_j^i)$ must tend to zero as well,
contradicting the fact that all coordinates of the limit point ${\bf
a}$ are positive, by condition (1).

Therefore, we can conclude that any ${\cal A}_m \in \partial
\mathbb{H}^3$ is the limit of a single vertex of the $P_i$.  Hence, the
vertices of $P_i$ that converge to points at infinity (in the original
normalization) converge to distinct points at infinity.

\vspace{.07in}\noindent
{\em \large It remains to show that if ${\cal A}_m \in \mathbb{H}^3$,
then there is a single vertex of $P_i$ converging to ${\cal A}_m$: }
\vspace{.05in}

First, we check that none of the faces of the $P_i$ can degenerate to
either a point, a line segment, or a ray.

Because we have already proven that each of the vertices of $P_i$ that
converge to points in $\partial \mathbb{H}^3$ converge to distinct
points, any face ${\cal F}_i$ of $P_i$ that degenerates to a point, a line segment,
 or a ray, must have area that limits to zero.
Hence, by the Gauss-Bonnet formula, the sum of the face angles would
have to converge to $\pi(n-2)$, if the face has $n$ sides.  The fact
that the face angles are non-obtuse, allows us to see that $\pi(n-2)
\leq n\pi/2$ implying that $n \leq 4$. This restricts such a
degenerating face ${\cal F}_i$ to
either a triangle or a quadrilateral, the only Euclidean polygons
having non-obtuse angles.

If ${\cal F}_i$ is a triangle, then the three edges leading to ${\cal F}_i$
form a prismatic 3-circuit, because our hypothesis $N>4$ implies that
$C$ is not the simplex. If ${\cal F}_i$ degenerates to a point, the three
faces adjacent to ${\cal F}_i$ would meet at a vertex, in the
limit.  
Therefore, by Lemma 3.2, the sum of the dihedral angles at the
edges leading to ${\cal F}_i$ would limit to a value $\geq \pi$, contrary to
condition (3). Otherwise, if ${\cal F}_i$ is a triangle which degenerates
to a line segment or ray in the limit, then two of its face angles will tend
to $\pi/2$.  Then, by Lemma \ref{FACEANG}, the dihedral angles at the
edges opposite of these face angles become $\pi/2$.  However, these
edges are part of the prismatic 3-circuit of edges leading to
${\cal F}_i$, resulting in an angle sum $\geq \pi$, contrary to condition
(3).

If, on the other hand, ${\cal F}_i$ is a quadrilateral, each of the face
angles would have to limit to $\pi/2$.  By Lemma \ref{FACEANG}, the
dihedral angles at each of the edges leading from ${\cal F}_i$ to the rest
of $P_i$ would limit to $\pi/2$, as well as at least one edge of
${\cal F}_i$ leading to each vertex of ${\cal F}_i$.  Therefore, the dihedral
angles at each of the edges leading from ${\cal F}_i$
to the rest of $P_i$
and at least one opposite pair of edges of ${\cal F}_i$ limit to $\pi/2$,
in violation of condition (5).

Since none of the faces of the $P_i$ can degenerate to a point a
line segment, or a ray, neither can the $P_i$.  Suppose that the $P_i$ degenerate
to a polygon, $\cal{G}$.  Because the dihedral angles are non-obtuse,
only two of the faces of the $P_i$ can limit to the polygon $\cal{G}$.
Therefore the rest of the faces of the $P_i$ must limit to points,
line segments, or rays, contrary to our reasoning above.

We are now ready to show that any ${\cal A}_m$ that is a point in
$\mathbb{H}^3$ is the limit of a single vertex of the $P_i$. Let
$v_1,\cdots,v_k$ be the distinct vertices that converge to the same
point ${\cal A}_m$.  Then, since the $P_i$ do not shrink to a point, a
line segment, a ray, or a polygon, there are at least three vertices $v_p^i,
v_q^i$ and $v_r^i$ that don't converge to ${\cal A}_m$ and that don't
converge to each other. Perform the appropriate isometry taking ${\cal
A}_m$ to the origin in the ball model.  Place sphere $S$ centered at
the origin, so small that $v_p^i, v_q^i$ and $v_r^i$ never enter $S$.

For all sufficiently large values of $i$, the intersection $P_i \cap S$
approximates a spherical polygon whose angles approximate the dihedral
angles between the faces of $P_i$ that enter $S$.  These spherical
polygons cannot degenerate to a point or a line segment because the
polyhedra $P_i$ do not degenerate to a line segment, a ray, or a polygon.  By
reasoning similar to that of Proposition 1.1, one can check that this
polygon must have only three sides and angle sum $> \pi$.  If $k>1$,
the edges of this triangle form a prismatic 3-circuit in $C^*$, since
for each $i$, $P_i$ has more than one vertex inside the sphere ($k >
1$) and at least the three vertices corresponding to the points
represented by $v_p^i, v_q^i$ and $v_r^i$. So, the $P_i$ would have a
prismatic 3-circuit whose angle sum limits to a value $> \pi$.
However, this contradicts our assumption that ${\bf a}$ satisfies condition (3).
Therefore, we must have $k=1$.

Therefore, we conclude that each ${\cal A}_m$ is the limit of a single
vertex of $P_i$.
\Endproof

\begin{cor} \label{PROPER}
The mapping $\alpha : {\cal P}_C^0 \rightarrow A_C$ is proper.
\end{cor}
\noindent {\bf Proof:}
Suppose that $P_i$ is a sequence of polyhedra in ${\cal P}_C^0$ with
$\alpha(P_i) = {\bf a_i} \in A_C$. If the sequence ${\bf a_i}$
converges to a point in ${\bf a} \in A_C$, we must show that there is a
subsequence of the $P_i$ that converges in ${\cal P}_C^0$.

Since ${\bf a}$ satisfies conditions (1) and (3-5), by Proposition
\ref{GENERAL}, there is a subsequence of the $P_i$ converging to a
polyhedron $P_0$ that realizes $C$ with dihedral angles ${\bf a}$. The
sum of the dihedral angles at each vertex of $P_0$ is $> \pi$ since
${\bf a}$ satisfies condition (2) as well.  Therefore, by Lemma 3.2,
each vertex of $P_0$ is at a finite point in $\mathbb{H}^3$, giving
that $P_0$ is compact. Hence $P_0$ realizes $C$, is compact, and has
non-obtuse dihedral angles. Therefore $P_0 \in {\cal P}_C^0$.
\Endproof

\begin{cor}\label{INFVERTEX}
Suppose that Andreev's Theorem holds for $C$ and suppose that there
is a sequence ${\bf a}_i \in A_C$ converging to ${\bf a} \in \partial
A_C$. If ${\bf a}$ satisfies conditions (1,3-5) and if condition (2) is
satisfied for vertices $v_1,\cdots,v_k$ of $C$, but not for
$v_{k+1},\cdots,v_n$ for which the dihedral angle sum is exactly $\pi$,
then there exists a non-compact polyhedron $P_0$ realizing $C$ with
dihedral angles ${\bf a}$. $P_0$ has vertices $v_1,\cdots,v_k$ at
distinct finite points and the $v_{k+1},\cdots,v_n$ at distinct points
at infinity.
\end{cor}

\noindent {\bf Proof: }
Because we assume that Andreev's Theorem holds for $C$, there exists
a sequence of polyhedra $P_i$ with $\alpha(P_i) = {\bf a}_i$. The proof
continues in the same way as that of Corollary \ref{PROPER}, except
that in this case, it then follows from Lemma 3.2 that $v_1,\cdots,v_k$
lie in $\mathbb{H}^3$, while $v_{k+1},\cdots,v_n$ lie in $\partial
\mathbb{H}^3.$
\Endproof

Notice that if an abstract polyhedron has no prismatic 3-circuits, it
has no triangular faces, so collapsing a single edge to a point
results in a cell complex satisfying all of the conditions of an
abstract polyhedron, except that one of the vertices is now 4-valent.

\begin{prop}\label{PREWHITEHEAD}
Let $C$ be an abstract polyhedron having no prismatic 3-circuits for
which Andreev's Theorem is satisfied.  For any edge $e_0$ of $C$, let
$C_0$ be the complex obtained by contracting $e_0$ to a point.  Then,
there exists a non-compact polyhedron $P_0$ realizing $C_0$ with the
edge $e_0$ contracted to a single vertex at infinity and the rest of
the vertices at finite points in $\mathbb{H}^3$.
\end{prop}

\noindent {\bf Proof of Proposition \ref{PREWHITEHEAD}:}
Let $v_1$ and $v_2$ be the vertices at the ends of $e_0$, let
$e_1,e_2,e_3,e_4$ the edges emanating from the ends of $e_0$, and
$f_1,f_2,f_3,f_4$ be the four faces containing either $v_1$ or $v_2$,
as illustrated below.

\vspace{.07in}
\begin{center}
\begin{picture}(0,0)%
\includegraphics{./whitehead_config.pstex}%
\end{picture}%
\setlength{\unitlength}{4144sp}%
\begingroup\makeatletter\ifx\SetFigFont\undefined%
\gdef\SetFigFont#1#2#3#4#5{%
  \reset@font\fontsize{#1}{#2pt}%
  \fontfamily{#3}\fontseries{#4}\fontshape{#5}%
  \selectfont}%
\fi\endgroup%
\begin{picture}(2023,584)(1,-152)
\put(  1,106){\makebox(0,0)[lb]{\smash{{\SetFigFont{7}{8.4}{\familydefault}{\mddefault}{\updefault}{\color[rgb]{0,0,0}$f_3$}%
}}}}
\put(852,191){\makebox(0,0)[lb]{\smash{{\SetFigFont{7}{8.4}{\familydefault}{\mddefault}{\updefault}{\color[rgb]{0,0,0}$e_0$}%
}}}}
\put(382,109){\makebox(0,0)[lb]{\smash{{\SetFigFont{7}{8.4}{\familydefault}{\mddefault}{\updefault}{\color[rgb]{0,0,0}$v_1$}%
}}}}
\put(1378,109){\makebox(0,0)[lb]{\smash{{\SetFigFont{7}{8.4}{\familydefault}{\mddefault}{\updefault}{\color[rgb]{0,0,0}$v_2$}%
}}}}
\put(1511,247){\makebox(0,0)[lb]{\smash{{\SetFigFont{7}{8.4}{\familydefault}{\mddefault}{\updefault}{\color[rgb]{0,0,0}$e_1$}%
}}}}
\put(207,274){\makebox(0,0)[lb]{\smash{{\SetFigFont{7}{8.4}{\familydefault}{\mddefault}{\updefault}{\color[rgb]{0,0,0}$e_4$}%
}}}}
\put(202,-39){\makebox(0,0)[lb]{\smash{{\SetFigFont{7}{8.4}{\familydefault}{\mddefault}{\updefault}{\color[rgb]{0,0,0}$e_3$}%
}}}}
\put(1796,104){\makebox(0,0)[lb]{\smash{{\SetFigFont{7}{8.4}{\familydefault}{\mddefault}{\updefault}{\color[rgb]{0,0,0}$f_1$}%
}}}}
\put(852,-113){\makebox(0,0)[lb]{\smash{{\SetFigFont{7}{8.4}{\familydefault}{\mddefault}{\updefault}{\color[rgb]{0,0,0}$f_2$}%
}}}}
\put(852,366){\makebox(0,0)[lb]{\smash{{\SetFigFont{7}{8.4}{\familydefault}{\mddefault}{\updefault}{\color[rgb]{0,0,0}$f_4$}%
}}}}
\put(1519,-52){\makebox(0,0)[lb]{\smash{{\SetFigFont{7}{8.4}{\familydefault}{\mddefault}{\updefault}{\color[rgb]{0,0,0}$e_2$}%
}}}}
\end{picture}%

\end{center}

For $\epsilon \in (0,\pi/2]$
the angles:  $\alpha(e_0) = \epsilon,$ $\alpha(e_1) = \alpha(e_2) =
\alpha(e_3) = \alpha(e_4)  = \pi/2,$ and $\alpha(e)  = 2\pi/5$, for all
other edges $e$, are in $A_C$ since $C$ has no prismatic 3-circuits. Therefore,
because we assume that Andreev's Theorem holds for $C$, there is a polyhedron
$P_{\epsilon} \in {\cal P}_C^0$ realizing these angles.  Choose a sequence
$\epsilon_i > 0$ converging to $0$.

As in the proof of Proposition \ref{GENERAL}, we choose three vertices
of $P_i$ that are on the same face (but not on $f_2$ or $f_4$) and
normalize the polyhedra so that each $P_i$ has these three vertices on
the $x$-axis, the $y$-axis, and the $z$-axis, respectively. The
vertices $v_1^i,\cdots,v_n^i$ of each $P_i$ are in the compact space
$\overline{\mathbb{H}^3}$ so, as before, we can take a subsequence so
that each vertex converges to some point in $\overline{\mathbb{H}^3}$.

The strategy of the proof is to consider the limit points ${\cal
A}_1,\cdots,{\cal A}_q$ of these vertices and to show that $v_1^i$ and
$v_2^i$ converge to the same limit point, say ${\cal A}_1$ and that
exactly one of each of the other vertices $v_j^i$ (for $j>2$) converges
to each of the other ${\cal A}_m$, with $m>1$.

We first check that that $v_1^i$ and $v_2^i$ converge to the same point
at infinity. Since the dihedral angle at edge $e_0$ converges to $0$,
the two planes carrying faces $f_2$ and $f_4$ will intersect at
dihedral angles converging to $0$.  In the limit, these planes will
intersect with dihedral angle $0$, therefore, at a single point in
$\partial {\mathbb{H}^3}$.  The edge $e_0$ is contained within the
intersection of these two planes, hence the entire edge $e_0$ converges
to a single point in $\partial \mathbb{H}^3$, and hence $v_1^i$ and
$v_2^i$ converge to this point. We label this point by ${\cal A}_1$.

We must now show $v_1^i$ and $v_2^i$ are the only vertices that
converge to ${\cal A}_1$. As in the proof of Proposition \ref{GENERAL},
one can do further normalizations so that all of the vertices that
converge to ${\cal A}_1$ are in an arbitrarily small neighborhood
${\cal N}$ of the north pole and all of the other vertices are in an
arbitrarily small neighborhood ${\cal S}$ of the south pole.  Since
${\bf a}$ satisfies conditions (3) and (4) we deduce that
there there are exactly four edges $e_1^i,e_2^i,e_3^i$, and $e_4^i$
that connect ${\cal N}$ to ${\cal
S}$ and that do not form a prismatic 4-circuit. Then, by Lemma \ref{4CIRC}, this 4-circuit separates exactly
two vertices from the remaining vertices of $P_i$, hence only $v_1^i$
and $v_2^i$ converge to ${\cal A}_1$.

The proof that a single vertex of $P_i$ converges to each ${\cal A}_m$,
for $m>1$ is almost the same as the proof of Proposition
\ref{GENERAL}, because the dihedral angles ${\bf a}_i$ are non-zero
for all edges other than $e_0$ and because conditions (3-5) are
satisfied.

The only difference is that one must directly
check that faces $f_2$ and $f_4$ do not collapse to line segments or rays, even
though on each of these faces, two vertices converge to the same point
at infinity.  A vertex of $f_2$ or $f_4$ that is not $v_1$ or $v_2$
must have face angle that is strictly acute.  Otherwise, Lemma \ref{FACEANG} would
give the dihedral angles at two of the edges meeting at this vertex
is $\pi/2$, contrary to the fact that at least two
of these dihedral angles are assigned to be $2\pi/5$.  From here,
the Gauss-Bonnet Theorem can be used to check that neither $f_2$ nor
$f_4$ collapses to a line segment or a ray.

\vspace{.05in}
Therefore, we conclude that the vertices $v_j^i$ ($i>3$) converge to
distinct points in $\overline{\mathbb{H}^3}$ away from the limit of
$v_1^i$ and $v_2^i$, which converge to the same point in $\partial
\mathbb{H}^3$.
Since condition (2) is satisfied at each vertex $v_j^i$
($i>3$), Lemma 3.2 guarantees that each of these ${\cal A}_m$ is at a
finite point.
Therefore, $v_1^i$ and $v_2^i$ converge to the same point at infinity,
and each of the other vertices converges to a distinct finite point.
\Endproof

\section{$A_C \neq \emptyset$ implies ${\cal P}_C^0 \neq \emptyset$}\label{SEC_NONEMPT}
At this point, we know the following result:

\begin{prop} \label{NONEMPTYIMPLIESAND}
If ${\cal P}^0_C \neq \emptyset$, then $\alpha : {\cal P}^0_C \rightarrow
A_C$ is a homeomorphism.
\end{prop}

\noindent {\bf Proof:}
We have shown in preceding sections that $\alpha: {\cal P}_C^1 \rightarrow
\mathbb{R}^E$ is a continuous and injective map whose domain and range are
manifolds (without boundary) of the same dimension, so it follows from invariance
of domain that $\alpha$ is a local homeomorphism.  Local homeomorphisms restrict
nicely to subspaces, giving that $\alpha: {\cal P}_C^0 \rightarrow \mathbb{R}^E$ is
a local homeomorphism, as well.  In fact, because $\alpha: {\cal P}_C^0
\rightarrow \mathbb{R}^E$ is also injective it is a homeomorphism onto its image,
which we showed in Section \ref{SEC_INEQSAT} is a subset of $A_C$.

Since $A_C$ is convex and therefore connected, it suffices to show that
$\alpha({\cal P}_C^0)$ is both open and closed in $A_C$, for it will then
follow (since the image is non-empty by hypothesis) that the image is the
entire set $A_C$.  But local homeomorphisms are open maps, so $\alpha({\cal
P}_C^0)$ is open in  $A_C$; and since Corollary \ref{PROPER} shows that
$\alpha: {\cal P}_C^0 \rightarrow A_C$ is proper, it immediately follows that
the limit of any sequence in the image of $\alpha$ which converges in $A_C$
must lie in the image of $\alpha$, so $\alpha({\cal P}_C^0)$ is closed in
$A_C$.  \Endproof

Indeed, any local homeomorphism between metric spaces which is also proper
will be a finite-sheeted covering map \cite[p. 23]{DOUADY} and 
\cite[p. 127]{LIMA}.  This 
gives an alternative route to proving the above result.

\vspace{.07in}

But what is left is absolutely not obvious, and is the hardest part of the
whole proof: proving that if $A_C \neq \emptyset$, then ${\cal P}^0_C \neq
\emptyset$.  We have no tools with which to approach it and must use bare
hands.  We follow the proof of Andreev, although the proof of his key lemma
contains a significant error.  We provide our own correction. 

First recall that in Corollary \ref{EX3APR}, we saw that if $C$ has no
prismatic 3-circuits, $A_C \neq \emptyset$.  We will call polyhedra that have
no prismatic 3-circuits {\it simple polyhedra}.   We will also say that the
dual graph $C^*$ is simple if it satisfies the dual condition, that every
$3$-cycle is the boundary of a single face.  (This usage of ``simple'' follows
Andreev, but differs from that used by others, including Vinberg \cite[p.
47]{VINREFL}, for polyhedra in all dimensions greater than two.)

We first prove that ${\cal P}^0_C \neq \emptyset$ for simple polyhedra, and
hence by Proposition \ref{NONEMPTYIMPLIESAND} that Andreev's Theorem holds for
simple polyhedra.  We then show that for any $C$ having prismatic 3-circuits,
if $A_C \neq \emptyset$, then ${\cal P}^0_C \neq \emptyset$ by making a
polyhedron realizing $C$ from (possibly many) simple polyhedra.  By Proposition
\ref{NONEMPTYIMPLIESAND}, this final step will complete the proof of Andreev's
theorem.

\vspace{.15in}
\noindent
{\large \it{Proof of Andreev's Theorem for Simple Polyhedra}}

\begin{prop} \label{APRISNONEMPTY} If $C$ is simple and has $N > 5$
faces, ${\cal P}^0_C \neq \emptyset$.  In words: every simple polyhedron
is realizable. \end{prop}

\noindent {\bf Proof:} The proof comprises three lemmas.  We will first
state the lemmas and prove this proposition using them.  Then we will
prove the lemmas.  

\begin{lem} \label{EXPRISM} Let $Pr_N$ and $D_N$ be the
abstract polyhedra corresponding to the $N$-faced prism and the $N$-faced
``split prism'', as illustrated below.  If $N > 4$, ${\cal P}^0_{Pr_N}$ is
nonempty and if $N > 7$, ${\cal P}^0_{D_N}$ is nonempty.
\end{lem}

\begin{center} 
\begin{picture}(0,0)%
\includegraphics{prism_comb.pstex}%
\end{picture}%
\setlength{\unitlength}{4144sp}%
\begingroup\makeatletter\ifx\SetFigFont\undefined%
\gdef\SetFigFont#1#2#3#4#5{%
  \reset@font\fontsize{#1}{#2pt}%
  \fontfamily{#3}\fontseries{#4}\fontshape{#5}%
  \selectfont}%
\fi\endgroup%
\begin{picture}(4224,2022)(1114,-1981)
\put(1669,-1968){\makebox(0,0)[lb]{\smash{\SetFigFont{7}{8.4}{\familydefault}{\mddefault}{\updefault}{\color[rgb]{0,0,0}Prism with 10 faces}%
}}}
\put(3769,-1950){\makebox(0,0)[lb]{\smash{\SetFigFont{7}{8.4}{\familydefault}{\mddefault}{\updefault}{\color[rgb]{0,0,0}Splitprism with 11 faces}%
}}}
\end{picture}

\end{center}

A Whitehead move on an edge $e$ of an abstract polyhedron is given by the local
change $Wh(e)$ described by the following diagram.  The Whitehead move in the
dual complex is dashed. (Sometimes we will find it convenient to describe the
Whitehead move entirely in terms of the dual complex, in which case we write
$Wh(f)$).

\vspace{.07in} \begin{center} 
\begin{picture}(0,0)%
\includegraphics{whitehead.pstex}%
\end{picture}%
\setlength{\unitlength}{3947sp}%
\begingroup\makeatletter\ifx\SetFigFont\undefined%
\gdef\SetFigFont#1#2#3#4#5{%
  \reset@font\fontsize{#1}{#2pt}%
  \fontfamily{#3}\fontseries{#4}\fontshape{#5}%
  \selectfont}%
\fi\endgroup%
\begin{picture}(2975,1424)(1339,-1412)
\put(2157,-1381){\makebox(0,0)[lb]{\smash{\SetFigFont{8}{9.6}{\familydefault}{\mddefault}{\updefault}{\color[rgb]{0,0,0}Whitehead move on edge $e$}%
}}}
\put(3613,-72){\makebox(0,0)[lb]{\smash{\SetFigFont{8}{9.6}{\familydefault}{\mddefault}{\updefault}{\color[rgb]{0,0,0}$e_1$}%
}}}
\put(4290,-83){\makebox(0,0)[lb]{\smash{\SetFigFont{8}{9.6}{\familydefault}{\mddefault}{\updefault}{\color[rgb]{0,0,0}$e_2$}%
}}}
\put(4314,-1085){\makebox(0,0)[lb]{\smash{\SetFigFont{8}{9.6}{\familydefault}{\mddefault}{\updefault}{\color[rgb]{0,0,0}$e_3$}%
}}}
\put(3614,-1096){\makebox(0,0)[lb]{\smash{\SetFigFont{8}{9.6}{\familydefault}{\mddefault}{\updefault}{\color[rgb]{0,0,0}$e_4$}%
}}}
\put(4004,-687){\makebox(0,0)[lb]{\smash{\SetFigFont{8}{9.6}{\familydefault}{\mddefault}{\updefault}{\color[rgb]{0,0,0}$e'$}%
}}}
\put(1799,-671){\makebox(0,0)[lb]{\smash{\SetFigFont{8}{9.6}{\familydefault}{\mddefault}{\updefault}{\color[rgb]{0,0,0}$e$}%
}}}
\put(1374,-340){\makebox(0,0)[lb]{\smash{\SetFigFont{8}{9.6}{\familydefault}{\mddefault}{\updefault}{\color[rgb]{0,0,0}$e_1$}%
}}}
\put(2426,-351){\makebox(0,0)[lb]{\smash{\SetFigFont{8}{9.6}{\familydefault}{\mddefault}{\updefault}{\color[rgb]{0,0,0}$e_2$}%
}}}
\put(2483,-1051){\makebox(0,0)[lb]{\smash{\SetFigFont{8}{9.6}{\familydefault}{\mddefault}{\updefault}{\color[rgb]{0,0,0}$e_3$}%
}}}
\put(1357,-1096){\makebox(0,0)[lb]{\smash{\SetFigFont{8}{9.6}{\familydefault}{\mddefault}{\updefault}{\color[rgb]{0,0,0}$e_4$}%
}}}
\put(1992,-556){\makebox(0,0)[lb]{\smash{\SetFigFont{8}{9.6}{\familydefault}{\mddefault}{\updefault}{\color[rgb]{0,0,0}$f$}%
}}}
\put(3793,-497){\makebox(0,0)[lb]{\smash{\SetFigFont{8}{9.6}{\familydefault}{\mddefault}{\updefault}{\color[rgb]{0,0,0}$f'$}%
}}}
\put(2832,-476){\makebox(0,0)[lb]{\smash{\SetFigFont{8}{9.6}{\familydefault}{\mddefault}{\updefault}{\color[rgb]{0,0,0}$Wh(e)$}%
}}}
\end{picture}

\end{center}
\vspace{.07in}

\begin{lem} \label{WHITEHEAD} Let the abstract polyhedron $C'$ be obtained
from the simple abstract polyhedron $C$ by a Whitehead move $Wh(e)$. 
Then if ${\cal P}^0_C$ is non-empty, so is ${\cal P}^0_{C'}.$ \end{lem}

\begin{lem} {\bf (Whitehead Sequence)} \label{ROEDER}
Let $C$ be a simple abstract polyhedron on $\mathbb{S}^2$ which is
not a prism.  If $C$ has $N > 7$ faces, one can simplify $C$ by a finite
sequence of Whitehead moves to $D_N$ such that {\it all of the intermediate
abstract polyhedra are simple.}

\end{lem}

\noindent {\bf Proof of Proposition \ref{APRISNONEMPTY}, assuming these
three lemmas:}
Given simple $C$ with $N > 5$ faces; if $C$ is a prism, the
statement is proven by Lemma \ref{EXPRISM}.  One can check that if $C$ has
$7$ or fewer faces (and is not the tetrahedron) it is a prism.  So, if $C$
is not a prism, we have $N > 7$.  Then, according to Lemma \ref{ROEDER},
one finds a reduction by (say $m$) Whitehead moves to $D_N$, with each
intermediate abstract polyhedron simple. Applying Lemma
\ref{WHITEHEAD} $m$ times, one sees that ${\cal P}^0_C$ is non-empty if
and only if ${\cal P}^0_{D_N}$ is non-empty.  However, ${\cal P}^0_{D_N}$
is non-empty by Lemma \ref{EXPRISM}.  \Endproof

Theorem 6 from Andreev's original paper corresponds to our Proposition
\ref{APRISNONEMPTY}.  The hard technical part of this is the proof of Lemma
\ref{ROEDER}.  Andreev's original proof of Theorem 6 in \cite{AND,AND2}
provides an algorithm giving the Whitehead moves needed for this lemma but the
algorithm {\it just doesn't work}.  It was implemented as a computer program by
the first author and failed on the first test case, $C$ being the dodecahedron.
On one of the final steps, it produced an abstract polyhedron which had a
prismatic 3-circuit.  This error was then traced back to a false statement in
Andreev's proof of the lemma.  We will explain the details of this error in the
proof of Lemma \ref{ROEDER}.

\vspace{0.05in}
\noindent 
{\bf Proof of Lemma \ref{EXPRISM}:}
We construct the $N$-faced prism explicitly.  First, construct a regular
polygon with $N-2$ sides centered at the origin in the disc model for
$\mathbb{H}^2$.  ($N-2 \geq 3$, since $N \geq 5$.)  We can do this with the
angles arbitrarily small.  Now view $\mathbb{H}^2$ as the equatorial plane the
ball model of $\mathbb{H}^3$; and consider the hyperbolic planes which are
perpendicular to the equatorial plane and contain one side of the polygon.  In
Euclidean geometry these are hemispheres with centers on the boundary of the
equatorial disc.  The dihedral angles between intersecting pairs of these
planes are the angles of the polygon. 

Now consider two hyperbolic planes close to the equatorial plane, one slightly
above and one slightly beneath, both perpendicular to the $z$-axis.  These will
intersect the previous planes at angles slightly smaller than $\pi/2$.  The
region defined by these $N$ planes makes a hyperbolic polyhedron realizing the
cell structure of the prism, so ${\cal P}_{Pr_N} \neq \emptyset$.    In
particular, using Proposition \ref{NONEMPTYIMPLIESAND}, Andreev's Theorem holds
for $C = Pr_N$, $N \geq 5$.

%We cut $D_N$ (here $N > 7$) into two prisms and use that Andreev's Theorem holds
%for $Pr_{N-1}$ to construct two prisms with appropriate angles that fit together to
%form $D_N$.  Consider the prism having $N-1$ faces with dihedral angles as
%labeled below.

When $N>7$, the split prism $D_N$ can be constructed
by gluing together a prism and its mirror image, each having
$N-1$ faces.  The dihedral angles are given in the figure
below.

\vspace{.07in} \begin{center} 
\begin{picture}(0,0)%
\epsfig{file=prism.pstex}%
\end{picture}%
\setlength{\unitlength}{4144sp}%
\begingroup\makeatletter\ifx\SetFigFont\undefined%
\gdef\SetFigFont#1#2#3#4#5{%
  \reset@font\fontsize{#1}{#2pt}%
  \fontfamily{#3}\fontseries{#4}\fontshape{#5}%
  \selectfont}%
\fi\endgroup%
\begin{picture}(1824,1711)(830,-1539)
\put(1715,-1036){\makebox(0,0)[lb]{\smash{{\SetFigFont{7}{8.4}{\familydefault}{\mddefault}{\updefault}{\color[rgb]{0,0,0}$\pi/3$}%
}}}}
\put(2654,-710){\makebox(0,0)[lb]{\smash{{\SetFigFont{7}{8.4}{\familydefault}{\mddefault}{\updefault}{\color[rgb]{0,0,0}$\pi/2$}%
}}}}
\put(2445,-1223){\makebox(0,0)[lb]{\smash{{\SetFigFont{7}{8.4}{\familydefault}{\mddefault}{\updefault}{\color[rgb]{0,0,0}$\pi/2$}%
}}}}
\put(1172,-1356){\makebox(0,0)[lb]{\smash{{\SetFigFont{7}{8.4}{\familydefault}{\mddefault}{\updefault}{\color[rgb]{0,0,0}$\pi/2$}%
}}}}
\put(1039,-140){\makebox(0,0)[lb]{\smash{{\SetFigFont{7}{8.4}{\familydefault}{\mddefault}{\updefault}{\color[rgb]{0,0,0}$\pi/2$}%
}}}}
\put(1762, 88){\makebox(0,0)[lb]{\smash{{\SetFigFont{7}{8.4}{\familydefault}{\mddefault}{\updefault}{\color[rgb]{0,0,0}$\pi/2$}%
}}}}
\put(2464,-159){\makebox(0,0)[lb]{\smash{{\SetFigFont{7}{8.4}{\familydefault}{\mddefault}{\updefault}{\color[rgb]{0,0,0}$\pi/2$}%
}}}}
\put(1362,-729){\makebox(0,0)[lb]{\smash{{\SetFigFont{7}{8.4}{\familydefault}{\mddefault}{\updefault}{\color[rgb]{0,0,0}$\pi/3$}%
}}}}
\put(1457,-444){\makebox(0,0)[lb]{\smash{{\SetFigFont{7}{8.4}{\familydefault}{\mddefault}{\updefault}{\color[rgb]{0,0,0}$\pi/3$}%
}}}}
\put(2122,-710){\makebox(0,0)[lb]{\smash{{\SetFigFont{7}{8.4}{\familydefault}{\mddefault}{\updefault}{\color[rgb]{0,0,0}$\pi/3$}%
}}}}
\put(1457,-975){\makebox(0,0)[lb]{\smash{{\SetFigFont{7}{8.4}{\familydefault}{\mddefault}{\updefault}{\color[rgb]{0,0,0}$\pi/3$}%
}}}}
\put(1576,-109){\makebox(0,0)[lb]{\smash{{\SetFigFont{7}{8.4}{\familydefault}{\mddefault}{\updefault}{\color[rgb]{0,0,0}$\pi/2$}%
}}}}
\put(1593,-1256){\makebox(0,0)[lb]{\smash{{\SetFigFont{7}{8.4}{\familydefault}{\mddefault}{\updefault}{\color[rgb]{0,0,0}$\pi/2$}%
}}}}
\put(2122,-1211){\makebox(0,0)[lb]{\smash{{\SetFigFont{7}{8.4}{\familydefault}{\mddefault}{\updefault}{\color[rgb]{0,0,0}$\pi/2$}%
}}}}
\put(1172,-402){\makebox(0,0)[lb]{\smash{{\SetFigFont{7}{8.4}{\familydefault}{\mddefault}{\updefault}{\color[rgb]{0,0,0}$\pi/2$}%
}}}}
\put(830,-710){\makebox(0,0)[lb]{\smash{{\SetFigFont{7}{8.4}{\familydefault}{\mddefault}{\updefault}{\color[rgb]{0,0,0}$\pi/2$}%
}}}}
\put(1187,-1070){\makebox(0,0)[lb]{\smash{{\SetFigFont{7}{8.4}{\familydefault}{\mddefault}{\updefault}{\color[rgb]{0,0,0}$\pi/2$}%
}}}}
\put(2178,-187){\makebox(0,0)[lb]{\smash{{\SetFigFont{7}{8.4}{\familydefault}{\mddefault}{\updefault}{\color[rgb]{0,0,0}$\pi/2$}%
}}}}
\put(1731,-1508){\makebox(0,0)[lb]{\smash{{\SetFigFont{7}{8.4}{\familydefault}{\mddefault}{\updefault}{\color[rgb]{0,0,0}$\pi/4$}%
}}}}
\put(2024,-938){\makebox(0,0)[lb]{\smash{{\SetFigFont{7}{8.4}{\familydefault}{\mddefault}{\updefault}{\color[rgb]{0,0,0}$\pi/3$}%
}}}}
\put(1976,-413){\makebox(0,0)[lb]{\smash{{\SetFigFont{7}{8.4}{\familydefault}{\mddefault}{\updefault}{\color[rgb]{0,0,0}$\pi/3$}%
}}}}
\put(2374,-856){\makebox(0,0)[lb]{\smash{{\SetFigFont{7}{8.4}{\familydefault}{\mddefault}{\updefault}{\color[rgb]{0,0,0}$\pi/2$}%
}}}}
\put(2385,-562){\makebox(0,0)[lb]{\smash{{\SetFigFont{7}{8.4}{\familydefault}{\mddefault}{\updefault}{\color[rgb]{0,0,0}$\pi/2$}%
}}}}
\put(1740,-299){\makebox(0,0)[lb]{\smash{{\SetFigFont{7}{8.4}{\familydefault}{\mddefault}{\updefault}{\color[rgb]{0,0,0}$\pi/3$}%
}}}}
\end{picture}%

\end{center}

These angles satisfy Andreev's conditions (1--5), and Andreev's Theorem holds
for $Pr_{N-1}$ since $N-1>6>5$, so there exists such a hyperbolic prism.  When
this prism is glued to its mirror image, along the $(N-3)$-gon given by the
outermost edges in the figure, the corresponding dihedral angles all double.
So the edges on the outside which were labeled $\pi/2$ ``disappear'' into the
interior of a ``merged'' face, and the edge which was labeled $\pi/4$ now
corresponds to a dihedral angle of $\pi/2$.  Hence, ${\cal P}^0_{D_N} \neq
\emptyset$, when $N>7$.  Notice that when $N = 7$, the construction yields
$Pr_7$ (which is combinatorially equivalent to $D_7$).  \Endproof

%These angles satisfy Andreev's conditions (1-5), so there exists a
%hyperbolic prism with the dihedral angles given in the figure.  When this prism
%is glued to its mirror image, the edges labeled $\pi/2$ on the outside
%disappear as edges, and the edges labeled on the outside by $\pi/4$ glue
%together becoming an edge with dihedral angle $\pi/2.$ Hence, we have
%constructed a polyhedron realizing $D_N$, assuming  $N > 7$.  
%Notice that when $N = 7$, the construction yields $Pr_7$ (which is
%combinatorially equivalent to $D_7$).
%\Endproof

\noindent {\bf Proof of Lemma \ref{WHITEHEAD}:} We are given $C$ and $C'$
simple with $C'$ obtained by a Whitehead move on the edge $e_0$ and we are
given that $P_C \neq \emptyset$.  By Proposition \ref{NONEMPTYIMPLIESAND},
since $P_C \neq \emptyset$, we conclude that Andreev's Theorem holds for $C$.
Let $C_0$ be the complex obtained from $C$ by shrinking the edge $e_0$ down to
a point.  By Proposition \ref{PREWHITEHEAD}, there exists a non-compact
polyhedron $P_0$ realizing $C_0$ since Andreev's Theorem holds for $C$.

We use the upper half-space model of $\mathbb{H}^3$, and normalize so
that $e_0$ has collapsed to the origin of $\mathbb{C} \subset \partial
\mathbb{H}^3$.  The faces incident to $e_0$ are carried by 4 planes $H_1,...,H_4$
each intersecting the adjacent ones at right angles, and all meeting at
the origin.  Their configuration will look like the center of the
following figure.  (Recall that planes in the upper half-space model of
$\mathbb{H}^3$ are hemispheres which intersect $\partial \mathbb{H}^3$ in
their boundary circles.  The dihedral angle between a pair of planes is
the angle between the corresponding pair of circles in $\partial
\mathbb{H}^3$.) 

\vspace{.07in} 
\begin{center} 
\begin{picture}(0,0)%
\epsfig{file=geom_whitehead.pstex}%
\end{picture}%
\setlength{\unitlength}{4144sp}%
\begingroup\makeatletter\ifx\SetFigFont\undefined%
\gdef\SetFigFont#1#2#3#4#5{%
  \reset@font\fontsize{#1}{#2pt}%
  \fontfamily{#3}\fontseries{#4}\fontshape{#5}%
  \selectfont}%
\fi\endgroup%
\begin{picture}(4943,1306)(270,-647)
\put(4906,333){\makebox(0,0)[lb]{\smash{{\SetFigFont{8}{9.6}{\familydefault}{\mddefault}{\updefault}{\color[rgb]{0,0,0}$e_4$}%
}}}}
\put(3151,344){\makebox(0,0)[lb]{\smash{{\SetFigFont{8}{9.6}{\familydefault}{\mddefault}{\updefault}{\color[rgb]{0,0,0}$e_4$}%
}}}}
\put(4771,-556){\makebox(0,0)[lb]{\smash{{\SetFigFont{8}{9.6}{\familydefault}{\mddefault}{\updefault}{\color[rgb]{0,0,0}$e_3$}%
}}}}
\put(342,209){\makebox(0,0)[lb]{\smash{{\SetFigFont{8}{9.6}{\familydefault}{\mddefault}{\updefault}{\color[rgb]{0,0,0}$e_1$}%
}}}}
\put(1053,-517){\makebox(0,0)[lb]{\smash{{\SetFigFont{8}{9.6}{\familydefault}{\mddefault}{\updefault}{\color[rgb]{0,0,0}$e_3$}%
}}}}
\put(544,-489){\makebox(0,0)[lb]{\smash{{\SetFigFont{8}{9.6}{\familydefault}{\mddefault}{\updefault}{\color[rgb]{0,0,0}$e_2$}%
}}}}
\put(816,-275){\makebox(0,0)[lb]{\smash{{\SetFigFont{8}{9.6}{\familydefault}{\mddefault}{\updefault}{\color[rgb]{0,0,0}$e_0$}%
}}}}
\put(2439,-484){\makebox(0,0)[lb]{\smash{{\SetFigFont{8}{9.6}{\familydefault}{\mddefault}{\updefault}{\color[rgb]{0,0,0}$e_2$}%
}}}}
\put(2271,215){\makebox(0,0)[lb]{\smash{{\SetFigFont{8}{9.6}{\familydefault}{\mddefault}{\updefault}{\color[rgb]{0,0,0}$e_1$}%
}}}}
\put(3005,-539){\makebox(0,0)[lb]{\smash{{\SetFigFont{8}{9.6}{\familydefault}{\mddefault}{\updefault}{\color[rgb]{0,0,0}$e_3$}%
}}}}
\put(4116,-494){\makebox(0,0)[lb]{\smash{{\SetFigFont{8}{9.6}{\familydefault}{\mddefault}{\updefault}{\color[rgb]{0,0,0}$e_2$}%
}}}}
\put(4009,197){\makebox(0,0)[lb]{\smash{{\SetFigFont{8}{9.6}{\familydefault}{\mddefault}{\updefault}{\color[rgb]{0,0,0}$e_1$}%
}}}}
\put(4310,-135){\makebox(0,0)[lb]{\smash{{\SetFigFont{8}{9.6}{\familydefault}{\mddefault}{\updefault}{\color[rgb]{0,0,0}$e_0$}%
}}}}
\put(1216,349){\makebox(0,0)[lb]{\smash{{\SetFigFont{8}{9.6}{\familydefault}{\mddefault}{\updefault}{\color[rgb]{0,0,0}$e_4$}%
}}}}
\end{picture}%

\end{center}

The pattern of circles in the center of the figure can by modified forming the
figures on the left and the right with each of the four circles intersecting
the adjacent two circles orthogonally.  If we leave the other faces of $P_0$
fixed we can make a small enough modification that the edges $e_1,e_2,e_3,e_4$
still have finite non-zero lengths.  Since each of the dihedral angles
corresponding to edges other than $e_0,e_1,e_2,e_3,$ and $e_4$ were chosen to
be $2\pi/5$, this small modification will not increase any of these angles past
$\pi/2$.  One of these modified patterns of intersecting circles will
correspond to a polyhedron in ${\cal P}^0_C$ and the other to a polyhedron in
${\cal P}^0_{C'}$.  \Endproof

\vspace{.1in}
\noindent {\bf Proof of Lemma \ref{ROEDER}:} 
We assume that $C \neq Pr_N$ is a simple abstract polyhedron with $N>7$ faces.
We will construct a sequence of Whitehead moves that change $C$ to $D_N$, so
that no intermediate complex has a prismatic 3-circuit.

Find a vertex $v_\infty$ of $C^*$ which is connected to the greatest number of
other vertices.  We will call the link of $v_\infty$, a cycle of $k$ vertices
and $k$ edges, the {\it outer polygon}.  Notice that $k \leq N-2$, with
equality precisely when $C = Pr_N$.  Therefore, since $C \neq Pr_N$ by
hypothesis, our first goal is to find Whitehead moves which increase $k$ to
$N-3$ without introducing any prismatic $3$-circuits along the way.  Once this
is completed, an easy sequence of Whitehead moves changes the resulting complex
to $D^*_N$.

Let us set up some notation.  Draw the dual complex $C^*$ in the
plane with the vertex $v_\infty$ at infinity, so that the outer polygon $P$
surrounding the remaining vertices and triangles.  We call the
vertices inside of $P$ {\it interior vertices}.  All of the
edges inside of $P$ which do not have an endpoint on $P$ are called {\it interior
edges}.  

Note that the graph of interior vertices and edges is connected, since
$C^*$ is simple.  An interior vertex which is connected to only one
other interior vertex will be called an {\it endpoint}.

\begin{center}
\begin{picture}(0,0)%
\epsfig{file=fig2.pstex}%
\end{picture}%
\setlength{\unitlength}{3947sp}%
\begingroup\makeatletter\ifx\SetFigFont\undefined%
\gdef\SetFigFont#1#2#3#4#5{%
  \reset@font\fontsize{#1}{#2pt}%
  \fontfamily{#3}\fontseries{#4}\fontshape{#5}%
  \selectfont}%
\fi\endgroup%
\begin{picture}(3387,1956)(1189,-2220)
\put(3578,-1771){\makebox(0,0)[lb]{\smash{{\SetFigFont{12}{14.4}{\familydefault}{\mddefault}{\updefault}{\color[rgb]{0,0,0}\small{$F^2_v$}}%
}}}}
\put(1223,-751){\makebox(0,0)[lb]{\smash{{\SetFigFont{12}{14.4}{\familydefault}{\mddefault}{\updefault}{\color[rgb]{0,0,0}\small{$F^1_w$}}%
}}}}
\put(1651,-1329){\makebox(0,0)[lb]{\smash{{\SetFigFont{12}{14.4}{\familydefault}{\mddefault}{\updefault}{\color[rgb]{0,0,0}\small{$w$}}%
}}}}
\put(3270,-638){\makebox(0,0)[lb]{\smash{{\SetFigFont{12}{14.4}{\familydefault}{\mddefault}{\updefault}{\color[rgb]{0,0,0}\small{$F^1_v$}}%
}}}}
\put(3256,-1179){\makebox(0,0)[lb]{\smash{{\SetFigFont{12}{14.4}{\familydefault}{\mddefault}{\updefault}{\color[rgb]{0,0,0}\small{$v$}}%
}}}}
\put(4576,-661){\makebox(0,0)[lb]{\smash{{\SetFigFont{12}{14.4}{\familydefault}{\mddefault}{\updefault}{\color[rgb]{0,0,0}\small{endpoint}}%
}}}}
\end{picture}%

\end{center}

Throughout this proof we will draw $P$ in black and we draw interior edges and
vertices of $C^*$ in black, as well.  The connections between $P$ and the
interior vertices will be in grey.  Connections between $P$ and $v_{\infty}$ will
be black, if shown at all.

The link of an interior vertex $v$ will intersect $P$ in a
number of components $F_v^1,\cdots,F_v^n$. (Possibly $n = 0$.)  See the
above figure.  We will say that $v$ is {\it connected to $P$ in these
components.}  Notice that since $C^*$ is simple, an endpoint is
always connected to $P$ in exactly one such component.

\begin{slem} \label{CHECK}
If a Whitehead move on the dual $C^*$ of an abstract polyhedron yields $C'^*$ (replacing $f$
by $f'$), and if $\delta$ is a simple closed path in $C^*$, which separates
one endpoint of $f'$ from the other, then any newly-created 3-circuit will
contain some vertex of $\delta$ which shares edges with both endpoints
of $f'$.
\end{slem}

\noindent {\bf Proof:}
A newly created 3-circuit $\gamma$ must contain the new edge $f'$
as well as two additional edges $e_1$ and $e_2$ connecting from a single vertex
$V$ to the two endpoints of $f'$.  By the Jordan Curve Theorem, since $\delta$
separates the endpoints of $f'$, the path $e_1e_2$
intersects $\delta$.  The path $e_1e_2$ is entirely in $C^*$
since $f'$ is the only new edge in $C'^*$,
so the vertex $V$ is in $\delta$. \Endproof

We now prove three additional sub-lemmas that specify certain Whitehead moves
that, when performed on a simple abstract polyhedron $C$ (which is not a prism
and has more than 7 faces), do not introduce any prismatic 3-circuits.  Hence
the resulting abstract polyhedron $C'^*$ is simple.  More specifically, we will use Sub-lemma
\ref{CHECK} to see that each Whitehead move introduces exactly two
newly-created 3-circuits in $C'^*$, the two triangles containing the new edge
$f'$.

\begin{slem} \label{SUBLEMMA1}
%Let $C$ satisfy the hypotheses of Lemma \ref{ROEDER},
%with outer polygon $P \subset C^*$.
Suppose that there is an interior vertex $A$ of $C^*$ which is connected to $P$ in exactly
one component consisting of exactly two consecutive vertices $Q$ and $R$.
The Whitehead move $Wh(QR)$ on $C^*$ increases the length of the
outer polygon by one, and introduces no prismatic $3$-circuit.
\end{slem}

\noindent {\bf Proof:}
\begin{center}
\begin{picture}(0,0)%
\epsfig{file=fig3.pstex}%
\end{picture}%
\setlength{\unitlength}{3947sp}%
\begingroup\makeatletter\ifx\SetFigFont\undefined%
\gdef\SetFigFont#1#2#3#4#5{%
  \reset@font\fontsize{#1}{#2pt}%
  \fontfamily{#3}\fontseries{#4}\fontshape{#5}%
  \selectfont}%
\fi\endgroup%
\begin{picture}(4176,1741)(1537,-2354)
\put(2316,-835){\makebox(0,0)[lb]{\smash{{\SetFigFont{12}{14.4}{\familydefault}{\mddefault}{\updefault}{\color[rgb]{0,0,0}\small{$v_\infty$}}%
}}}}
\put(5542,-2022){\makebox(0,0)[lb]{\smash{{\SetFigFont{8}{9.6}{\familydefault}{\mddefault}{\updefault}{\color[rgb]{0,0,0}\small{$E$}}%
}}}}
\put(4755,-769){\makebox(0,0)[lb]{\smash{{\SetFigFont{12}{14.4}{\familydefault}{\mddefault}{\updefault}{\color[rgb]{0,0,0}\small{$v_\infty$}}%
}}}}
\put(3402,-1387){\makebox(0,0)[lb]{\smash{{\SetFigFont{8}{9.6}{\familydefault}{\mddefault}{\updefault}{\color[rgb]{0,0,0}$Wh(QR)$}%
}}}}
\put(2828,-1223){\makebox(0,0)[lb]{\smash{{\SetFigFont{8}{9.6}{\familydefault}{\mddefault}{\updefault}{\color[rgb]{0,0,0}\small{$R$}}%
}}}}
\put(5315,-1216){\makebox(0,0)[lb]{\smash{{\SetFigFont{8}{9.6}{\familydefault}{\mddefault}{\updefault}{\color[rgb]{0,0,0}\small{$R$}}%
}}}}
\put(4509,-2236){\makebox(0,0)[lb]{\smash{{\SetFigFont{8}{9.6}{\familydefault}{\mddefault}{\updefault}{\color[rgb]{0,0,0}\small{interior stuff}}%
}}}}
\put(1704,-2243){\makebox(0,0)[lb]{\smash{{\SetFigFont{8}{9.6}{\familydefault}{\mddefault}{\updefault}{\color[rgb]{0,0,0}\small{other interior stuff}}%
}}}}
\put(2546,-2018){\makebox(0,0)[lb]{\smash{{\SetFigFont{8}{9.6}{\familydefault}{\mddefault}{\updefault}{\color[rgb]{0,0,0}\small{$A$}}%
}}}}
\put(3106,-2041){\makebox(0,0)[lb]{\smash{{\SetFigFont{8}{9.6}{\familydefault}{\mddefault}{\updefault}{\color[rgb]{0,0,0}\small{$E$}}%
}}}}
\put(1621,-2033){\makebox(0,0)[lb]{\smash{{\SetFigFont{8}{9.6}{\familydefault}{\mddefault}{\updefault}{\color[rgb]{0,0,0}\small{$D$}}%
}}}}
\put(1789,-1207){\makebox(0,0)[lb]{\smash{{\SetFigFont{8}{9.6}{\familydefault}{\mddefault}{\updefault}{\color[rgb]{0,0,0}\small{$Q$}}%
}}}}
\put(4267,-1214){\makebox(0,0)[lb]{\smash{{\SetFigFont{8}{9.6}{\familydefault}{\mddefault}{\updefault}{\color[rgb]{0,0,0}\small{$Q$}}%
}}}}
\put(5057,-2019){\makebox(0,0)[lb]{\smash{{\SetFigFont{8}{9.6}{\familydefault}{\mddefault}{\updefault}{\color[rgb]{0,0,0}\small{$A$}}%
}}}}
\put(4114,-2034){\makebox(0,0)[lb]{\smash{{\SetFigFont{8}{9.6}{\familydefault}{\mddefault}{\updefault}{\color[rgb]{0,0,0}\small{$D$}}%
}}}}
\end{picture}%

\end{center}
\vspace{.05in}

Clearly this Whitehead move increases the length of $P$ by one.  We apply
Sub-lemma \ref{CHECK} to see that this move introduces no prismatic 3-circuits.
We let $\delta = P$, the outer polygon, which clearly separates the interior
vertex $A$ from $v_\infty$ in $C^*$. Any new 3-circuit would consist of a point
on $P$ connected to both $A$ and $v_\infty$.  By hypothesis, there were only
the two points $Q$ and $R$ on $P$ connected to $A$.  These two points result in
the new triangles $QAv_\infty$ and $RAv_\infty$ in ${C'}^*$.  Therefore  $Wh(QR)$
result in no prismatic 3-circuits.  \Endproof

In the above sub-lemma, the condition that $A$ is connected to exactly two
consecutive vertices of $P$ prevents $A$ from being an endpoint.  For if $A$ is
an endpoint, let $B$ denote the unique interior vertex connected to $A$.  Then
$BQR$ would be a $3$-circuit in $C^*$ separating $A$ from the other vertices
and hence would contradict the hypothesis that $C$ is simple.  Therefore any
endpoint must be connected to $P$ in a single component having three or more
vertices.

\begin{slem} \label{SUBLEMMA2}
%Let $C$ satisfy the hypotheses of Lemma \ref{ROEDER},
%with outer polygon $P \subset C^*$.
Suppose that there is an interior vertex $A$ that is connected to $P$ in a
component consisting of $M$ consecutive vertices $Q_1,\cdots,Q_M$ of $P$ (and
possibly other components).

\noindent (a)
If $A$ is not an endpoint and $M > 2$, the sequence of Whitehead moves
$Wh(AQ_M),\ldots,Wh(AQ_3)$ results in a complex in which $A$ is connected
to the same component of $P$ in only $Q_1$ and $Q_2$.  These moves leave $P$
unchanged, and introduce no prismatic 3-circuit.

\vspace{.07in}
\begin{picture}(0,0)%
\epsfig{file=fig4.pstex}%
\end{picture}%
\setlength{\unitlength}{3947sp}%
\begingroup\makeatletter\ifx\SetFigFont\undefined%
\gdef\SetFigFont#1#2#3#4#5{%
  \reset@font\fontsize{#1}{#2pt}%
  \fontfamily{#3}\fontseries{#4}\fontshape{#5}%
  \selectfont}%
\fi\endgroup%
\begin{picture}(5193,1438)(1104,-1372)
\put(1902,-1341){\makebox(0,0)[lb]{\smash{{\SetFigFont{9}{10.8}{\familydefault}{\mddefault}{\updefault}{\color[rgb]{0,0,0}$A$}%
}}}}
\put(2319,-30){\makebox(0,0)[lb]{\smash{{\SetFigFont{9}{10.8}{\familydefault}{\mddefault}{\updefault}{\color[rgb]{0,0,0}$Q_{M-1}$}%
}}}}
\put(2797,-30){\makebox(0,0)[lb]{\smash{{\SetFigFont{9}{10.8}{\familydefault}{\mddefault}{\updefault}{\color[rgb]{0,0,0}$Q_M$}%
}}}}
\put(4348,-1341){\makebox(0,0)[lb]{\smash{{\SetFigFont{9}{10.8}{\familydefault}{\mddefault}{\updefault}{\color[rgb]{0,0,0}$D$}%
}}}}
\put(5959,-1341){\makebox(0,0)[lb]{\smash{{\SetFigFont{9}{10.8}{\familydefault}{\mddefault}{\updefault}{\color[rgb]{0,0,0}$E$}%
}}}}
\put(5124,-1341){\makebox(0,0)[lb]{\smash{{\SetFigFont{9}{10.8}{\familydefault}{\mddefault}{\updefault}{\color[rgb]{0,0,0}$A$}%
}}}}
\put(4348,-30){\makebox(0,0)[lb]{\smash{{\SetFigFont{9}{10.8}{\familydefault}{\mddefault}{\updefault}{\color[rgb]{0,0,0}$Q_1 Q_2 Q_3$}%
}}}}
\put(5481,-30){\makebox(0,0)[lb]{\smash{{\SetFigFont{9}{10.8}{\familydefault}{\mddefault}{\updefault}{\color[rgb]{0,0,0}$Q_{M-1}$}%
}}}}
\put(5898,-30){\makebox(0,0)[lb]{\smash{{\SetFigFont{9}{10.8}{\familydefault}{\mddefault}{\updefault}{\color[rgb]{0,0,0}$Q_M$}%
}}}}
\put(3035,-387){\makebox(0,0)[lb]{\smash{{\SetFigFont{9}{10.8}{\familydefault}{\mddefault}{\updefault}{\color[rgb]{0,0,0}$Wh(AQ_M)$}%
}}}}
\put(1126,-1341){\makebox(0,0)[lb]{\smash{{\SetFigFont{9}{10.8}{\familydefault}{\mddefault}{\updefault}{\color[rgb]{0,0,0}$D$}%
}}}}
\put(2737,-1341){\makebox(0,0)[lb]{\smash{{\SetFigFont{9}{10.8}{\familydefault}{\mddefault}{\updefault}{\color[rgb]{0,0,0}$E$}%
}}}}
\put(1126,-30){\makebox(0,0)[lb]{\smash{{\SetFigFont{9}{10.8}{\familydefault}{\mddefault}{\updefault}{\color[rgb]{0,0,0}$Q_1 Q_2 Q_3$}%
}}}}
\end{picture}%

\vspace{.07in}

\noindent (b)
If $A$ is an endpoint and $M > 3$, the sequence of Whitehead moves
\newline $Wh(AQ_M),\ldots,Wh(AQ_4)$ results in a complex in which $A$ is connected
to the same component of $P$ in only $Q_1,Q_2$, and $Q_3$.  These moves leave $P$
unchanged and introduce no prismatic 3-circuits.

\vspace{.07in}
\begin{picture}(0,0)%
\epsfig{file=./fig8.pstex}%
\end{picture}%
\setlength{\unitlength}{3947sp}%
\begingroup\makeatletter\ifx\SetFigFont\undefined%
\gdef\SetFigFont#1#2#3#4#5{%
  \reset@font\fontsize{#1}{#2pt}%
  \fontfamily{#3}\fontseries{#4}\fontshape{#5}%
  \selectfont}%
\fi\endgroup%
\begin{picture}(5002,1664)(571,-1179)
\put(4918,-1121){\makebox(0,0)[lb]{\smash{{\SetFigFont{12}{14.4}{\familydefault}{\mddefault}{\updefault}{\color[rgb]{0,0,0}\small{$Q_1$}}%
}}}}
\put(4801,-511){\makebox(0,0)[lb]{\smash{{\SetFigFont{12}{14.4}{\familydefault}{\mddefault}{\updefault}{\color[rgb]{0,0,0}\small{$E$}}%
}}}}
\put(4471,-495){\makebox(0,0)[lb]{\smash{{\SetFigFont{12}{14.4}{\familydefault}{\mddefault}{\updefault}{\color[rgb]{0,0,0}\small{$A$}}%
}}}}
\put(4171,171){\makebox(0,0)[lb]{\smash{{\SetFigFont{12}{14.4}{\familydefault}{\mddefault}{\updefault}{\color[rgb]{0,0,0}\small{$Q_{M-1}$}}%
}}}}
\put(4884,329){\makebox(0,0)[lb]{\smash{{\SetFigFont{12}{14.4}{\familydefault}{\mddefault}{\updefault}{\color[rgb]{0,0,0}\small{$Q_M$}}%
}}}}
\put(4051,-811){\makebox(0,0)[lb]{\smash{{\SetFigFont{12}{14.4}{\familydefault}{\mddefault}{\updefault}{\color[rgb]{0,0,0}\small{$Q_2$}}%
}}}}
\put(3721,-422){\makebox(0,0)[lb]{\smash{{\SetFigFont{12}{14.4}{\familydefault}{\mddefault}{\updefault}{\color[rgb]{0,0,0}\small{$Q_4$}}%
}}}}
\put(3796,-615){\makebox(0,0)[lb]{\smash{{\SetFigFont{12}{14.4}{\familydefault}{\mddefault}{\updefault}{\color[rgb]{0,0,0}\small{$Q_3$}}%
}}}}
\put(3751, 14){\makebox(0,0)[lb]{\smash{{\SetFigFont{12}{14.4}{\familydefault}{\mddefault}{\updefault}{\color[rgb]{0,0,0}\small{$Q_{M-2}$}}%
}}}}
\put(2626,-286){\makebox(0,0)[lb]{\smash{{\SetFigFont{12}{14.4}{\familydefault}{\mddefault}{\updefault}{\color[rgb]{0,0,0}\small{$Wh(AQ_M)$}}%
}}}}
\put(1651,-511){\makebox(0,0)[lb]{\smash{{\SetFigFont{12}{14.4}{\familydefault}{\mddefault}{\updefault}{\color[rgb]{0,0,0}\small{$E$}}%
}}}}
\put(1321,-495){\makebox(0,0)[lb]{\smash{{\SetFigFont{12}{14.4}{\familydefault}{\mddefault}{\updefault}{\color[rgb]{0,0,0}\small{$A$}}%
}}}}
\put(901,-811){\makebox(0,0)[lb]{\smash{{\SetFigFont{12}{14.4}{\familydefault}{\mddefault}{\updefault}{\color[rgb]{0,0,0}\small{$Q_2$}}%
}}}}
\put(571,-422){\makebox(0,0)[lb]{\smash{{\SetFigFont{12}{14.4}{\familydefault}{\mddefault}{\updefault}{\color[rgb]{0,0,0}\small{$Q_4$}}%
}}}}
\put(1021,171){\makebox(0,0)[lb]{\smash{{\SetFigFont{12}{14.4}{\familydefault}{\mddefault}{\updefault}{\color[rgb]{0,0,0}\small{$Q_{M-1}$}}%
}}}}
\put(1734,329){\makebox(0,0)[lb]{\smash{{\SetFigFont{12}{14.4}{\familydefault}{\mddefault}{\updefault}{\color[rgb]{0,0,0}\small{$Q_M$}}%
}}}}
\put(646,-615){\makebox(0,0)[lb]{\smash{{\SetFigFont{12}{14.4}{\familydefault}{\mddefault}{\updefault}{\color[rgb]{0,0,0}\small{$Q_3$}}%
}}}}
\put(601, 14){\makebox(0,0)[lb]{\smash{{\SetFigFont{12}{14.4}{\familydefault}{\mddefault}{\updefault}{\color[rgb]{0,0,0}\small{$Q_{M-2}$}}%
}}}}
\put(1767,-1112){\makebox(0,0)[lb]{\smash{{\SetFigFont{12}{14.4}{\familydefault}{\mddefault}{\updefault}{\color[rgb]{0,0,0}\small{$Q_1$}}%
}}}}
\end{picture}%
~                                                                          
\vspace{.07in}
\end{slem}

{\bf Proof:}
Part (a) $A$ is not an endpoint.  Clearly the move $Wh(AQ_M)$ decreases $M$ by
one.  We check that if $M>2$, this move introduces no prismatic 3-circuits.  We
let $\delta$ be the path $v_\infty Q_{M-2}AQ_M$ which separates $Q_{M-1}$ from
$E$ in $C^*$.  By Sub-lemma \ref{CHECK}, any new  3-circuit contains a vertex
on $\delta$ connected to both $E$ and $Q_{M-1}$.  Clearly $v_\infty$ is not
connected to the interior vertex $E$.   If $M > 3$, $Q_{M-2}$ is connected only
to $Q_{M-1}$, $A$, $Q_{M-3}$, and $v_\infty$.  Otherwise, when $M=3$, a
connection of $Q_1$ to $E$ would mean that $C^*$ had a $3$-cycle $EQ_1A$
which would have to separate $D$ (in the figure above corresponding to part (a))
from $Q_2$, contrary to the hypothesis that $C^*$ is simple.

Hence, the only two vertices on $\delta$ that are connected to both $E$ and
$Q_{M-1}$ are $A$ and $Q_M$, forming the two triangles $AQ_{M-1}E$ and
$Q_MQ_{M-1}E$ in $C'^*$.  Hence there are no new prismatic 3-circuits, so we
can reduce $M$ by one, when $M > 2$.

Part (b) $A$ is an endpoint.  We again use $\delta = v_\infty Q_{M-2}AQ_M$ to
check that the move $Wh(AQ_M)$ introduces no prismatic $3$-circuits.  The proof is
identical to Part (a), except that $M > 3$ is needed to conclude that $Q_{M-2}$
is not connected to $E$ since $A$ is now an endpoint.  

Therefore, as long as $M
> 3$ we can reduce $M$ by one without introducing prismatic 3-circuits.  Recall
that an endpoint of a simple complex cannot be connected to fewer than three
points of $P$, so this is optimal. \Endproof

Note: In both parts (a) and (b), each of the Whitehead moves $Wh(AQ_M)$
transfers the connection between $A$ and $Q_M$ to a connection between the
neighboring interior vertex $E$ and $Q_{M-1}$.  In fact $Q_{M-1}$ gets added to
the component containing $Q_M$ in which $E$ is connected to $P$.  This is
helpful later on, in Case 2 of Lemma \ref{ROEDER}.

\begin{slem}\label{SUBLEMMA3}
%Let $C$ satisfy the hypotheses of Lemma \ref{ROEDER},
%with outer polygon $P \subset C^*$.
Suppose that there is an interior vertex $A$ whose link contains two distinct
vertices $X$ and $Y$ of $P$.  Then there are Whitehead moves which eliminate
any component in which $A$ is connected to $P$, if that component does not
contain $X$ or $Y$.  $P$ is unchanged, and no prismatic $3$-circuits will be
introduced.
\end{slem}

\noindent
{\bf Example:}  
\vspace{.07in}
\begin{center}
\begin{picture}(0,0)%
\includegraphics{fig5.pstex}%
\end{picture}%
\setlength{\unitlength}{3947sp}%
\begingroup\makeatletter\ifx\SetFigFont\undefined%
\gdef\SetFigFont#1#2#3#4#5{%
  \reset@font\fontsize{#1}{#2pt}%
  \fontfamily{#3}\fontseries{#4}\fontshape{#5}%
  \selectfont}%
\fi\endgroup%
\begin{picture}(2959,1013)(879,-1352)
\put(1126,-811){\makebox(0,0)[lb]{\smash{\SetFigFont{12}{14.4}{\familydefault}{\mddefault}{\updefault}{\color[rgb]{0,0,0}\small{$A$}}%
}}}
\put(2994,-849){\makebox(0,0)[lb]{\smash{\SetFigFont{12}{14.4}{\familydefault}{\mddefault}{\updefault}{\color[rgb]{0,0,0}\small{$A$}}%
}}}
\put(2026,-511){\makebox(0,0)[lb]{\smash{\SetFigFont{6}{7.2}{\familydefault}{\mddefault}{\updefault}{\color[rgb]{0,0,0}\small{$X$}}%
}}}
\put(3826,-511){\makebox(0,0)[lb]{\smash{\SetFigFont{6}{7.2}{\familydefault}{\mddefault}{\updefault}{\color[rgb]{0,0,0}\small{$X$}}%
}}}
\put(2044,-1089){\makebox(0,0)[lb]{\smash{\SetFigFont{6}{7.2}{\familydefault}{\mddefault}{\updefault}{\color[rgb]{0,0,0}\small{$Y$}}%
}}}
\put(3809,-1094){\makebox(0,0)[lb]{\smash{\SetFigFont{6}{7.2}{\familydefault}{\mddefault}{\updefault}{\color[rgb]{0,0,0}\small{$Y$}}%
}}}
\end{picture}

\end{center}
\vspace{.07in}

\noindent
Here $A$ is connected to $P$ in four components containing six vertices.  We can
eliminate connections of $A$ to all of the components except for the
single-point components $X$ and $Y$.

{\bf Proof:} Let $O$ be a component not containing $X$ or $Y$.  
If $O$ contains more than two vertices, we can reduce it to two vertices
by Sub-lemma \ref{SUBLEMMA2}(a). 

Suppose that $O$ contains exactly two vertices, $V$ and $W$.  We check that the move
$Wh(AW)$ eliminates the connection from $A$ to $W$ without introducing
prismatic $3$-circuits.  Let $D$ be the unique interior vertex forming triangle
$ADW$, as in the figure below.  The move $Wh(AW)$ creates the new edge $DV$.
Let $\delta$ be the loop $v_\infty Y A W$ which separates $D$ from $V$ in
$C^*$.  See the dashed curve in the figure below.  By Sub-lemma
\ref{CHECK}, any new 3-circuit contains a point on $\delta$ that is connected
to both $D$ and $V$.  Clearly $v_\infty$ is not connected to the interior
vertex $D$.  Since $Y$ and $V$ are in different components of connection
between $A$ and $P$, $Y$ is not connected to $V$.  Therefore, only $A$ and $W$
are connected to both $D$ and $V$, forming the triangles $ADV$ and $WVD$ in
$C'^*$.  Therefore, $Wh(AW)$ results in no prismatic 3-circuits.

\begin{center}
\vspace{.07in}
\begin{picture}(0,0)%
\epsfig{file=./new_fig6a.pstex}%
\end{picture}%
\setlength{\unitlength}{3947sp}%
\begingroup\makeatletter\ifx\SetFigFont\undefined%
\gdef\SetFigFont#1#2#3#4#5{%
  \reset@font\fontsize{#1}{#2pt}%
  \fontfamily{#3}\fontseries{#4}\fontshape{#5}%
  \selectfont}%
\fi\endgroup%
\begin{picture}(5347,1010)(453,-816)
\put(811,-631){\makebox(0,0)[lb]{\smash{{\SetFigFont{6}{7.2}{\familydefault}{\mddefault}{\updefault}{\color[rgb]{0,0,0}\small{$X$}}%
}}}}
\put(2450,-786){\makebox(0,0)[lb]{\smash{{\SetFigFont{6}{7.2}{\familydefault}{\mddefault}{\updefault}{\color[rgb]{0,0,0}\small{$v_\infty$}}%
}}}}
\put(2119,-554){\makebox(0,0)[lb]{\smash{{\SetFigFont{6}{7.2}{\familydefault}{\mddefault}{\updefault}{\color[rgb]{0,0,0}\small{$V$}}%
}}}}
\put(4435,-325){\makebox(0,0)[lb]{\smash{{\SetFigFont{6}{7.2}{\familydefault}{\mddefault}{\updefault}{\color[rgb]{0,0,0}\small{$A$}}%
}}}}
\put(4915, 95){\makebox(0,0)[lb]{\smash{{\SetFigFont{6}{7.2}{\familydefault}{\mddefault}{\updefault}{\color[rgb]{0,0,0}\small{$D$}}%
}}}}
\put(4150,-685){\makebox(0,0)[lb]{\smash{{\SetFigFont{6}{7.2}{\familydefault}{\mddefault}{\updefault}{\color[rgb]{0,0,0}\small{$X$}}%
}}}}
\put(4135,-250){\makebox(0,0)[lb]{\smash{{\SetFigFont{6}{7.2}{\familydefault}{\mddefault}{\updefault}{\color[rgb]{0,0,0}\small{$Y$}}%
}}}}
\put(2403,-70){\makebox(0,0)[lb]{\smash{{\SetFigFont{6}{7.2}{\familydefault}{\mddefault}{\updefault}{\color[rgb]{0,0,0}\small{$v_\infty$}}%
}}}}
\put(5695,-122){\makebox(0,0)[lb]{\smash{{\SetFigFont{6}{7.2}{\familydefault}{\mddefault}{\updefault}{\color[rgb]{0,0,0}\small{$v_\infty$}}%
}}}}
\put(5800,-789){\makebox(0,0)[lb]{\smash{{\SetFigFont{6}{7.2}{\familydefault}{\mddefault}{\updefault}{\color[rgb]{0,0,0}\small{$v_\infty$}}%
}}}}
\put(453,-595){\makebox(0,0)[lb]{\smash{{\SetFigFont{6}{7.2}{\familydefault}{\mddefault}{\updefault}{\color[rgb]{0,0,0}\small{$v_\infty$}}%
}}}}
\put(461,-145){\makebox(0,0)[lb]{\smash{{\SetFigFont{6}{7.2}{\familydefault}{\mddefault}{\updefault}{\color[rgb]{0,0,0}\small{$v_\infty$}}%
}}}}
\put(3768,-205){\makebox(0,0)[lb]{\smash{{\SetFigFont{6}{7.2}{\familydefault}{\mddefault}{\updefault}{\color[rgb]{0,0,0}\small{$v_\infty$}}%
}}}}
\put(3767,-640){\makebox(0,0)[lb]{\smash{{\SetFigFont{6}{7.2}{\familydefault}{\mddefault}{\updefault}{\color[rgb]{0,0,0}\small{$v_\infty$}}%
}}}}
\put(2725,-259){\makebox(0,0)[lb]{\smash{{\SetFigFont{6}{7.2}{\familydefault}{\mddefault}{\updefault}{\color[rgb]{0,0,0}\small{$Wh(AW)$}}%
}}}}
\put(1591,134){\makebox(0,0)[lb]{\smash{{\SetFigFont{6}{7.2}{\familydefault}{\mddefault}{\updefault}{\color[rgb]{0,0,0}\small{$D$}}%
}}}}
\put(796,-196){\makebox(0,0)[lb]{\smash{{\SetFigFont{6}{7.2}{\familydefault}{\mddefault}{\updefault}{\color[rgb]{0,0,0}\small{$Y$}}%
}}}}
\put(1096,-271){\makebox(0,0)[lb]{\smash{{\SetFigFont{6}{7.2}{\familydefault}{\mddefault}{\updefault}{\color[rgb]{0,0,0}\small{$A$}}%
}}}}
\put(1581,-206){\makebox(0,0)[lb]{\smash{{\SetFigFont{6}{7.2}{\familydefault}{\mddefault}{\updefault}{\color[rgb]{0,0,0}\small{$\delta$}}%
}}}}
\put(2121,-331){\makebox(0,0)[lb]{\smash{{\SetFigFont{6}{7.2}{\familydefault}{\mddefault}{\updefault}{\color[rgb]{0,0,0}\small{$W$}}%
}}}}
\put(5465,-560){\makebox(0,0)[lb]{\smash{{\SetFigFont{6}{7.2}{\familydefault}{\mddefault}{\updefault}{\color[rgb]{0,0,0}\small{$V$}}%
}}}}
\put(5434,-385){\makebox(0,0)[lb]{\smash{{\SetFigFont{6}{7.2}{\familydefault}{\mddefault}{\updefault}{\color[rgb]{0,0,0}\small{$W$}}%
}}}}
\end{picture}%

\vspace{.07in}
\end{center}

Thus, we can suppose that $O$ contains a single vertex $V$.  We check that the
move $Wh(AV)$, which eliminates this connection, does not introduce any
prismatic $3$-circuits.  Let $D$ and $E$ be the unique interior vertices
forming the triangles $ADV$ and $AEV$.  Let $\delta_1$ be the curve $v_\infty
YAV$ and $\delta_2$ be the curve $v_\infty XAV$ in $C^*$.  See the two dashed
curves in the figure below.  Both of these curves separate $D$ and $E$ in
$C^*$.  Applying Sub-lemma \ref{CHECK} twice, we conclude that any newly
created 3-circuit contains a point that is {\it both on $\delta_1$ and on
$\delta_2$} and that connects to both $D$ and $E$.  The only points on both
$\delta_1$ and $\delta_2$ are $v_\infty, A$, and $V$.  Since $D$ and $E$ are
interior, $v_\infty$ cannot connect to either of them.  The connections from
$A$ and from $V$ to $D$ and $E$ result in the triangles $ADE$ and $VDE$ in
$C'^*$.  Therefore, $Wh(AV)$ results in no prismatic 3-circuits. 

\vspace{.07in}
\begin{center}
\begin{picture}(0,0)%
\epsfig{file=./new_fig6b.pstex}%
\end{picture}%
\setlength{\unitlength}{3947sp}%
\begingroup\makeatletter\ifx\SetFigFont\undefined%
\gdef\SetFigFont#1#2#3#4#5{%
  \reset@font\fontsize{#1}{#2pt}%
  \fontfamily{#3}\fontseries{#4}\fontshape{#5}%
  \selectfont}%
\fi\endgroup%
\begin{picture}(5355,1210)(453,-1016)
\put(5507,-607){\makebox(0,0)[lb]{\smash{{\SetFigFont{6}{7.2}{\familydefault}{\mddefault}{\updefault}{\color[rgb]{0,0,0}\small{$V$}}%
}}}}
\put(3864,-614){\makebox(0,0)[lb]{\smash{{\SetFigFont{6}{7.2}{\familydefault}{\mddefault}{\updefault}{\color[rgb]{0,0,0}\small{$v_\infty$}}%
}}}}
\put(3872,-164){\makebox(0,0)[lb]{\smash{{\SetFigFont{6}{7.2}{\familydefault}{\mddefault}{\updefault}{\color[rgb]{0,0,0}\small{$v_\infty$}}%
}}}}
\put(4222,-650){\makebox(0,0)[lb]{\smash{{\SetFigFont{6}{7.2}{\familydefault}{\mddefault}{\updefault}{\color[rgb]{0,0,0}\small{$X$}}%
}}}}
\put(4207,-215){\makebox(0,0)[lb]{\smash{{\SetFigFont{6}{7.2}{\familydefault}{\mddefault}{\updefault}{\color[rgb]{0,0,0}\small{$Y$}}%
}}}}
\put(4507,-290){\makebox(0,0)[lb]{\smash{{\SetFigFont{6}{7.2}{\familydefault}{\mddefault}{\updefault}{\color[rgb]{0,0,0}\small{$A$}}%
}}}}
\put(5002,115){\makebox(0,0)[lb]{\smash{{\SetFigFont{6}{7.2}{\familydefault}{\mddefault}{\updefault}{\color[rgb]{0,0,0}\small{$D$}}%
}}}}
\put(5056,-994){\makebox(0,0)[lb]{\smash{{\SetFigFont{6}{7.2}{\familydefault}{\mddefault}{\updefault}{\color[rgb]{0,0,0}\small{$E$}}%
}}}}
\put(453,-595){\makebox(0,0)[lb]{\smash{{\SetFigFont{6}{7.2}{\familydefault}{\mddefault}{\updefault}{\color[rgb]{0,0,0}\small{$v_\infty$}}%
}}}}
\put(461,-145){\makebox(0,0)[lb]{\smash{{\SetFigFont{6}{7.2}{\familydefault}{\mddefault}{\updefault}{\color[rgb]{0,0,0}\small{$v_\infty$}}%
}}}}
\put(811,-631){\makebox(0,0)[lb]{\smash{{\SetFigFont{6}{7.2}{\familydefault}{\mddefault}{\updefault}{\color[rgb]{0,0,0}\small{$X$}}%
}}}}
\put(796,-196){\makebox(0,0)[lb]{\smash{{\SetFigFont{6}{7.2}{\familydefault}{\mddefault}{\updefault}{\color[rgb]{0,0,0}\small{$Y$}}%
}}}}
\put(1096,-271){\makebox(0,0)[lb]{\smash{{\SetFigFont{6}{7.2}{\familydefault}{\mddefault}{\updefault}{\color[rgb]{0,0,0}\small{$A$}}%
}}}}
\put(1591,134){\makebox(0,0)[lb]{\smash{{\SetFigFont{6}{7.2}{\familydefault}{\mddefault}{\updefault}{\color[rgb]{0,0,0}\small{$D$}}%
}}}}
\put(1656,-960){\makebox(0,0)[lb]{\smash{{\SetFigFont{6}{7.2}{\familydefault}{\mddefault}{\updefault}{\color[rgb]{0,0,0}\small{$E$}}%
}}}}
\put(2898,-265){\makebox(0,0)[lb]{\smash{{\SetFigFont{6}{7.2}{\familydefault}{\mddefault}{\updefault}{\color[rgb]{0,0,0}\small{$Wh(AV)$}}%
}}}}
\put(2405,-400){\makebox(0,0)[lb]{\smash{{\SetFigFont{6}{7.2}{\familydefault}{\mddefault}{\updefault}{\color[rgb]{0,0,0}\small{$v_\infty$}}%
}}}}
\put(5808,-420){\makebox(0,0)[lb]{\smash{{\SetFigFont{6}{7.2}{\familydefault}{\mddefault}{\updefault}{\color[rgb]{0,0,0}\small{$v_\infty$}}%
}}}}
\put(1576,-261){\makebox(0,0)[lb]{\smash{{\SetFigFont{6}{7.2}{\familydefault}{\mddefault}{\updefault}{\color[rgb]{0,0,0}\small{$\delta_1$}}%
}}}}
\put(1582,-557){\makebox(0,0)[lb]{\smash{{\SetFigFont{6}{7.2}{\familydefault}{\mddefault}{\updefault}{\color[rgb]{0,0,0}\small{$\delta_2$}}%
}}}}
\put(2087,-562){\makebox(0,0)[lb]{\smash{{\SetFigFont{6}{7.2}{\familydefault}{\mddefault}{\updefault}{\color[rgb]{0,0,0}\small{$V$}}%
}}}}
\end{picture}%

\end{center}
\Endproof

The proof that this move does not introduce any prismatic 3-circuit
depends essentially on the fact that $A$ is connected to
$P$ in at least two other vertices $X$ and $Y$.  Andreev describes a nearly
identical process to Sublemma \ref{SUBLEMMA3} in his paper \cite{AND} on pages
333-334.  However, he merely assumes that $A$ is connected to $P$ in at least
one component in addition to the components being eliminated.  He does not
require that $A$ is connected to $P$ in at least {\it two vertices} outside of
the components being eliminated. Andreev then asserts: ``It is readily seen
that all of the polyhedra obtained in this way are simple...'' In fact, the
Whitehead move demonstrated below clearly creates a prismatic
3-circuit.  (Here, $M$ and $N$ lie in $P$.)

\vspace{.07in}
\begin{center}
\begin{picture}(0,0)%
\epsfig{file=prismatic.pstex}%
\end{picture}%
\setlength{\unitlength}{3947sp}%
\begingroup\makeatletter\ifx\SetFigFont\undefined%
\gdef\SetFigFont#1#2#3#4#5{%
  \reset@font\fontsize{#1}{#2pt}%
  \fontfamily{#3}\fontseries{#4}\fontshape{#5}%
  \selectfont}%
\fi\endgroup%
\begin{picture}(4552,1449)(398,-719)
\put(4218, 77){\makebox(0,0)[lb]{\smash{{\SetFigFont{12}{14.4}{\familydefault}{\mddefault}{\updefault}{\color[rgb]{0,0,0}\small{$A$}}%
}}}}
\put(2254,136){\makebox(0,0)[lb]{\smash{{\SetFigFont{12}{14.4}{\familydefault}{\mddefault}{\updefault}{\color[rgb]{0,0,0}\small{$Wh(AN)$}}%
}}}}
\put(4354,-653){\makebox(0,0)[lb]{\smash{{\SetFigFont{12}{14.4}{\familydefault}{\mddefault}{\updefault}{\color[rgb]{0,0,0}\small{$N$}}%
}}}}
\put(4392,574){\makebox(0,0)[lb]{\smash{{\SetFigFont{12}{14.4}{\familydefault}{\mddefault}{\updefault}{\color[rgb]{0,0,0}\small{$M$}}%
}}}}
\put(4943, 29){\makebox(0,0)[lb]{\smash{{\SetFigFont{12}{14.4}{\familydefault}{\mddefault}{\updefault}{\color[rgb]{0,0,0}\small{$E$}}%
}}}}
\put(3496, 29){\makebox(0,0)[lb]{\smash{{\SetFigFont{12}{14.4}{\familydefault}{\mddefault}{\updefault}{\color[rgb]{0,0,0}\small{$D$}}%
}}}}
\put(1172, 39){\makebox(0,0)[lb]{\smash{{\SetFigFont{12}{14.4}{\familydefault}{\mddefault}{\updefault}{\color[rgb]{0,0,0}\small{$A$}}%
}}}}
\put(1187,566){\makebox(0,0)[lb]{\smash{{\SetFigFont{12}{14.4}{\familydefault}{\mddefault}{\updefault}{\color[rgb]{0,0,0}\small{$M$}}%
}}}}
\put(1187,-661){\makebox(0,0)[lb]{\smash{{\SetFigFont{12}{14.4}{\familydefault}{\mddefault}{\updefault}{\color[rgb]{0,0,0}\small{$N$}}%
}}}}
\put(1740, 14){\makebox(0,0)[lb]{\smash{{\SetFigFont{12}{14.4}{\familydefault}{\mddefault}{\updefault}{\color[rgb]{0,0,0}\small{$E$}}%
}}}}
\put(398, 14){\makebox(0,0)[lb]{\smash{{\SetFigFont{12}{14.4}{\familydefault}{\mddefault}{\updefault}{\color[rgb]{0,0,0}\small{$D$}}%
}}}}
\end{picture}%

\end{center}

Having assumed this stronger (and incorrect) version of Sub-lemma
\ref{SUBLEMMA3}, the remainder of Andreev's proof is relatively easy.
Unfortunately, the situation pictured above is not uncommon (as we will see in
Case 3 below!) Restricted to the weaker hypotheses of Sub-lemma \ref{SUBLEMMA3}
we will have to work a little bit harder.

We will now use these three sub-lemmas to show that if the length of $P$ is
less than $N-3$ (so that there are at least 3 interior vertices), then we
can do Whitehead moves to increase the length of $P$ by one, without
introducing any prismatic $3$-circuits.

{\bf Case 1:} An interior point which isn't an endpoint connects to $P$ in
a component with two or more vertices, and possibly in other components.

Apply Sub-lemma \ref{SUBLEMMA2}(a) decreasing this component to two vertices.  We can
then apply Sub-lemma \ref{SUBLEMMA3}, eliminating any other components since this component
contains two vertices. Finally, apply Sub-lemma \ref{SUBLEMMA1} to increase the length of
the outer polygon by 1.

{\bf Case 2:}
An interior vertex that is an endpoint is connected to more than three
vertices of $P$.

We assume that each of the interior points that are not endpoints are connected
to $P$ in components consisting of single points, otherwise we are in Case 1.
Let $A$ be the endpoint which is connected to more than three vertices of $P$.
By Sub-lemma \ref{SUBLEMMA2}(b), there is a Whitehead move that
transfers one of these connections to the interior vertex $E$ that is next to
$A$.  The point $E$ is not an endpoint since the interior graph is connected
and the assumption $N-k > 3$ implies that there are at least three interior
vertices.  Now, one of the components in which $E$ is connected to $P$ has
exactly two vertices, so we can then apply Case 1 for vertex $E$.

{\bf Case 3:}  Each interior vertex which is an endpoint is connected to
exactly $3$ points of $P$ and every other interior vertex is 
connected to $P$ in components each consisting of a single vertex.

First, notice that if the interior vertices and edges form a line segment, this
restriction on how interior points are connected to $P$ results in the prism,
contrary to hypothesis of this lemma.  However, there are many
complexes satisfying the hypotheses of this case which have interior vertices
and edges forming a graph more complicated than a line segment:

\vspace{.05in}
\begin{center}
\begin{picture}(0,0)%
\includegraphics{fig9.pstex}%
\end{picture}%
\setlength{\unitlength}{3947sp}%
\begingroup\makeatletter\ifx\SetFigFont\undefined%
\gdef\SetFigFont#1#2#3#4#5{%
  \reset@font\fontsize{#1}{#2pt}%
  \fontfamily{#3}\fontseries{#4}\fontshape{#5}%
  \selectfont}%
\fi\endgroup%
\begin{picture}(4724,1469)(879,-1283)
\end{picture}

\end{center}

\noindent
For such complexes we need a very special sequence of Whitehead moves to
increase the length of $P$.

Pick an interior vertex which is an endpoint and label it $I_1$.  Denote by
$P_1$, $P_2$, and $P_3$
the three vertices of $P$ to which $I_1$ connects.  $I_1$
will be connected to a sequence of interior vertices $I_2, I_3, \cdots
I_m, m \ge 2$, with $I_m$ the first interior vertex in the sequence that is connected to
more than two other interior vertices.  Vertex $I_m$ must exist by the assumption
that the interior vertices don't form a line segment, the configuration that we
ruled out above.  By hypothesis, $I_2,\cdots,I_m$ can only connect to $P$ in
components which each consist of a vertex, hence each must be connected to $P_1$
and to $P_3$.  Similarly, there is an interior vertex (call it $X$) which
connects both to $I_m$ and to $P_1$ and another vertex $Y$ which connects to
$I_m$ and $P_3$.
Vertex $I_m$ may connect to other vertices of $P$ and other interior vertices, as
shown on the left side of the following diagram.

\vspace{.07in}
\begin{center}
\begin{picture}(0,0)%
\epsfig{file=fig10.pstex}%
\end{picture}%
\setlength{\unitlength}{3947sp}%
\begingroup\makeatletter\ifx\SetFigFont\undefined%
\gdef\SetFigFont#1#2#3#4#5{%
  \reset@font\fontsize{#1}{#2pt}%
  \fontfamily{#3}\fontseries{#4}\fontshape{#5}%
  \selectfont}%
\fi\endgroup%
\begin{picture}(4185,1639)(16,-1844)
\put(1861,-915){\makebox(0,0)[lb]{\smash{{\SetFigFont{5}{6.0}{\familydefault}{\mddefault}{\updefault}{\color[rgb]{0,0,0}\small{$I_{m-1}$}}%
}}}}
\put(2851,-361){\makebox(0,0)[lb]{\smash{{\SetFigFont{12}{14.4}{\familydefault}{\mddefault}{\updefault}{\color[rgb]{0,0,0}\small{$P_1$}}%
}}}}
\put(4201,-961){\makebox(0,0)[lb]{\smash{{\SetFigFont{12}{14.4}{\familydefault}{\mddefault}{\updefault}{\color[rgb]{0,0,0}\small{$P_2$}}%
}}}}
\put(2701,-1786){\makebox(0,0)[lb]{\smash{{\SetFigFont{12}{14.4}{\familydefault}{\mddefault}{\updefault}{\color[rgb]{0,0,0}\small{$P_3$}}%
}}}}
\put(1051,-361){\makebox(0,0)[lb]{\smash{{\SetFigFont{12}{14.4}{\familydefault}{\mddefault}{\updefault}{\color[rgb]{0,0,0}\small{$X$}}%
}}}}
\put(1051,-1711){\makebox(0,0)[lb]{\smash{{\SetFigFont{12}{14.4}{\familydefault}{\mddefault}{\updefault}{\color[rgb]{0,0,0}\small{$Y$}}%
}}}}
\put(136,-961){\makebox(0,0)[lb]{\smash{{\SetFigFont{12}{14.4}{\familydefault}{\mddefault}{\updefault}{\color[rgb]{0,0,0}\tiny{other vertices}}%
}}}}
\put(1201,-833){\makebox(0,0)[lb]{\smash{{\SetFigFont{12}{14.4}{\familydefault}{\mddefault}{\updefault}{\color[rgb]{0,0,0}\small{$I_m$}}%
}}}}
\put(3601,-916){\makebox(0,0)[lb]{\smash{{\SetFigFont{12}{14.4}{\familydefault}{\mddefault}{\updefault}{\color[rgb]{0,0,0}\small{$I_1$}}%
}}}}
\put(2784,-909){\makebox(0,0)[lb]{\smash{{\SetFigFont{12}{14.4}{\familydefault}{\mddefault}{\updefault}{\color[rgb]{0,0,0}\small{$I_3$}}%
}}}}
\put(3084,-908){\makebox(0,0)[lb]{\smash{{\SetFigFont{12}{14.4}{\familydefault}{\mddefault}{\updefault}{\color[rgb]{0,0,0}\small{$I_2$}}%
}}}}
\put(2194,-914){\makebox(0,0)[lb]{\smash{{\SetFigFont{5}{6.0}{\familydefault}{\mddefault}{\updefault}{\color[rgb]{0,0,0}\small{$I_{m-2}$}}%
}}}}
\end{picture}%

\end{center}
\vspace{.07in}

Now we describe a sequence of Whitehead moves that can be used to connect $I_m$
to $P$ in only $P_1$ and $P_2$.  This will allow us to use Sub-lemma
\ref{SUBLEMMA1} to increase the length of $P$ by one.

First, using  Sub-lemma \ref{SUBLEMMA3}, one can eliminate
all possible connections of $I_m$ to $P$ in places other than $P_1$ and $P_3$.
Next, we do the move $Wh(I_mP_3)$ so that $I_m$ connects to $P$ only in
$P_1$.  We check that this Whitehead move does not create any prismatic
3-circuits.  Let $\delta$ be the curve $v_\infty P_1 I_m P_3$ separating $I_{m-1}$
from $Y$.  By Sub-lemma \ref{CHECK}, any newly created prismatic 3-circuit would
contain a point on $\delta$ connected to both $I_{m-1}$ and $Y$.  Since $Y$ and
$I_{m-1}$ are interior, $v_\infty$ does not connect to them.  Also, $P_1$ is not
connected to $Y$ as this would correspond to a pre-existing prismatic 3-circuit
$P_1I_{m}Y$, contrary to assumption.  So, the only vertices of $\delta$ that are
connected to both $I_{m-1}$ and $Y$ are $I_m$ and $P_3$, which result in the triangles
$I_mI_{m-1}Y$ and $P_3I_{m-1}Y$, hence do not correspond to newly created prismatic
3-circuits.  Therefore $Wh(I_mP_3)$ introduces no prismatic 3-circuits.

\vspace{.07in}
\begin{center}
\begin{picture}(0,0)%
\epsfig{file=fig11.pstex}%
\end{picture}%
\setlength{\unitlength}{3947sp}%
\begingroup\makeatletter\ifx\SetFigFont\undefined%
\gdef\SetFigFont#1#2#3#4#5{%
  \reset@font\fontsize{#1}{#2pt}%
  \fontfamily{#3}\fontseries{#4}\fontshape{#5}%
  \selectfont}%
\fi\endgroup%
\begin{picture}(4185,1639)(16,-1844)
\put(1605,-1179){\makebox(0,0)[lb]{\smash{{\SetFigFont{12}{14.4}{\familydefault}{\mddefault}{\updefault}{\color[rgb]{0,0,0}\small{$I_{m-1}$}}%
}}}}
\put(2851,-361){\makebox(0,0)[lb]{\smash{{\SetFigFont{12}{14.4}{\familydefault}{\mddefault}{\updefault}{\color[rgb]{0,0,0}\small{$P_1$}}%
}}}}
\put(4201,-961){\makebox(0,0)[lb]{\smash{{\SetFigFont{12}{14.4}{\familydefault}{\mddefault}{\updefault}{\color[rgb]{0,0,0}\small{$P_2$}}%
}}}}
\put(2701,-1786){\makebox(0,0)[lb]{\smash{{\SetFigFont{12}{14.4}{\familydefault}{\mddefault}{\updefault}{\color[rgb]{0,0,0}\small{$P_3$}}%
}}}}
\put(3601,-886){\makebox(0,0)[lb]{\smash{{\SetFigFont{12}{14.4}{\familydefault}{\mddefault}{\updefault}{\color[rgb]{0,0,0}\small{$I_1$}}%
}}}}
\put(3076,-886){\makebox(0,0)[lb]{\smash{{\SetFigFont{12}{14.4}{\familydefault}{\mddefault}{\updefault}{\color[rgb]{0,0,0}\small{$I_2$}}%
}}}}
\put(2776,-886){\makebox(0,0)[lb]{\smash{{\SetFigFont{12}{14.4}{\familydefault}{\mddefault}{\updefault}{\color[rgb]{0,0,0}\small{$I_3$}}%
}}}}
\put(1201,-886){\makebox(0,0)[lb]{\smash{{\SetFigFont{12}{14.4}{\familydefault}{\mddefault}{\updefault}{\color[rgb]{0,0,0}\small{$I_m$}}%
}}}}
\put(1051,-361){\makebox(0,0)[lb]{\smash{{\SetFigFont{12}{14.4}{\familydefault}{\mddefault}{\updefault}{\color[rgb]{0,0,0}\small{$X$}}%
}}}}
\put(1051,-1711){\makebox(0,0)[lb]{\smash{{\SetFigFont{12}{14.4}{\familydefault}{\mddefault}{\updefault}{\color[rgb]{0,0,0}\small{$Y$}}%
}}}}
\put(2228,-901){\makebox(0,0)[lb]{\smash{{\SetFigFont{12}{14.4}{\familydefault}{\mddefault}{\updefault}{\color[rgb]{0,0,0}\small{$I_{m-2}$}}%
}}}}
\put(151,-961){\makebox(0,0)[lb]{\smash{{\SetFigFont{12}{14.4}{\familydefault}{\mddefault}{\updefault}{\color[rgb]{0,0,0}\tiny{other vertices}}%
}}}}
\end{picture}%

\end{center}
\vspace{.07in}

Next, one must do the moves $Wh(I_{m-1}P_1)$,...,$Wh(I_1P_1)$, in that order.
(See the figure below).  We check that each of these moves creates no prismatic
3-circuits.  Fix $1 \le k \le m-1$, and let $\delta$ be the loop $v_\infty P_1
I_k P_3$.  $Wh(I_k P_1)$ creates a new edge $I_{k-1}I_m$ if $k>1$, or $P_2 I_m$
if $k=1$, the vertices of which are separated by $\delta$.  Since $I_m$ is
interior, $v_\infty$ does not connect to $I_m$.  Also, $I_m$ is no longer
connected to $P_3$.  Therefore the only points of $\delta$ that are both
connected to $I_m$ and $I_{k-1}$ are $I_k$ and $P_1$.  
Those connections form the new triangles $P_1I_mI_{k-1}$
and $I_kI_{k-1}I_m$ (when $k=1$, replace $I_{k-1}$ with $P_2$).
Hence no prismatic $3$-circuits were created, as claimed.

%The connections form the
%new triangles $P_1I_mI_{k-1}$ and $I_k I_{k-1}I_m$, hence no prismatic
%3-circuits (when $k = 1$, the above is true with $P_2$ in place of $I_{k-1}$).
%So the move $Wh(I_kP_1)$ introduces no prismatic 3-circuits, so we can the
%sequence of Whitehead moves listed above.

\vspace{.07in}
\begin{center}
\begin{picture}(0,0)%
\epsfig{file=./fig11b.pstex}%
\end{picture}%
\setlength{\unitlength}{3947sp}%
\begingroup\makeatletter\ifx\SetFigFont\undefined%
\gdef\SetFigFont#1#2#3#4#5{%
  \reset@font\fontsize{#1}{#2pt}%
  \fontfamily{#3}\fontseries{#4}\fontshape{#5}%
  \selectfont}%
\fi\endgroup%
\begin{picture}(4215,3607)(-14,-3608)
\put(2851,-2161){\makebox(0,0)[lb]{\smash{{\SetFigFont{6}{7.2}{\familydefault}{\mddefault}{\updefault}{\color[rgb]{0,0,0}\small{$P_1$}}%
}}}}
\put(2723,-2919){\makebox(0,0)[lb]{\smash{{\SetFigFont{6}{7.2}{\familydefault}{\mddefault}{\updefault}{\color[rgb]{0,0,0}\small{$I_k$}}%
}}}}
\put(1576,-2986){\makebox(0,0)[lb]{\smash{{\SetFigFont{6}{7.2}{\familydefault}{\mddefault}{\updefault}{\color[rgb]{0,0,0}\small{$I_{m-1}$}}%
}}}}
\put(151,-2761){\makebox(0,0)[lb]{\smash{{\SetFigFont{6}{7.2}{\familydefault}{\mddefault}{\updefault}{\color[rgb]{0,0,0}\tiny{other vertices}}%
}}}}
\put(1051,-3511){\makebox(0,0)[lb]{\smash{{\SetFigFont{6}{7.2}{\familydefault}{\mddefault}{\updefault}{\color[rgb]{0,0,0}\small{$Y$}}%
}}}}
\put(1051,-2161){\makebox(0,0)[lb]{\smash{{\SetFigFont{6}{7.2}{\familydefault}{\mddefault}{\updefault}{\color[rgb]{0,0,0}\small{$X$}}%
}}}}
\put(1201,-2611){\makebox(0,0)[lb]{\smash{{\SetFigFont{6}{7.2}{\familydefault}{\mddefault}{\updefault}{\color[rgb]{0,0,0}\small{$I_m$}}%
}}}}
\put(4201,-2761){\makebox(0,0)[lb]{\smash{{\SetFigFont{6}{7.2}{\familydefault}{\mddefault}{\updefault}{\color[rgb]{0,0,0}\small{$P_2$}}%
}}}}
\put(2191,-1816){\makebox(0,0)[lb]{\smash{{\SetFigFont{6}{7.2}{\familydefault}{\mddefault}{\updefault}{\color[rgb]{0,0,0}\small{$Wh(I_kP_1)$}}%
}}}}
\put(2851,-61){\makebox(0,0)[lb]{\smash{{\SetFigFont{6}{7.2}{\familydefault}{\mddefault}{\updefault}{\color[rgb]{0,0,0}\small{$P_1$}}%
}}}}
\put(2723,-819){\makebox(0,0)[lb]{\smash{{\SetFigFont{6}{7.2}{\familydefault}{\mddefault}{\updefault}{\color[rgb]{0,0,0}\small{$I_k$}}%
}}}}
\put(1576,-886){\makebox(0,0)[lb]{\smash{{\SetFigFont{6}{7.2}{\familydefault}{\mddefault}{\updefault}{\color[rgb]{0,0,0}\small{$I_{m-1}$}}%
}}}}
\put(151,-661){\makebox(0,0)[lb]{\smash{{\SetFigFont{6}{7.2}{\familydefault}{\mddefault}{\updefault}{\color[rgb]{0,0,0}\tiny{other vertices}}%
}}}}
\put(1051,-1411){\makebox(0,0)[lb]{\smash{{\SetFigFont{6}{7.2}{\familydefault}{\mddefault}{\updefault}{\color[rgb]{0,0,0}\small{$Y$}}%
}}}}
\put(1051,-61){\makebox(0,0)[lb]{\smash{{\SetFigFont{6}{7.2}{\familydefault}{\mddefault}{\updefault}{\color[rgb]{0,0,0}\small{$X$}}%
}}}}
\put(1201,-511){\makebox(0,0)[lb]{\smash{{\SetFigFont{6}{7.2}{\familydefault}{\mddefault}{\updefault}{\color[rgb]{0,0,0}\small{$I_m$}}%
}}}}
\put(4201,-661){\makebox(0,0)[lb]{\smash{{\SetFigFont{6}{7.2}{\familydefault}{\mddefault}{\updefault}{\color[rgb]{0,0,0}\small{$P_2$}}%
}}}}
\put(2701,-1486){\makebox(0,0)[lb]{\smash{{\SetFigFont{6}{7.2}{\familydefault}{\mddefault}{\updefault}{\color[rgb]{0,0,0}\small{$P_3$}}%
}}}}
\put(3076,-586){\makebox(0,0)[lb]{\smash{{\SetFigFont{6}{7.2}{\familydefault}{\mddefault}{\updefault}{\color[rgb]{0,0,0}\small{$I_{k-1}$}}%
}}}}
\put(3076,-2686){\makebox(0,0)[lb]{\smash{{\SetFigFont{6}{7.2}{\familydefault}{\mddefault}{\updefault}{\color[rgb]{0,0,0}\small{$I_{k-1}$}}%
}}}}
\put(2414,-819){\makebox(0,0)[lb]{\smash{{\SetFigFont{6}{7.2}{\familydefault}{\mddefault}{\updefault}{\color[rgb]{0,0,0}\small{$I_{k+1}$}}%
}}}}
\put(2701,-3586){\makebox(0,0)[lb]{\smash{{\SetFigFont{6}{7.2}{\familydefault}{\mddefault}{\updefault}{\color[rgb]{0,0,0}\small{$P_3$}}%
}}}}
\put(2414,-2919){\makebox(0,0)[lb]{\smash{{\SetFigFont{6}{7.2}{\familydefault}{\mddefault}{\updefault}{\color[rgb]{0,0,0}\small{$I_{k+1}$}}%
}}}}
\put(3587,-798){\makebox(0,0)[lb]{\smash{{\SetFigFont{6}{7.2}{\familydefault}{\mddefault}{\updefault}{\color[rgb]{0,0,0}\small{$I_1$}}%
}}}}
\put(3592,-2893){\makebox(0,0)[lb]{\smash{{\SetFigFont{6}{7.2}{\familydefault}{\mddefault}{\updefault}{\color[rgb]{0,0,0}\small{$I_1$}}%
}}}}
\end{picture}%

\end{center}
\vspace{.07in}

After this sequence of Whitehead moves, we obtain the diagram below,
with $I_m$ connected to $P$ exactly at $P_1$ and $P_2$.  We can then apply
Sub-lemma \ref{SUBLEMMA1} to increase the length of $P$ by the move $Wh(P_1P_2)$, as
shown below.

\begin{center}
\begin{picture}(0,0)%
\epsfig{file=./fig12.pstex}%
\end{picture}%
\setlength{\unitlength}{3947sp}%
\begingroup\makeatletter\ifx\SetFigFont\undefined%
\gdef\SetFigFont#1#2#3#4#5{%
  \reset@font\fontsize{#1}{#2pt}%
  \fontfamily{#3}\fontseries{#4}\fontshape{#5}%
  \selectfont}%
\fi\endgroup%
\begin{picture}(4185,1639)(16,-1844)
\put(3603,-1092){\makebox(0,0)[lb]{\smash{{\SetFigFont{12}{14.4}{\familydefault}{\mddefault}{\updefault}{\color[rgb]{0,0,0}\small{$I_1$}}%
}}}}
\put(2851,-361){\makebox(0,0)[lb]{\smash{{\SetFigFont{12}{14.4}{\familydefault}{\mddefault}{\updefault}{\color[rgb]{0,0,0}\small{$P_1$}}%
}}}}
\put(4201,-961){\makebox(0,0)[lb]{\smash{{\SetFigFont{12}{14.4}{\familydefault}{\mddefault}{\updefault}{\color[rgb]{0,0,0}\small{$P_2$}}%
}}}}
\put(2701,-1786){\makebox(0,0)[lb]{\smash{{\SetFigFont{12}{14.4}{\familydefault}{\mddefault}{\updefault}{\color[rgb]{0,0,0}\small{$P_3$}}%
}}}}
\put(1051,-361){\makebox(0,0)[lb]{\smash{{\SetFigFont{12}{14.4}{\familydefault}{\mddefault}{\updefault}{\color[rgb]{0,0,0}\small{$X$}}%
}}}}
\put(1051,-1711){\makebox(0,0)[lb]{\smash{{\SetFigFont{12}{14.4}{\familydefault}{\mddefault}{\updefault}{\color[rgb]{0,0,0}\small{$Y$}}%
}}}}
\put(151,-961){\makebox(0,0)[lb]{\smash{{\SetFigFont{12}{14.4}{\familydefault}{\mddefault}{\updefault}{\color[rgb]{0,0,0}\tiny{other vertices}}%
}}}}
\put(3076,-1118){\makebox(0,0)[lb]{\smash{{\SetFigFont{12}{14.4}{\familydefault}{\mddefault}{\updefault}{\color[rgb]{0,0,0}\small{$I_2$}}%
}}}}
\put(2723,-1119){\makebox(0,0)[lb]{\smash{{\SetFigFont{12}{14.4}{\familydefault}{\mddefault}{\updefault}{\color[rgb]{0,0,0}\small{$I_3$}}%
}}}}
\put(2415,-1126){\makebox(0,0)[lb]{\smash{{\SetFigFont{12}{14.4}{\familydefault}{\mddefault}{\updefault}{\color[rgb]{0,0,0}\small{$I_4$}}%
}}}}
\put(1201,-811){\makebox(0,0)[lb]{\smash{{\SetFigFont{12}{14.4}{\familydefault}{\mddefault}{\updefault}{\color[rgb]{0,0,0}\small{$I_m$}}%
}}}}
\put(1576,-1186){\makebox(0,0)[lb]{\smash{{\SetFigFont{12}{14.4}{\familydefault}{\mddefault}{\updefault}{\color[rgb]{0,0,0}\small{$I_{m-1}$}}%
}}}}
\end{picture}%

\end{center}
\vspace{.1in}

\begin{center}
\begin{picture}(0,0)%
\epsfig{file=fig13.pstex}%
\end{picture}%
\setlength{\unitlength}{3947sp}%
\begingroup\makeatletter\ifx\SetFigFont\undefined%
\gdef\SetFigFont#1#2#3#4#5{%
  \reset@font\fontsize{#1}{#2pt}%
  \fontfamily{#3}\fontseries{#4}\fontshape{#5}%
  \selectfont}%
\fi\endgroup%
\begin{picture}(4185,1639)(16,-1844)
\put(3596,-1094){\makebox(0,0)[lb]{\smash{{\SetFigFont{12}{14.4}{\familydefault}{\mddefault}{\updefault}{\color[rgb]{0,0,0}\small{$I_1$}}%
}}}}
\put(2851,-361){\makebox(0,0)[lb]{\smash{{\SetFigFont{12}{14.4}{\familydefault}{\mddefault}{\updefault}{\color[rgb]{0,0,0}\small{$P_1$}}%
}}}}
\put(2701,-1786){\makebox(0,0)[lb]{\smash{{\SetFigFont{12}{14.4}{\familydefault}{\mddefault}{\updefault}{\color[rgb]{0,0,0}\small{$P_3$}}%
}}}}
\put(1051,-361){\makebox(0,0)[lb]{\smash{{\SetFigFont{12}{14.4}{\familydefault}{\mddefault}{\updefault}{\color[rgb]{0,0,0}\small{$X$}}%
}}}}
\put(1051,-1711){\makebox(0,0)[lb]{\smash{{\SetFigFont{12}{14.4}{\familydefault}{\mddefault}{\updefault}{\color[rgb]{0,0,0}\small{$Y$}}%
}}}}
\put(151,-961){\makebox(0,0)[lb]{\smash{{\SetFigFont{12}{14.4}{\familydefault}{\mddefault}{\updefault}{\color[rgb]{0,0,0}\tiny{other vertices}}%
}}}}
\put(1576,-1186){\makebox(0,0)[lb]{\smash{{\SetFigFont{12}{14.4}{\familydefault}{\mddefault}{\updefault}{\color[rgb]{0,0,0}\small{$I_{m-1}$}}%
}}}}
\put(3076,-1118){\makebox(0,0)[lb]{\smash{{\SetFigFont{12}{14.4}{\familydefault}{\mddefault}{\updefault}{\color[rgb]{0,0,0}\small{$I_2$}}%
}}}}
\put(2723,-1119){\makebox(0,0)[lb]{\smash{{\SetFigFont{12}{14.4}{\familydefault}{\mddefault}{\updefault}{\color[rgb]{0,0,0}\small{$I_3$}}%
}}}}
\put(2415,-1126){\makebox(0,0)[lb]{\smash{{\SetFigFont{12}{14.4}{\familydefault}{\mddefault}{\updefault}{\color[rgb]{0,0,0}\small{$I_4$}}%
}}}}
\put(1201,-811){\makebox(0,0)[lb]{\smash{{\SetFigFont{12}{14.4}{\familydefault}{\mddefault}{\updefault}{\color[rgb]{0,0,0}\small{$I_m$}}%
}}}}
\put(4201,-961){\makebox(0,0)[lb]{\smash{{\SetFigFont{12}{14.4}{\familydefault}{\mddefault}{\updefault}{\color[rgb]{0,0,0}\small{$P_2$}}%
}}}}
\end{picture}%
\end{center}
\vspace{.07in}

This concludes Case 3.

\vspace{.07in}

Since $C^*$ must belong to one of these cases, we have seen that if the
length of $P$ is less than $N-3$, we can do Whitehead moves to increase it
to $N-3$ without creating prismatic $3$-circuits.  Hence we can
reduce to the case of two interior vertices, both of which must be
endpoints.  Then we can apply Sublemma \ref{SUBLEMMA2}(b) to decrease the
number of connections between one of these two interior vertices and $P$
to exactly $3$.  The result is the complex $D_N$, as shown to the right below.

\vspace{.07in}
\begin{picture}(0,0)%
\epsfig{file=fig14.pstex}%
\end{picture}%
\setlength{\unitlength}{3947sp}%
\begingroup\makeatletter\ifx\SetFigFont\undefined%
\gdef\SetFigFont#1#2#3#4#5{%
  \reset@font\fontsize{#1}{#2pt}%
  \fontfamily{#3}\fontseries{#4}\fontshape{#5}%
  \selectfont}%
\fi\endgroup%
\begin{picture}(4074,1599)(439,-1798)
\put(1189,-641){\makebox(0,0)[lb]{\smash{{\SetFigFont{6}{7.2}{\familydefault}{\mddefault}{\updefault}{\color[rgb]{0,0,0}\small{$w$}}%
}}}}
\put(3751,-1186){\makebox(0,0)[lb]{\smash{{\SetFigFont{6}{7.2}{\familydefault}{\mddefault}{\updefault}{\color[rgb]{0,0,0}\small{$v$}}%
}}}}
\put(1240,-1148){\makebox(0,0)[lb]{\smash{{\SetFigFont{6}{7.2}{\familydefault}{\mddefault}{\updefault}{\color[rgb]{0,0,0}\small{$v$}}%
}}}}
\put(3684,-661){\makebox(0,0)[lb]{\smash{{\SetFigFont{6}{7.2}{\familydefault}{\mddefault}{\updefault}{\color[rgb]{0,0,0}\small{$w$}}%
}}}}
\end{picture}%

\Endproof Lemma \ref{ROEDER}.

\vspace{.15in}

\noindent
{\large \it{Proof of Andreev's Theorem for general polyhedra}}
\vspace{.1in}

\noindent
We have seen that Andreev's Theorem holds for every simple abstract
polyhedron $C$.  Now we consider the case of $C$ having prismatic 3-circuits.
So far, the only example we have seen has been the triangular prism.  Recall
that there are some such $C$ for which $A_C = \emptyset$, so we can only hope
to prove that ${\cal P}^0_C \neq \emptyset$ when $A_C \neq \emptyset$.

\begin{lem} \label{PI3} If $A_C \neq \emptyset$, then there are points in
$A_C$ arbitrarily close to $(\pi/3,\pi/3,\cdots,\pi/3)$.  \end{lem}

\noindent {\bf Proof:} Let ${\bf a} \in A_C$ and let ${\bf a_t} = {\bf a}(1-t)
+ (\pi/3,\pi/3,\cdots,\pi/3)t$.  To see that for each $t \in [0,1)$, ${\bf a_t}
\in A_C$, we check conditions (1-5):  each coordinate is clearly $>0$
so condition (1) is satisfied.  Given edges $e_i, e_j, e_k$
meeting at a vertex we have $({\bf a_i + a_j + a_k})(1-t) + \pi t > \pi(1-t) +
\pi t = \pi$ for $t < 1$, since $({\bf a_i + a_j + a_k}) > \pi$.  So, condition
(2) is satisfied.  Similarly, given a prismatic 3-circuit intersecting edges
$e_i, e_j, e_k$ we have $({\bf a_i + a_j + a_k})(1-t) + \pi t < \pi(1-t) + \pi
t = \pi$ for $t < 1$, so condition (3) is satisfied.  Conditions (4) and (5)
are satisfied since each component of ${\bf a_t}$ is $< \pi/2$ for $t > 0$ and
since ${\bf a}$ satisfies these conditions for $t = 0$. \Endproof

\noindent
The polyhedra corresponding to these points in $A_C$, if they exist,
will have ``spiky'' vertices and think ``necks'',
wherever there is a prismatic 3-circuit.

\vspace{.05in}

We will distinguish two types of prismatic 3-circuits.  If a prismatic
3-circuit in $C^*$ separates one point from the rest of $C^*$, we will
call it a {\it truncated triangle}, otherwise we will call it an {\it
essential 3-circuit}.  The name truncated triangle comes from the fact
that such a 3-circuit in $C^*$ corresponds geometrically to the truncation
of a vertex, forming a triangular face.  We will first prove the following
sub-case: 

\begin{prop} \label{TRUNCEXIST} Let $C$ be an abstract polyhedron (with $N > 4$
faces) in which
every prismatic 3-circuit in $C^*$ is a truncated triangle.  If $A_C$ is
non-empty, then ${\cal P}_C^0$ is non-empty.  \end{prop}

We will need the following three elementary lemmas in the proof: 

\begin{lem} \label{PERPPLANE} Given three planes in $\mathbb{H}^3$ that
intersect pairwise, but which do not intersect at a point in
$\overline{\mathbb{H}^3}$, there is a fourth plane that intersects each of
these planes at right angles.  \end{lem}

\noindent {\bf Proof:}
Suppose that the three planes are given by $P_{\bf v_1},P_{\bf v_2},$ and $P_{\bf v_3}$. 
Since there is no common point of intersection in $\overline{\mathbb{H}^3}$, the
line $P_{\bf v_1} \cap P_{\bf v_2} \cap P_{\bf v_3}$ in $E^{3,1}$ is outside of the
light-cone,
so the hyperplane $(P_{\bf v_1} \cap P_{\bf v_2} \cap P_{\bf v_3})^\perp$
intersects $\mathbb{H}^3$ and hence defines a plane orthogonal to
$P_{\bf v_1},P_{\bf v_2},$ and $P_{\bf v_3}$. \Endproof

\begin{lem} \label{DECREASE_ANGLE}  
Given two halfspaces $H_1$ and $H_2$ intersecting with dihedral angle $a \in
(0,\pi/2]$ and a point $p$ in the interior of $H_1 \cap H_2$.  Let $l_1$ be the
ray from $p$ perpendicular to $\partial H_1$ and let $\widetilde{H_1}$ be the
half-space obtained by translating $\partial H_1$ a distance $\delta$ further
from $p$ perpendicularly along $l_1$.  For $\delta > 0$ and sufficiently small
$\widetilde{H_1}$ and $H_2$ intersect with dihedral angle $\alpha(\delta)$
where $\alpha$ is a decreasing continuous function of $\delta$.
\end{lem}

\noindent {\bf Proof:}
In this and the following lemma, we use the fact that in a polyhedron $P$ with
non-obtuse dihedral angles (here $P = H_1 \cap H_2$), the foot of the
perpendicular from an arbitrary interior point of $P$ to a plane containing a
face of $P$ will lie in that face (and indeed will be an interior point of that
face).  See, for example, \cite[p. 48]{VINREFL}.

We assume that $\delta$ is sufficiently small so that $\partial
\widetilde{H_1}$ and $\partial H_2$ intersect.  Let $l_2$ be the ray from $p$
perpendicular to $\partial H_2$ and let $Q$ be the plane containing $l_1$ and
$l_2$.  By construction, $Q$ intersects $\partial H_1$, $\partial H_2$ and
$\partial \widetilde{H_1}$ each perpendicularly so that $Q \cap H_1$ and  $Q
\cap H_2$ are half-planes in $Q$ intersecting with angle $a$ and $Q \cap
\widetilde{H_1}$ and $Q \cap H_2$ are half-planes in $Q$ intersecting with
angle $\alpha(\delta)$.  See the figure below:

\vspace{0.05in}
\begin{center}
\begin{picture}(0,0)%
\epsfig{file=perturbation.pstex}%
\end{picture}%
\setlength{\unitlength}{3947sp}%
\begingroup\makeatletter\ifx\SetFigFont\undefined%
\gdef\SetFigFont#1#2#3#4#5{%
  \reset@font\fontsize{#1}{#2pt}%
  \fontfamily{#3}\fontseries{#4}\fontshape{#5}%
  \selectfont}%
\fi\endgroup%
\begin{picture}(3114,2129)(1170,-6787)
\put(2020,-5829){\makebox(0,0)[lb]{\smash{{\SetFigFont{8}{9.6}{\familydefault}{\mddefault}{\updefault}{\color[rgb]{0,0,0}$R$}%
}}}}
\put(2138,-6145){\makebox(0,0)[lb]{\smash{{\SetFigFont{8}{9.6}{\familydefault}{\mddefault}{\updefault}{\color[rgb]{0,0,0}$H_2$}%
}}}}
\put(1247,-5295){\makebox(0,0)[lb]{\smash{{\SetFigFont{8}{9.6}{\familydefault}{\mddefault}{\updefault}{\color[rgb]{0,0,0}$H_1$}%
}}}}
\put(1924,-5395){\makebox(0,0)[lb]{\smash{{\SetFigFont{8}{9.6}{\familydefault}{\mddefault}{\updefault}{\color[rgb]{0,0,0}$\delta$}%
}}}}
\put(1360,-5072){\makebox(0,0)[lb]{\smash{{\SetFigFont{8}{9.6}{\familydefault}{\mddefault}{\updefault}{\color[rgb]{0,0,0}$\widetilde{H_1}$}%
}}}}
\put(2720,-5631){\makebox(0,0)[lb]{\smash{{\SetFigFont{8}{9.6}{\familydefault}{\mddefault}{\updefault}{\color[rgb]{0,0,0}$\alpha(\delta)$}%
}}}}
\put(1312,-6022){\makebox(0,0)[lb]{\smash{{\SetFigFont{8}{9.6}{\familydefault}{\mddefault}{\updefault}{\color[rgb]{0,0,0}$p$}%
}}}}
\put(1958,-4799){\makebox(0,0)[lb]{\smash{{\SetFigFont{8}{9.6}{\familydefault}{\mddefault}{\updefault}{\color[rgb]{0,0,0}$l_1$}%
}}}}
\put(2289,-6582){\makebox(0,0)[lb]{\smash{{\SetFigFont{8}{9.6}{\familydefault}{\mddefault}{\updefault}{\color[rgb]{0,0,0}$l_2$}%
}}}}
\put(2494,-5854){\makebox(0,0)[lb]{\smash{{\SetFigFont{8}{9.6}{\familydefault}{\mddefault}{\updefault}{\color[rgb]{0,0,0}$a$}%
}}}}
\end{picture}%

\end{center}
\vspace{0.05in}

Let $R$ be the compact polygon in $Q$ bounded by $\partial \widetilde{H_1}$,
$\partial H_2$, $l_1$ and $l_2$.  (In the figure above, $R$ is the shaded
region.)  Because $R$ has non-obtuse interior angles it is a parallelogram.  If
we denote the outward pointing normal vectors to $\widetilde{H_1}$ and $H_2$,
by ${\bf v}_1(\delta)$ and ${\bf v}_2$, then  Lemma \ref{AUX} gives that
$\langle {\bf v}_1(\delta), {\bf v}_2 \rangle$ is a decreasing
continuous function of $\delta$.  Since $\alpha(\delta) = -\cos \left( \langle {\bf
v}_1(\delta), {\bf v}_2 \rangle \right)$, we see that $\alpha$ is a 
decreasing continuous function of $\delta$, as well.  \Endproof

\begin{lem}\label{TRUNCATION} Given a finite volume hyperbolic polyhedron $P$
with dihedral angles in $(0,\pi/2]$ and with trivalent ideal vertices.  Suppose
that the vertices $v_1,...,v_n$ are at distinct points at infinity and the remaining
are at finite points in $\mathbb{H}^3$.  Then there exists a
polyhedron $P'$ which is combinatorially equivalent to the result of truncating
$P$ at its ideal vertices, such that the new triangular faces of $P'$ are
orthogonal to each adjacent face, and each of the remaining dihedral angles of
$P'$ lies in $(0,\pi/2]$.
\end{lem}

\noindent {\bf Proof:} 
Let $p$ be an arbitrary point in the interior of $P = \bigcap_{i=0}^N H_{i}$ and let
$l_i$ be the ray from $p$ to $\partial H_i$ that is perpendicular
to $\partial H_i$.

By Lemma \ref{DECREASE_ANGLE}, we can decrease the dihedral angles between the
face carried by $\partial H_i$ and all adjacent faces an arbitrarily small
non-zero amount by translating  $\partial H_i$ a sufficiently small distance further
from $p$.  Appropriately repeating for each $i=1,\ldots,N$, we can shift each
of the half-spaces an appropriate small distance further away from $p$ in order
to decrease all of the dihedral angles by some non-zero
amount bounded above by any given $\epsilon > 0$.

We choose $\epsilon$ sufficiently small that the sum dihedral angles at finite
vertices of $P$ remains $> \pi$  and dihedral angles between each pair of faces
remains $> 0$.  At each infinite vertex $v_i$ of $P$  the sum of dihedral
angles becomes $< \pi$ because each dihedral angle is decreased by a non-zero
amount.  Consequently, Lemma \ref{PERPPLANE} gives a fourth plane perpendicular
to each of the three planes previously meeting at $v_i$.  This resulting
polyhedron $P'$ has the same combinatorics as $P$, except that each of the
infinite vertices of $P$ is replaced by a triangular face perpendicular to its
three adjacent faces.  By construction the dihedral angles of $P'$ are in
$(0,\pi/2]$.  \Endproof

\noindent {\bf Proof of Proposition \ref{TRUNCEXIST}:}
By the hypothesis that $N>4$ from Andreev's Theorem, $C$ is not the tetrahedron
and since we have already seen that Andreev's Theorem holds for the triangular
prism, we assume that $C$ has more than 5 faces.  In this case, one can replace
all (or all but one) of the truncated triangles by single vertices to reduce
$C$ to a simple abstract polyhedron (or to a triangular prism).  The latter
case is necessary when replacing all of the truncated triangles would lead to a
tetrahedron, ($Pr_5$ is a truncated tetrahedron.) In either case, we call the
resulting abstract polyhedron $C^0$.  

\vspace{.05in} The idea of the proof quite simple.   Since Andreev's
Theorem holds for $C^0$,  we construct a polyhedron $P^0$ realizing $C^0$ with
appropriately chosen dihedral angles.  We then decrease the dihedral angles of
$P^0$, using Lemmas \ref{TRUNCATION} and \ref{PERPPLANE} to truncate vertices
of $P^0$ as they ``go past $\infty$,'' eventually obtaining a compact
polyhedron realizing $C$ with non-obtuse dihedral angles.  \vspace{.05in}

Using that $A_C \neq \emptyset$ and Lemma \ref{PI3}, choose a point $\beta \in
A_C$ so that each coordinate of $\beta$ is within $\delta$ of $\pi/3$, with
$\delta < \pi/18$.  It will be convenient to number the edges of $C$ and $C^0$
in the following way: If there is a prismatic 3-circuit in $C^0$ (i.e.,
$C^0=Pr_5$), we number these edges $1,2,$ and $3$ in $C$ and $C^0$.

Otherwise, we can use the facts that $C^0$ is trivalent and is not a
tetrahedron, to find three edges whose endpoints are six distinct vertices.
Next, we number the remaining edges common to $C$ and $C^0$ by
$4,5,\cdots,k$.  Finally, we number the edges of $C$ that were removed to
form $C^0$ by $k+1,\cdots,n$ so that the edges surrounded by prismatic
3-circuits of $C$ with smaller angle sum (given by $\beta$) come
before those surrounded by prismatic 3-circuits with larger angle
sum.

To see that the point $\gamma =
(\beta_1,\beta_2,\beta_3,\beta_4+2\delta,\beta_5+2\delta,...,\beta_k+2\delta)$
is an element of $A_{C^0}$, we check conditions (1-5).  Each of the
dihedral angles specified by $\gamma$ is in $(0,\pi/2)$ because $0 <
\beta_i + 2\delta < \pi/3 + 3\delta < \pi/3 + \pi/6 = \pi/2$.  Therefore,
condition (1) is satisfied, as well as conditions (4) and (5) because the
angles are acute.  Two of the edges labeled $4$ and higher will enter any
vertex of $C^0$ so the sum of the three dihedral angles at each vertex is
at least $4\delta$ greater than the sum of the three dihedral angles given
by $\beta$, which is $> \pi -3\delta$.  Therefore condition (2) is
satisfied.  If there is a prismatic 3-circuit in $C^0$, it crosses the
first three edges of $C^0$ and is also a prismatic 3-circuit in $C$.  
Since $\beta \in A_C$, $\beta_1+\beta_2+\beta_3 < \pi$, and it follows 
that condition (3) is satisfied by $\gamma$.

Now define $\alpha(t) = (1-t){\gamma} + t(\beta_1,...\beta_k)$.  Let
$T_0,...,T_{l-1} \in
(0,1)$ be the values of $t$ at which there is at least one vertex of $C^0$
 for which
$\alpha(t)$ gives an angle sum of $\pi$.  We also define
$T_{-1} = 0, T_l = 1$.  
We will label the vertices that have angle sum $\pi$ at
$T_i$ by $v^i_1,\cdots,v_{n(i)}^i$. Let $C^{i+1}$ be $C^{i}$ with
$v^i_1,\cdots,v_{n(i)}^i$ truncated for $i=0,\cdots,l-1$.  Hence
$C^l$ is combinatorially equivalent to $C$.

Since the number of edges increases as we move from $C^0$ toward $C$, we will
redefine $\alpha(t)$, appending $\sum_{j=1}^i n(j)$ coordinates, all constant
and equal to $\pi/2$, for values of $t$ between $T_{i-1}$ and $T_i$.

We know that Andreev's Theorem holds for $C^0$ because $C^0$ is either
simple, or the triangular prism.  So, it is sufficient to show that if
Andreev's Theorem is satisfied for $C^i$ then it is satisfied for
$C^{i+1}$, for each $i=0,\cdots,l-1$.  To do this, we must generate a
polyhedron realizing $C^{i}$ with the vertices $v^i_1,...,v_{n(i)}^i$ at
infinity and the other vertices at finite points in $\mathbb{H}^3$.  This
will be easy with our definition of $\alpha(t)$ and Corollary
\ref{INFVERTEX}.  We will then use Lemma \ref{TRUNCATION} to truncate
these vertices.   Details follow.

To use Corollary \ref{INFVERTEX}, we must check that $\alpha(t) \in A_{C^i}$
when $t \in (T_{i-1},T_i)$.  This follows directly from the definition of
$\alpha(t)$.  To check condition (1), notice that both $\beta_j$ and $\gamma_j$
are non-zero and non-obtuse, for each $j (1 \le j \le k)$, so $\alpha_j(t)$
must be as well.  To check condition (2), note first that if a vertex of $C^i$
belongs also to both $C$ and $C^0$, then it will have dihedral angle sum
greater than $\pi$ because both the corresponding sums in both $\beta$ and
$\gamma$ have that property.  If the vertex corresponds to a truncated triangle
of $C$, but not of $C^i$, then by the definition of $T_i$, the angle sum is
greater than $\pi$.  If the vertex lies on one of the truncated triangles of
$C^i$, then at least two of the incident edges will have angles equal to
$\pi/2$.  In each case, condition (2) is satisfied.

Any prismatic 3-circuit in $C^i$ is is either a prismatic
3-circuit in both $C^0$ and $C$ (the special case where $C^0$ is the
triangular prism) or is a prismatic circuit of $C^i$ which wasn't a
prismatic circuit of $C^0$. In the first case, the dihedral angle sum is
$<\pi$ because condition (3) is satisfied by both $\beta$ and $\gamma$ and
in the second case, the angle sum is $< \pi$ by definition of the $T_i$.

For each $j=1,\cdots,k$ we have $\beta_j, \gamma_j \in (0,\pi/2)$, so
$\alpha_j(t) \in (0,\pi/2)$.  However, $\alpha_j(t) = \pi/2$ for $j > k$,
corresponding to the edges of the added triangular faces.  Fortunately, a
prismatic 4-circuit cannot cross edges of these triangular faces, since it
would have to cross two edges from the same triangle, which meet at a vertex,
which is contrary to the definition of a prismatic circuit.  Thus, a
prismatic $4$-circuit can only cross edges numbered less than or equal to $k$,
each of which is assigned an acute dihedral angle, giving that condition (4) is
satisfied.

Lemma \ref{FIVE} gives that condition (5) is a consequence of conditions (3)
and (4) for $C^1,\cdots,C^l$, because they cannot be triangular prisms.
If $C^0$ happens to be the triangular prism, condition (5) holds, since
each of its edges $e_1,\cdots,e_k$ is assigned an acute dihedral angle.

\vspace{.05in}
Consider the sequence of dihedral angles $\alpha_{m,i} = \alpha(T_{i-1} +
(1-1/m)(T_i-T_{i-1}))$.  By our above analysis, $\alpha_{m,i} \in A_{C^i}$
for each $m,i$.  In fact, by definition $\alpha_{m,i}$ limits to the point
$\alpha(T_{i}) \in \partial A_{C^i}$ (as $m \rightarrow \infty$), which satisfies conditions (1-5)  to
be in $A_{C^i}$, except that the sum of the dihedral angles at each vertex
$v^i_1,\cdots,v_{n(i)}^i$ is exactly $\pi$.  We assume that Andreev's
Theorem holds for $C^i$, so by Corollary \ref{INFVERTEX}, there exists a
non-compact polyhedron $P^i$ realizing $C^i$ with each of the vertices
$v^i_1,\cdots,v_{n(i)}^i$ at infinity and the rest of the vertices at
finite points.

By Lemma \ref{TRUNCATION}, the existence of $P^{i}$ implies that there
is a polyhedron realizing $C^{i+1}$ and therefore, by Proposition
\ref{NONEMPTYIMPLIESAND}, that Andreev's Theorem holds for the
abstract polyhedron $C^{i+1}$.  Repeating this process until $i+1 = l$,
we see that Andreev's Theorem holds for $C^l$, which is our original
abstract polyhedron $C$.  \Endproof

\begin{prop}\label{FINAL} If $A_C \neq \emptyset$, then ${\cal P}_C^0 \neq
\emptyset$.  \end{prop} Combined with Proposition
\ref{NONEMPTYIMPLIESAND}, this proposition concludes the proof of Andreev's Theorem. 

\noindent {\bf Proof:}
By Proposition \ref{APRISNONEMPTY} and Proposition \ref{TRUNCEXIST} we we are
left with the case that there are $k>0$ essential 3-circuits.  We will show
that a polyhedron realizing $C$ with dihedral angles ${\bf a} \in A_C$ can be
formed by gluing together $k+1$ appropriate sub-polyhedra, each of which has
only truncated triangles.

We will work entirely within the dual complex $C^*$.  Label
the essential 3-circuits $\gamma_1,...,\gamma_k$. The idea will be to replace
$C^*$ with $k+1$ separate abstract polyhedra $C_1^*,...,C_{k+1}^*$ each of
which has no essential 3-circuits.  The $\gamma_i$ separate the sphere into
exactly $k+1$ components.  Let $C_i^*$ be the $i$-th of these components.  To
make $C_i^*$ a simplicial complex on the sphere we must fill in the holes, 
each of which is bounded by 3 edges (some $\gamma_l$).  We glue in
the following figure (the dark outer edge is $\gamma_l$).

\begin{center} 
\begin{picture}(0,0)%
\includegraphics{./truncated_triangle.pstex}%
\end{picture}%
\setlength{\unitlength}{4144sp}%
\begingroup\makeatletter\ifx\SetFigFont\undefined%
\gdef\SetFigFont#1#2#3#4#5{%
  \reset@font\fontsize{#1}{#2pt}%
  \fontfamily{#3}\fontseries{#4}\fontshape{#5}%
  \selectfont}%
\fi\endgroup%
\begin{picture}(945,676)(866,-522)
\put(1611,-113){\makebox(0,0)[lb]{\smash{{\SetFigFont{9}{10.8}{\familydefault}{\mddefault}{\updefault}{\color[rgb]{0,0,0}$\gamma_l$}%
}}}}
\put(1369,-255){\makebox(0,0)[lb]{\smash{{\SetFigFont{9}{10.8}{\familydefault}{\mddefault}{\updefault}{\color[rgb]{0,0,0}$1$}%
}}}}
\end{picture}%
\end{center}

\noindent
so that where there was an essential prismatic 3-circuit in $C^*$ there are
now two truncated triangles in distinct $C_i^*$.
Notice that none of the $C_i$ is a triangular prism, since we have divided
up $C$ along essential prismatic 3-circuits. 

In $C_i$, we will call each vertex, edge, or face obtained by such filling
in a {\it new vertex, new edge,} or {\it new face} respectively.  We will
call all of the other edges {\it old edges}.  We label each such new
vertex with the number $l$ corresponding to the 3-circuit $\gamma_l$ that
was filled in.  For each $l$, there will be exactly 2 new vertices
labeled $l$ which are in two different $C_i^*$, $C_j^*$.
See the following diagram for an example.

\vspace{.07in} \begin{center} 
\begin{picture}(0,0)%
\epsfig{file=compound.pstex}%
\end{picture}%
\setlength{\unitlength}{4144sp}%
\begingroup\makeatletter\ifx\SetFigFont\undefined%
\gdef\SetFigFont#1#2#3#4#5{%
  \reset@font\fontsize{#1}{#2pt}%
  \fontfamily{#3}\fontseries{#4}\fontshape{#5}%
  \selectfont}%
\fi\endgroup%
\begin{picture}(3621,3518)(574,-2892)
\put(3923,-1854){\makebox(0,0)[lb]{\smash{{\SetFigFont{9}{10.8}{\familydefault}{\mddefault}{\updefault}{\color[rgb]{0,0,0}$\gamma_1$}%
}}}}
\put(1872,-12){\makebox(0,0)[lb]{\smash{{\SetFigFont{9}{10.8}{\familydefault}{\mddefault}{\updefault}{\color[rgb]{0,0,0}$\gamma_1$}%
}}}}
\put(1628,-1819){\makebox(0,0)[lb]{\smash{{\SetFigFont{9}{10.8}{\familydefault}{\mddefault}{\updefault}{\color[rgb]{0,0,0}$1$}%
}}}}
\put(1733,-950){\makebox(0,0)[lb]{\smash{{\SetFigFont{9}{10.8}{\familydefault}{\mddefault}{\updefault}{\color[rgb]{0,0,0}$C^*$}%
}}}}
\put(1594,-2861){\makebox(0,0)[lb]{\smash{{\SetFigFont{9}{10.8}{\familydefault}{\mddefault}{\updefault}{\color[rgb]{0,0,0}$C_1^*$}%
}}}}
\put(3471,-2687){\makebox(0,0)[lb]{\smash{{\SetFigFont{9}{10.8}{\familydefault}{\mddefault}{\updefault}{\color[rgb]{0,0,0}$C_2^*$}%
}}}}
\put(3367,-1124){\makebox(0,0)[lb]{\smash{{\SetFigFont{9}{10.8}{\familydefault}{\mddefault}{\updefault}{\color[rgb]{0,0,0}$1$ (at infinity)}%
}}}}
\put(1802,-1714){\makebox(0,0)[lb]{\smash{{\SetFigFont{9}{10.8}{\familydefault}{\mddefault}{\updefault}{\color[rgb]{0,0,0}$\gamma_1$}%
}}}}
\end{picture}%

\end{center}
\vspace{.07in}

The choice of angles ${\bf a} \in A_C$ gives dihedral angles assigned to
each old edge in each $C_i^*$.  Assign to each of the new edges $\pi/2$. 
This gives a choice of angles ${\bf a_i}$ for for which it is easy to check that
${\bf a_i} \in A_{C_i}$ for each $i$:

Clearly condition (1) is satisfied since these angles are non-zero and none of
them obtuse.  The angles along each triangle of old edges in $C_i^*$ already
satisfy condition (2) since ${\bf a} \in A_C$.  For each of the new triangles
added, two of the edges are assigned $\pi/2$ and the third was already assigned
a non-zero angle, according to ${\bf a}$, so condition (2) is satisfied for
these triangles, too.  None of the new edges in $C_i^*$ can be in a prismatic
3-circuit or a prismatic 4-circuit since such a circuit would have to involve
two such new edges, which form two sides of a triangle.  Therefore, each
prismatic 3- or 4-circuit in $C_i^*$ has come from such a circuit in $C^*$, so
the choice of angles ${\bf a_i}$ will satisfy (3) and (4).  Since none of the
$C_i^*$ is a triangular prism, condition (5) is a consequence of condition (4),
and hence is satisfied. 

Therefore by Proposition \ref{TRUNCEXIST} there exist polyhedra $P_i$ realizing
the data $(C_i,{\bf a_i})$, for $1 \le i \le k+1$.  For each pair of new
vertices labeled $l$, the two faces dual to them are isomorphic, since by
Proposition \ref{TRIVALENT}, the face angles are the same (giving congruent
triangular faces).  So one can glue all of the $P_i$ together according to
the labeling by $l$.  Since the edges of these triangles were assigned dihedral
angles of $\pi/2$, the faces coming together from opposite sides of such a
glued pair fit together smoothly. The result is a polyhedron $P$ realizing $C$
and angles ${\bf a}$.  See the following diagram. 

\vspace{0.08in}
\begin{center} 
\begin{picture}(0,0)%
\epsfig{file=./gluing.pstex}%
\end{picture}%
\setlength{\unitlength}{4144sp}%
\begingroup\makeatletter\ifx\SetFigFont\undefined%
\gdef\SetFigFont#1#2#3#4#5{%
  \reset@font\fontsize{#1}{#2pt}%
  \fontfamily{#3}\fontseries{#4}\fontshape{#5}%
  \selectfont}%
\fi\endgroup%
\begin{picture}(4812,835)(361,-1239)
\put(361,-770){\makebox(0,0)[lb]{\smash{{\SetFigFont{6}{7.2}{\familydefault}{\mddefault}{\updefault}{\color[rgb]{0,0,0}The rest of $P_1$}%
}}}}
\put(4288,-872){\makebox(0,0)[lb]{\smash{{\SetFigFont{6}{7.2}{\familydefault}{\mddefault}{\updefault}{\color[rgb]{0,0,0}1}%
}}}}
\put(4509,-827){\makebox(0,0)[lb]{\smash{{\SetFigFont{6}{7.2}{\familydefault}{\mddefault}{\updefault}{\color[rgb]{0,0,0}The rest of $P_2$}%
}}}}
\put(3314,-776){\makebox(0,0)[lb]{\smash{{\SetFigFont{6}{7.2}{\familydefault}{\mddefault}{\updefault}{\color[rgb]{0,0,0}The rest of $P_1$}%
}}}}
\put(2058,-864){\makebox(0,0)[lb]{\smash{{\SetFigFont{6}{7.2}{\familydefault}{\mddefault}{\updefault}{\color[rgb]{0,0,0}The rest of $P_2$}%
}}}}
\put(1360,-887){\makebox(0,0)[lb]{\smash{{\SetFigFont{6}{7.2}{\familydefault}{\mddefault}{\updefault}{\color[rgb]{0,0,0}1}%
}}}}
\put(1822,-910){\makebox(0,0)[lb]{\smash{{\SetFigFont{6}{7.2}{\familydefault}{\mddefault}{\updefault}{\color[rgb]{0,0,0}1}%
}}}}
\end{picture}%

\end{center}
\Endproof

\noindent
That concludes the proof that $\alpha: {\cal P}_C^0 \rightarrow A_C$ is a 
homeomorphism for every abstract polyhedron $C$ having more than four faces
and hence concludes the proof of Andreev's Theorem.

\section{Example of the combinatorial algorithm from Lemma \ref{ROEDER}}
We include an example of the combinatorial algorithm described in Lemma
\ref{ROEDER}, which gives a sequence of Whitehead moves to reduce the dual
complex of a simple abstract polyhedron, $C^*$, to $D_N^*$.  We then explain
how the sequence of Whitehead moves described in Andreev's paper \cite{AND}
would result in prismatic 3-circuits for this $C$ and for many others. 

\begin{center}
\begin{picture}(0,0)%
\epsfig{file=wh_seq4a.pstex}%
\end{picture}%
\setlength{\unitlength}{4144sp}%
\begingroup\makeatletter\ifx\SetFigFont\undefined%
\gdef\SetFigFont#1#2#3#4#5{%
  \reset@font\fontsize{#1}{#2pt}%
  \fontfamily{#3}\fontseries{#4}\fontshape{#5}%
  \selectfont}%
\fi\endgroup%
\begin{picture}(5370,4167)(994,-3872)
\put(1068,-3818){\makebox(0,0)[lb]{\smash{{\SetFigFont{12}{14.4}{\familydefault}{\mddefault}{\updefault}{\color[rgb]{0,0,0}{\small Case 3 from the proof of Lemma \ref{ROEDER} is done in subfigures (1)-(7).}}%
}}}}
\put(3490,-3455){\makebox(0,0)[lb]{\smash{{\SetFigFont{12}{14.4}{\familydefault}{\mddefault}{\updefault}{\color[rgb]{0,0,0}{\small (5)}}%
}}}}
\put(5381,-3455){\makebox(0,0)[lb]{\smash{{\SetFigFont{12}{14.4}{\familydefault}{\mddefault}{\updefault}{\color[rgb]{0,0,0}{\small (6)}}%
}}}}
\put(1571,-1460){\makebox(0,0)[lb]{\smash{{\SetFigFont{12}{14.4}{\familydefault}{\mddefault}{\updefault}{\color[rgb]{0,0,0}{\small (1)}}%
}}}}
\put(3464,-1429){\makebox(0,0)[lb]{\smash{{\SetFigFont{12}{14.4}{\familydefault}{\mddefault}{\updefault}{\color[rgb]{0,0,0}{\small (2)}}%
}}}}
\put(5446,-1474){\makebox(0,0)[lb]{\smash{{\SetFigFont{12}{14.4}{\familydefault}{\mddefault}{\updefault}{\color[rgb]{0,0,0}{\small (3)}}%
}}}}
\put(1582,-3470){\makebox(0,0)[lb]{\smash{{\SetFigFont{12}{14.4}{\familydefault}{\mddefault}{\updefault}{\color[rgb]{0,0,0}{\small (4)}}%
}}}}
\end{picture}%

\begin{picture}(0,0)%
\epsfig{file=./wh_seq4b.pstex}%
\end{picture}%
\setlength{\unitlength}{4144sp}%
\begingroup\makeatletter\ifx\SetFigFont\undefined%
\gdef\SetFigFont#1#2#3#4#5{%
  \reset@font\fontsize{#1}{#2pt}%
  \fontfamily{#3}\fontseries{#4}\fontshape{#5}%
  \selectfont}%
\fi\endgroup%
\begin{picture}(5409,4173)(129,-3730)
\put(226,299){\makebox(0,0)[lb]{\smash{{\SetFigFont{12}{14.4}{\familydefault}{\mddefault}{\updefault}{\color[rgb]{0,0,0}{\small Case 1 follows in subfigures (7)-(9) and again in subfigures(9)-(12).}}%
}}}}
\put(727,-1649){\makebox(0,0)[lb]{\smash{{\SetFigFont{12}{14.4}{\familydefault}{\mddefault}{\updefault}{\color[rgb]{0,0,0}{\small (7)}}%
}}}}
\put(2679,-1635){\makebox(0,0)[lb]{\smash{{\SetFigFont{12}{14.4}{\familydefault}{\mddefault}{\updefault}{\color[rgb]{0,0,0}{\small (8)}}%
}}}}
\put(4556,-1635){\makebox(0,0)[lb]{\smash{{\SetFigFont{12}{14.4}{\familydefault}{\mddefault}{\updefault}{\color[rgb]{0,0,0}{\small (9)}}%
}}}}
\put(752,-3661){\makebox(0,0)[lb]{\smash{{\SetFigFont{12}{14.4}{\familydefault}{\mddefault}{\updefault}{\color[rgb]{0,0,0}{\small (10)}}%
}}}}
\put(2660,-3676){\makebox(0,0)[lb]{\smash{{\SetFigFont{12}{14.4}{\familydefault}{\mddefault}{\updefault}{\color[rgb]{0,0,0}{\small (11)}}%
}}}}
\put(4612,-3646){\makebox(0,0)[lb]{\smash{{\SetFigFont{12}{14.4}{\familydefault}{\mddefault}{\updefault}{\color[rgb]{0,0,0}{\small (12)}}%
}}}}
\end{picture}%

\begin{picture}(0,0)%
\epsfig{file=wh_seq4c.pstex}%
\end{picture}%
\setlength{\unitlength}{4144sp}%
\begingroup\makeatletter\ifx\SetFigFont\undefined%
\gdef\SetFigFont#1#2#3#4#5{%
  \reset@font\fontsize{#1}{#2pt}%
  \fontfamily{#3}\fontseries{#4}\fontshape{#5}%
  \selectfont}%
\fi\endgroup%
\begin{picture}(5651,4067)(959,-3233)
\put(959,690){\makebox(0,0)[lb]{\smash{{\SetFigFont{12}{14.4}{\familydefault}{\mddefault}{\updefault}{\color[rgb]{0,0,0}{\small Case 1 is repeated three times, first in subfigures (12)-(15), then in (15)-(17)}}%
}}}}
\put(3601,-3164){\makebox(0,0)[lb]{\smash{{\SetFigFont{12}{14.4}{\familydefault}{\mddefault}{\updefault}{\color[rgb]{0,0,0}{\small (17)}}%
}}}}
\put(1595,-1246){\makebox(0,0)[lb]{\smash{{\SetFigFont{12}{14.4}{\familydefault}{\mddefault}{\updefault}{\color[rgb]{0,0,0}{\small (13)}}%
}}}}
\put(3542,-1231){\makebox(0,0)[lb]{\smash{{\SetFigFont{12}{14.4}{\familydefault}{\mddefault}{\updefault}{\color[rgb]{0,0,0}{\small (14)}}%
}}}}
\put(5636,-1231){\makebox(0,0)[lb]{\smash{{\SetFigFont{12}{14.4}{\familydefault}{\mddefault}{\updefault}{\color[rgb]{0,0,0}{\small (15)}}%
}}}}
\put(1668,-3179){\makebox(0,0)[lb]{\smash{{\SetFigFont{12}{14.4}{\familydefault}{\mddefault}{\updefault}{\color[rgb]{0,0,0}{\small (16)}}%
}}}}
\put(5680,-3164){\makebox(0,0)[lb]{\smash{{\SetFigFont{12}{14.4}{\familydefault}{\mddefault}{\updefault}{\color[rgb]{0,0,0}{\small (18)}}%
}}}}
\end{picture}%

\begin{picture}(0,0)%
\epsfig{file=wh_seq4d.pstex}%
\end{picture}%
\setlength{\unitlength}{4144sp}%
\begingroup\makeatletter\ifx\SetFigFont\undefined%
\gdef\SetFigFont#1#2#3#4#5{%
  \reset@font\fontsize{#1}{#2pt}%
  \fontfamily{#3}\fontseries{#4}\fontshape{#5}%
  \selectfont}%
\fi\endgroup%
\begin{picture}(5170,4471)(218,-3743)
\put(421,-3689){\makebox(0,0)[lb]{\smash{{\SetFigFont{12}{14.4}{\familydefault}{\mddefault}{\updefault}{\color[rgb]{0,0,0}{\small subfigures (21) and (22), then Case 1 is done in (22)-(29)}.}%
}}}}
\put(314,584){\makebox(0,0)[lb]{\smash{{\SetFigFont{12}{14.4}{\familydefault}{\mddefault}{\updefault}{\color[rgb]{0,0,0}{\small and finally in (17)-(21).  The diagram is straightened out between}}%
}}}}
\put(840,-1409){\makebox(0,0)[lb]{\smash{{\SetFigFont{12}{14.4}{\familydefault}{\mddefault}{\updefault}{\color[rgb]{0,0,0}{\small (19)}}%
}}}}
\put(840,-3283){\makebox(0,0)[lb]{\smash{{\SetFigFont{12}{14.4}{\familydefault}{\mddefault}{\updefault}{\color[rgb]{0,0,0}{\small (22)}}%
}}}}
\put(2056,-3283){\makebox(0,0)[lb]{\smash{{\SetFigFont{12}{14.4}{\familydefault}{\mddefault}{\updefault}{\color[rgb]{0,0,0}{\small (23)}}%
}}}}
\put(3315,-3269){\makebox(0,0)[lb]{\smash{{\SetFigFont{12}{14.4}{\familydefault}{\mddefault}{\updefault}{\color[rgb]{0,0,0}{\small (24)}}%
}}}}
\put(4575,-3269){\makebox(0,0)[lb]{\smash{{\SetFigFont{12}{14.4}{\familydefault}{\mddefault}{\updefault}{\color[rgb]{0,0,0}{\small (25)}}%
}}}}
\put(4273,-1380){\makebox(0,0)[lb]{\smash{{\SetFigFont{12}{14.4}{\familydefault}{\mddefault}{\updefault}{\color[rgb]{0,0,0}{\small (21)}}%
}}}}
\put(2369,-1381){\makebox(0,0)[lb]{\smash{{\SetFigFont{12}{14.4}{\familydefault}{\mddefault}{\updefault}{\color[rgb]{0,0,0}{\small (20)}}%
}}}}
\end{picture}%

\begin{picture}(0,0)%
\epsfig{file=./wh_seq4e.pstex}%
\end{picture}%
\setlength{\unitlength}{4144sp}%
\begingroup\makeatletter\ifx\SetFigFont\undefined%
\gdef\SetFigFont#1#2#3#4#5{%
  \reset@font\fontsize{#1}{#2pt}%
  \fontfamily{#3}\fontseries{#4}\fontshape{#5}%
  \selectfont}%
\fi\endgroup%
\begin{picture}(4568,4708)(1389,-4454)
\put(1456,-3946){\makebox(0,0)[lb]{\smash{{\SetFigFont{12}{14.4}{\familydefault}{\mddefault}{\updefault}{\color[rgb]{0,0,0}{\small Apply Sublemma 6.8(b) in (30)-(33) so that one of the two}}%
}}}}
\put(5478,-3702){\makebox(0,0)[lb]{\smash{{\SetFigFont{12}{14.4}{\familydefault}{\mddefault}{\updefault}{\color[rgb]{0,0,0}{\small (33)}}%
}}}}
\put(1456,-4171){\makebox(0,0)[lb]{\smash{{\SetFigFont{12}{14.4}{\familydefault}{\mddefault}{\updefault}{\color[rgb]{0,0,0}{\small interior vertices is only connected to three points on the outer}}%
}}}}
\put(1456,-4396){\makebox(0,0)[lb]{\smash{{\SetFigFont{12}{14.4}{\familydefault}{\mddefault}{\updefault}{\color[rgb]{0,0,0}{\small polygon.  This reduces the complex to $D_{18}^*$}}%
}}}}
\put(1741,-1846){\makebox(0,0)[lb]{\smash{{\SetFigFont{12}{14.4}{\familydefault}{\mddefault}{\updefault}{\color[rgb]{0,0,0}{\small (26)}}%
}}}}
\put(3046,-1816){\makebox(0,0)[lb]{\smash{{\SetFigFont{12}{14.4}{\familydefault}{\mddefault}{\updefault}{\color[rgb]{0,0,0}{\small (27)}}%
}}}}
\put(4216,-1861){\makebox(0,0)[lb]{\smash{{\SetFigFont{12}{14.4}{\familydefault}{\mddefault}{\updefault}{\color[rgb]{0,0,0}{\small (28)}}%
}}}}
\put(5446,-1831){\makebox(0,0)[lb]{\smash{{\SetFigFont{12}{14.4}{\familydefault}{\mddefault}{\updefault}{\color[rgb]{0,0,0}{\small (29)}}%
}}}}
\put(1396,-2161){\makebox(0,0)[lb]{\smash{{\SetFigFont{12}{14.4}{\familydefault}{\mddefault}{\updefault}{\color[rgb]{0,0,0}{\small The diagram is striaghtened out between subfigures (29) and (30).}}%
}}}}
\put(1389,110){\makebox(0,0)[lb]{\smash{{\SetFigFont{12}{14.4}{\familydefault}{\mddefault}{\updefault}{\color[rgb]{0,0,0}{\small Subfigures (22)-(29) are another instance of Case 1.}}%
}}}}
\put(1668,-3672){\makebox(0,0)[lb]{\smash{{\SetFigFont{12}{14.4}{\familydefault}{\mddefault}{\updefault}{\color[rgb]{0,0,0}{\small (30)}}%
}}}}
\put(3018,-3687){\makebox(0,0)[lb]{\smash{{\SetFigFont{12}{14.4}{\familydefault}{\mddefault}{\updefault}{\color[rgb]{0,0,0}{\small (31)}}%
}}}}
\put(4188,-3687){\makebox(0,0)[lb]{\smash{{\SetFigFont{12}{14.4}{\familydefault}{\mddefault}{\updefault}{\color[rgb]{0,0,0}{\small (32)}}%
}}}}
\end{picture}%

\end{center}

It is interesting to note that Andreev's version of our Lemma \ref{ROEDER} (his
Theorem 6 in \cite{AND}) would fail for this abstract polyhedron $C^*$.  The
major difficulty is to achieve the first increase in the length of the outer
polygon $P$.  We carefully chose the vertex $v$ where the graph of interior
vertices and edges branches and then did Whitehead moves to reduce the number
of components where this vertex is connected to $P$.  This is done in
sub-figures (1)-(4).  If we had started with any other interior vertex and
tried to decrease the number of components where it is connected to $P$,
prismatic 3-circuits would develop as shown below (the dashed curve) in
sub-figures (1) and (2). 

\begin{center}
\vspace{.07in}
\begin{picture}(0,0)%
\includegraphics{wh_seq2d.pstex}%
\end{picture}%
\setlength{\unitlength}{4144sp}%
\begingroup\makeatletter\ifx\SetFigFont\undefined%
\gdef\SetFigFont#1#2#3#4#5{%
  \reset@font\fontsize{#1}{#2pt}%
  \fontfamily{#3}\fontseries{#4}\fontshape{#5}%
  \selectfont}%
\fi\endgroup%
\begin{picture}(3410,1809)(998,-1514)
\put(1571,-1460){\makebox(0,0)[lb]{\smash{\SetFigFont{12}{14.4}{\familydefault}{\mddefault}{\updefault}{\color[rgb]{0,0,0}{\small (1)}}%
}}}
\put(3464,-1429){\makebox(0,0)[lb]{\smash{\SetFigFont{12}{14.4}{\familydefault}{\mddefault}{\updefault}{\color[rgb]{0,0,0}{\small (2)}}%
}}}
\end{picture}

\end{center}

\noindent
Andreev states
that one must do Whitehead moves until each interior vertex is connected to
$P$ in a single component (creating what he calls the inner polygon), but 
he {\it does not} indicate that one must do this for certain interior vertices
before others.  Using a very specific order of Whitehead moves to reduce the 
number of components of connection from each interior vertex to $P$ to be
either one or zero would work, but Andreev does not prove this.  Instead 
of doing this, we find that it is simpler just to do Whitehead moves so that
$v$ is connected to $P$ in a single component consisting of two vertices and
an edge, instead of creating the whole ``inner polygon'' as Andreev would.
Once this is done, doing the Whitehead move on this edge of $P$ increases
the length of $P$ by one.

There are cases where the graph of interior vertices and edges has no branching 
points and Andreev's proof could not work, even having chosen to do
Whitehead moves to decrease the components of connection between each interior
vertex and $P$ to one in some specific order.  This is true for the following
abstract polyhedron $C^*$, shown in sub-figure (1) below:

\begin{center}
\vspace{.07in}
\begin{picture}(0,0)%
\includegraphics{wh_seq3a.pstex}%
\end{picture}%
\setlength{\unitlength}{4144sp}%
\begingroup\makeatletter\ifx\SetFigFont\undefined%
\gdef\SetFigFont#1#2#3#4#5{%
  \reset@font\fontsize{#1}{#2pt}%
  \fontfamily{#3}\fontseries{#4}\fontshape{#5}%
  \selectfont}%
\fi\endgroup%
\begin{picture}(3552,972)(1426,-928)
\put(1812,-900){\makebox(0,0)[lb]{\smash{\SetFigFont{7}{8.4}{\familydefault}{\mddefault}{\updefault}{\color[rgb]{0,0,0}{\small (1)}}%
}}}
\put(3896,-874){\makebox(0,0)[lb]{\smash{\SetFigFont{7}{8.4}{\familydefault}{\mddefault}{\updefault}{\color[rgb]{0,0,0}{\small (2)}}%
}}}
\end{picture}

\end{center}

\noindent
Each interior vertex of $C^*$ is either an endpoint, or connected to $P$
in exactly two components, each of which is a single point.  Doing a Whitehead
move to eliminate any of these connections would result in a prismatic 3-circuit
like the dashed one in sub-figure (2) above.

Instead, Case 2 of our Lemma \ref{ROEDER} ``borrows'' a point of connection
from one of the endpoints (sub-figures (1) and (2) below), making one of the interior
vertices connected to $P$ in two components, one consisting of two vertices and
an edge, and the other consisting of one vertex.  One then can eliminate the
single point of connection as in the Whitehead move (2) to (3).  Then, it is
simple to increase the length of $P$ by doing the Whitehead move on the single
edge in the single component of connection between this interior vertex and
$P$, as shown in the change from sub-figure (3) to (4).

\begin{center}
\vspace{.07in}
\begin{picture}(0,0)%
\includegraphics{wh_seq3b.pstex}%
\end{picture}%
\setlength{\unitlength}{4144sp}%
\begingroup\makeatletter\ifx\SetFigFont\undefined%
\gdef\SetFigFont#1#2#3#4#5{%
  \reset@font\fontsize{#1}{#2pt}%
  \fontfamily{#3}\fontseries{#4}\fontshape{#5}%
  \selectfont}%
\fi\endgroup%
\begin{picture}(3322,1918)(1416,-1874)
\put(3720,-812){\makebox(0,0)[lb]{\smash{\SetFigFont{7}{8.4}{\familydefault}{\mddefault}{\updefault}{\color[rgb]{0,0,0}{\small (2)}}%
}}}
\put(1776,-1848){\makebox(0,0)[lb]{\smash{\SetFigFont{7}{8.4}{\familydefault}{\mddefault}{\updefault}{\color[rgb]{0,0,0}{\small (3)}}%
}}}
\put(1783,-836){\makebox(0,0)[lb]{\smash{\SetFigFont{7}{8.4}{\familydefault}{\mddefault}{\updefault}{\color[rgb]{0,0,0}{\small (1)}}%
}}}
\put(3732,-1847){\makebox(0,0)[lb]{\smash{\SetFigFont{7}{8.4}{\familydefault}{\mddefault}{\updefault}{\color[rgb]{0,0,0}{\small (4)}}%
}}}
\end{picture}

\end{center}

\vspace{.2in}
\noindent
{\bf \Large Acknowledgments}\\
\noindent
This paper is a abridged version of the thesis written by the first author
for the Universit\'e de Provence \cite{ROE}.
He would like to thank first of all John Hamal Hubbard, his thesis director, for
getting him interested and involved in the beautiful subject of hyperbolic
geometry.  

In 2002, J. H. Hubbard suggested that the first author write a
computer program that implements Andreev's Theorem using the ideas from
Andreev's proof.  He was given a version of this proof written by J. H. Hubbard
and William D. Dunbar as a starting point.  From the ideas in that manuscript,
he wrote a program which, while faithfully implementing Andreev's algorithm,
was unsuccessful in computing the polyhedra with given dihedral angles and
combinatorial structure.  From there, he found the error in Andreev's proof,
which led him to first to a way to correct the program so that the polyhedra
for which it was tested were actually computed correctly, and later to the
arguments needed close the gap in Andreev's proof.

All of the work that the first author has done has developed out of ideas that
he learned from the manuscript written by J. H. Hubbard and W. D. Dunbar.
Therefore, he thanks both J. H. Hubbard and Bill Dunbar for their extensive
work on that preliminary manuscript, and for their significant subsequent help
and interest in his corrections of the proof by including them as co-authors.

The authors would also like to thank Adrien Douady, Michel Boileau, Jean-Pierre Otal,
Hamish Short, and J\'er\^ome Los for their helpful comments and their interest
in our work.

\bibliographystyle{plain}
\bibliography{andreev.bib}

\begin{thebibliography}{10}

\bibitem{GEO}
www.geomview.org, Developed by The Geometry Center at the University of
  Minnesota in the late 1990's.

\bibitem{PFB}
Martin Aigner and G{\"u}nter~M. Ziegler.
\newblock {\em Proofs from {T}he {B}ook}.
\newblock Springer-Verlag, Berlin, third edition, 2004.
\newblock Including illustrations by Karl H. Hofmann.

\bibitem{AVS}
D.~V. Alekseevskij, {\`E}.~B. Vinberg, and A.~S. Solodovnikov.
\newblock Geometry of spaces of constant curvature.
\newblock In {\em Geometry, II}, volume~29 of {\em Encyclopaedia Math. Sci.},
  pages 1--138. Springer, Berlin, 1993.

\bibitem{AND}
E.~M. Andreev.
\newblock On convex polyhedra in {L}obacevskii spaces ({E}nglish
  {T}ranslation).
\newblock {\em Math. USSR Sbornik}, 10:413--440, 1970.

\bibitem{AND2}
E.~M. Andreev.
\newblock On convex polyhedra in {L}obacevskii spaces (in {R}ussian).
\newblock {\em Mat. Sb.}, 81(123):445--478, 1970.

\bibitem{BAO}
Xiliang Bao and Francis Bonahon.
\newblock Hyperideal polyhedra in hyperbolic 3-space.
\newblock {\em Bull. Soc. Math. France}, 130(3):457--491, 2002.

\bibitem{BOI}
Michel Boileau.
\newblock Uniformisation en dimension trois.
\newblock {\em Ast\'erisque}, (266):Exp.\ No.\ 855, 4, 137--174, 2000.
\newblock S\'eminaire Bourbaki, Vol. 1998/99.

\bibitem{BP}
Michel Boileau and Joan Porti.
\newblock Geometrization of 3-orbifolds of cyclic type.
\newblock {\em Ast\'erisque}, (272):208, 2001.
\newblock Appendix A by Michael Heusener and Porti.

\bibitem{STEVE}
Phil Bowers and Kenneth Stephenson.
\newblock A branched {A}ndreev-{T}hurston theorem for circle packings of the
  sphere.
\newblock {\em Proc. London Math. Soc. (3)}, 73(1):185--215, 1996.

\bibitem{LUO}
Bennett Chow and Feng Luo.
\newblock Combinatorial {R}icci flows on surfaces.
\newblock {\em J. Differential Geom.}, 63(1):97--129, 2003.

\bibitem{CHK}
Daryl Cooper, Craig~D. Hodgson, and Steven~P. Kerckhoff.
\newblock {\em Three-dimensional orbifolds and cone-manifolds}, volume~5 of
  {\em MSJ Memoirs}.
\newblock Mathematical Society of Japan, Tokyo, 2000.
\newblock With a postface by Sadayoshi Kojima.

\bibitem{DIAZ}
Raquel D{\'{\i}}az.
\newblock Non-convexity of the space of dihedral angles of hyperbolic
  polyhedra.
\newblock {\em C. R. Acad. Sci. Paris S\'er. I Math.}, 325(9):993--998, 1997.

\bibitem{DIAZ_ANDREEV}
Raquel D{\'{\i}}az.
\newblock A generalization of {A}ndreev's theorem.
\newblock {\em J. Math. Soc. Japan}, 58(2):333--349, 2006.

\bibitem{DOUADY}
R{\'e}gine Douady and Adrien Douady.
\newblock {\em Alg\`ebre et th\'eories galoisiennes. 2}.
\newblock CEDIC, Paris, 1979.

\bibitem{GUE}
Fran{\c{c}}ois Gu{\'e}ritaud.
\newblock On an elementary proof of {R}ivin's characterization of convex ideal
  hyperbolic polyhedra by their dihedral angles.
\newblock {\em Geom. Dedicata}, 108:111--124, 2004.

\bibitem{H}
C.~D. Hodgson.
\newblock Deduction of {A}ndreev's theorem from {R}ivin's characterization of
  convex hyperbolic polyhedra.
\newblock In {\em Topology 90}, pages 185--193. de Gruyter, 1992.

\bibitem{KAP}
Michael Kapovich.
\newblock {\em Hyperbolic manifolds and discrete groups}, volume 183 of {\em
  Progress in Mathematics}.
\newblock Birkh\"auser Boston Inc., Boston, MA, 2001.

\bibitem{LIMA}
Elon~Lages Lima.
\newblock {\em Fundamental groups and covering spaces}.
\newblock A K Peters Ltd., Natick, MA, 2003.
\newblock Translated from the Portuguese by Jonas Gomes.

\bibitem{MR}
A.~Marden and B.~Rodin.
\newblock On {T}hurston's formulation and proof of {A}ndreev's {T}heorem.
\newblock In {\em Computational Methods and Function Theory}, volume 1435 of
  {\em Lecture Notes in Mathematics}, pages 103--115. Springer-Verlag, 1990.

\bibitem{OT}
Jean-Pierre Otal.
\newblock Thurston's hyperbolization of {H}aken manifolds.
\newblock In {\em Surveys in differential geometry, Vol. III (Cambridge, MA,
  1996)}, pages 77--194. Int. Press, 1998.

\bibitem{RH}
I.~Rivin and C.~D. Hodgson.
\newblock A characterization of compact convex polyhedra in hyperbolic 3-space.
\newblock {\em Invent. Math.}, 111:77--111, 1993.

\bibitem{RIV_IDEAL1}
Igor Rivin.
\newblock On geometry of convex ideal polyhedra in hyperbolic {$3$}-space.
\newblock {\em Topology}, 32(1):87--92, 1993.

\bibitem{RIV_IDEAL2}
Igor Rivin.
\newblock A characterization of ideal polyhedra in hyperbolic {$3$}-space.
\newblock {\em Ann. of Math. (2)}, 143(1):51--70, 1996.

\bibitem{RIV_IDEAL3}
Igor Rivin.
\newblock Combinatorial optimization in geometry.
\newblock {\em Adv. in Appl. Math.}, 31(1):242--271, 2003.

\bibitem{ROE}
Roland K.~W. Roeder.
\newblock Le th\'eor\`eme d'andreev sur poly\`edres hyperboliques.
\newblock Doctoral thesis (In English), May 2004.
\newblock Universit\'e de Provence, Aix-Marseille 1.

\bibitem{ROE_TET}
Roland K.~W. Roeder.
\newblock Compact hyperbolic tetrahedra with non-obtuse dihedral angles.
\newblock {\em Publicacions Matem\`atiques}, 50(1):211--227, 2006.

\bibitem{SCH2}
J.-M. Schlenker.
\newblock Dihedral angles of convex polyhedra.
\newblock {\em Discrete Comput. Geom.}, 23(3):409--417, 2000.

\bibitem{SCH1}
Jean-Marc Schlenker.
\newblock M\'etriques sur les poly\`edres hyperboliques convexes.
\newblock {\em J. Differential Geom.}, 48(2):323--405, 1998.

\bibitem{SCH3}
Jean-Marc Schlenker.
\newblock Hyperbolic manifolds with convex boundary.
\newblock {\em Invent. Math.}, 163(1):109--169, 2006.

\bibitem{T_NOTES}
W.~P. Thurston.
\newblock Geometry and topology of 3-manifolds.
\newblock Princeton University lecture notes.

\bibitem{THURSTON_BOOK}
William~P. Thurston.
\newblock {\em Three-dimensional geometry and topology. {V}ol. 1}, volume~35 of
  {\em Princeton Mathematical Series}.
\newblock Princeton University Press, Princeton, NJ, 1997.
\newblock Edited by Silvio Levy.

\bibitem{VIN}
{\`E}.~B. Vinberg.
\newblock Discrete groups generated by reflections in {L}oba\v cevski\u\i \
  spaces.
\newblock {\em Mat. Sb. (N.S.)}, 72 (114):471--488; correction, ibid. { 73
  (115) (1967), 303}, 1967.

\bibitem{VINREFL}
{\`E}.~B. Vinberg.
\newblock Hyperbolic groups of reflections.
\newblock {\em Russian Math. Surveys}, 40(1):31--75, 1985.

\bibitem{VINVOL}
{\`E}.~B. Vinberg.
\newblock The volume of polyhedra on a sphere and in {L}obachevsky space.
\newblock In {\em Algebra and analysis (Kemerovo, 1988)}, volume 148 of {\em
  Amer. Math. Soc. Transl. Ser. 2}, pages 15--27. Amer. Math. Soc., Providence,
  RI, 1991.

\bibitem{VS}
{\`E}.~B. Vinberg and O.~V. Shvartsman.
\newblock Discrete groups of motions of spaces of constant curvature.
\newblock In {\em Geometry, II}, volume~29 of {\em Encyclopaedia Math. Sci.},
  pages 139--248. Springer, Berlin, 1993.

\end{thebibliography}
\end{document}